\documentclass[12pt]{amsart}
\usepackage{amssymb}
\usepackage{amsmath}
\usepackage{longtable}
\newcommand\g{{\mathfrak g}}
\newcommand{\Mod}{\operatorname{Mod}}
\newcommand{\Rad}{\operatorname{R}}

\newcommand\GL{\operatorname{GL}}
\newcommand\Sk{\mathcal{S}}
\newcommand\Q{\mathbb Q}
\newcommand\m{\mathfrak m}
\newcommand\M{\mathcal{M}}
\newcommand\id{\operatorname{id}}

\newcommand\Verm{\mathcal{M}}
\newcommand\n{\mathfrak n}
\newcommand\z{\mathfrak z}
\newcommand\OCat{\mathcal{O}}
\renewcommand\t{\mathfrak t}
\newcommand\q{\mathfrak q}
\newcommand\h{\mathfrak h}
\newcommand\vf{\mathfrak v}
\newcommand\p{\mathfrak p}

\newcommand\skewperp{\angle}
\newcommand\SO{\operatorname{SO}}
\newcommand\codim{\operatorname{codim}}
\newcommand\Spec{\operatorname{Spec}}

\newcommand\F{\operatorname{F}}
\newcommand\W{{\bf A}}
\newcommand\K{\mathbb K}
\newcommand\U{\mathcal U}

\newcommand\D{\mathcal D}
\newcommand\Ann{\operatorname{Ann}}
\newcommand\Id{\mathfrak{Id}}
\newcommand\Ind{\operatorname{Ind}}
\newcommand\Walg{\mathcal W}
\newcommand\Z{\mathbb Z}
\newcommand\A{\mathcal A}
\newcommand\B{\mathcal B}
\newcommand\N{\mathbb N}
\newcommand\gr{\operatorname{gr}}

\newcommand\Orb{\mathbb{O}}
\newcommand\I{\mathcal I}
\newcommand\J{\mathcal J}
\renewcommand\sl{\mathfrak{sl}}
\newcommand\Sp{\mathop{\rm Sp}\nolimits}
\newcommand\Span{\operatorname{Span}}
\newcommand\Hom{\operatorname{Hom}}
\newcommand{\ad}{\mathop{\rm ad}\nolimits}
\newcommand{\Ad}{\mathop{\rm Ad}\nolimits}
\newcommand{\red}{/\!/\!/}
\newcommand\Centr{\mathcal Z}

\newcommand\mult{\operatorname{mult}}
\newcommand{\VA}{\operatorname{V}}

\newcommand{\Ocat}{\mathcal{O}}
\newcommand\quo{/\!/}
\newcommand\Mat{\operatorname{Mat}}
\newcommand\Rep{\operatorname{Rep}}
\newcommand\Vect{\operatorname{Vect}}
\newcommand\eu{\mathbf{eu}}
\newcommand\SL{\operatorname{SL}}
\newtheorem{Thm}{Theorem}[subsection]
\newtheorem{Prop}[Thm]{Proposition}
\newtheorem{Cor}[Thm]{Corollary}
\newtheorem{Lem}[Thm]{Lemma}
\theoremstyle{definition}

\newtheorem{defi}[Thm]{Definition}
\newtheorem{Rem}[Thm]{Remark}
\newtheorem{Conj}[Thm]{Conjecture}
\unitlength=1mm   \numberwithin{equation}{section}
\numberwithin{table}{section} \oddsidemargin=0cm
\author{Ivan Losev}
\title{1-dimensional representations and parabolic induction for W-algebras}
\thanks{Address: MIT, Department of Mathematics, 77 Massachusetts Avenue, Cambridge, MA 02139, United States}
\thanks{e-mail: ivanlosev@math.mit.edu}
\begin{document}
\begin{abstract}
A W-algebra is an associative algebra constructed from a reductive Lie algebra and
its nilpotent element. This paper concentrates on the study of 1-dimensional representations
of  W-algebras. Under some conditions on a nilpotent element (satisfied by all rigid elements) we obtain a criterium for a finite dimensional module to have dimension 1.
It is stated in terms of the Brundan-Goodwin-Kleshchev highest weight theory.  This criterium allows to compute highest weights for certain completely prime primitive ideals in universal enveloping algebras. We make an
explicit computation in a special case in type $E_8$. Our second principal result is a version
of a parabolic induction for W-algebras. In this case, the parabolic induction is an exact
functor between the categories of finite dimensional modules for two different W-algebras.
The most important feature of the functor is that it preserves dimensions. In particular, it preserves
one-dimensional representations. A closely related result was obtained previously by Premet.
We also establish some other properties of the parabolic induction functor.
\end{abstract}
\maketitle

\section{Introduction}
\subsection{W-algebras and their 1-dimensional representations}
Our base field is an algebraically closed field $\K$ of
characteristic $0$.  A W-algebra $U(\g,e)$ (of finite type) is a
certain finitely generated associative algebra constructed from a
reductive Lie algebra $\g$ and its nilpotent element $e$. There are
several definitions of W-algebras. The definitions introduced  in
\cite{Premet1} and \cite{Wquant} will be given in Subsection \ref{SUBSECTION_Walg_def}.

During the last decade W-algebras were extensively studied starting from Premet's paper
\cite{Premet1}\footnote{We should mention here that there are affine counterparts
of finite W-algebras and they were studied even before finite ones. The general
definition of affine W-algebras was sketched in Kac's ICM talk \cite{Kac_ICM}, and
fully developed in \cite{KRW}.}, see, for instance, \cite{BK1},\cite{BGK},\cite{GG},\cite{Ginzburg},\cite{Wquant}-\cite{LOCat},\cite{Premet2}-\cite{Premet4}.
One of the reasons why W-algebras are interesting is that they
are closely related to the universal enveloping algebra $U(\g)$. For instance, there is a map
$\I\mapsto \I^\dagger$ from the set of two-sided ideals of $U(\g,e)$ to the set of two-sided
ideals  of $U(\g)$ having many nice properties, see \cite{Ginzburg},\cite{Wquant},\cite{HC},\cite{Premet2},\cite{Premet3}, for details.

The central topic in this paper  is the study of 1-dimensional representations of W-algebras.
A key conjecture here was stated by Premet in \cite{Premet2}.

\begin{Conj}\label{Conj:2}
For any $e$ the algebra  $U(\g,e)$ has at least one 1-dimensional representation.
\end{Conj}

There are two major results towards this conjecture obtained previously.

First, the author proved Conjecture \ref{Conj:2} in \cite{Wquant} provided
the algebra $\g$ is classical. Actually, the proof given there can be easily
deduced from earlier results: by Moeglin \cite{Moeglin2}, and
Brylinski \cite{Brylinski}.

Second, in \cite{Premet4} Premet reduced Conjecture \ref{Conj:2} to the case when
the nilpotent orbit is rigid, that is, cannot be induced (in the sense of
Lusztig and Spaltenstein) from a proper Levi subalgebra. Namely, let $\underline{\g}$
be a Levi subalgebra in $\g$. To a nilpotent orbit $\underline{\Orb}\subset\underline{\g}$
Lusztig and Spaltenstein assigned a nilpotent orbit $\Orb=\Ind_{\underline{\g}}^\g(\underline{\Orb})$
in $\g$. See Subsection \ref{SUBSECTION_LSind} for a precise definition.
Premet checked in \cite[Theorem 1.1]{Premet4}, that $U(\g,e)$ has a one-dimensional
representation provided $U(\underline{\g},\underline{e})$ does, where $\underline{e}\in \underline{\Orb}$.
His proof used a connection with modular representations of semisimple Lie algebras.

There are several reasons to be interested in one-dimensional $U(\g,e)$-modules.
They give rise to completely prime primitive ideals in $U(\g)$. This was first
proved by Moeglin in \cite{Moeglin1} (her Whittaker models correspond to
one-dimensional $U(\g,e)$-modules via the Skryabin equivalence from \cite{Skryabin}).
Moreover, as Moeglin proved in \cite{Moeglin2}, any one-dimensional $U(\g,e)$-module
(=Whittaker model) leads to an (appropriately understood)
{\it quantization} of a covering of $Ge$, see Introduction to \cite{Moeglin2}.

On the other hand, in \cite{Premet4} Premet proved that the existence of a one-dimensional
$U(\g,e)$-module implies the Humphreys conjecture on the existence   of  a {\it small}
non-restricted representation in  characteristic $p\gg 0$.\footnote{Actually, this conjecture
that appeared in Humphreys review \cite{Humphreys} was also stated before (in one form or
another). It appeared in the same form in Kac's review, \cite{Kac_review}, to \cite{Premet_KW}
and in \cite{Premet_KW} itself but in a slightly different form.}

\subsection{Main results}
One of the  main results of this paper strengthens Theorem 1.1 from \cite{Premet4}.

\begin{Thm}\label{Thm:1}
Let $\g$ be a semisimple Lie algebra, $\underline{\g}$ its Levi subalgebra, $\underline{e}$ an element of
a nilpotent orbit $\underline{\Orb}$ in
$\underline{\g}$, and $e$ an element from the induced nilpotent orbit
$\Ind_{\underline{\g}}^\g(\underline{\Orb})$.
Then there is a dimension preserving exact   functor $\rho$  from the category of finite dimensional
$U(\underline{\g},\underline{e})$-modules to the category of finite dimensional $U(\g,e)$-modules.
\end{Thm}

The functor $\rho$ mentioned in the theorem will be called the {\it parabolic induction} functor.
The proof of this theorem  given in Subsection \ref{SUBSECTION_proof_main} follows from a stronger result:
we show that $U(\g,e)$ can be embedded
into a certain completion of $U(\underline{\g},\underline{e})$. This completion has the property that any finite dimensional representation of $U(\underline{\g},\underline{e})$   extends to it. Now our functor is just the
pull-back. To get an embedding we use techniques of quantum Hamiltonian reduction.

In Subsection \ref{SUBSECTION_ideals} we will relate our parabolic induction functor to the
parabolic induction map between the sets of ideals of $U(\underline{\g})$ and $U(\g)$.

So to complete the proof of Conjecture \ref{Conj:2} it remains to deal with rigid orbits in exceptional
Lie algebras. The first result here is also due to Premet. In \cite{Premet2} he proved that $U(\g,e)$ has a one-dimensional representation provided $e$ is a minimal nilpotent element (outside of type $A$ such an element is always rigid)\footnote{When the present paper was almost ready to be made public, another result
towards Premet's conjecture appeared: in \cite{GRU} Goodwin, R\"{o}hrle and Ubly checked that $U(\g,e)$
has a one-dimensional representation provided $\g$ is $G_2,F_4,E_6$ or $E_7$. They use Premet's approach
based on the analysis of relations in $U(\g,e)$,
\cite{Premet2}, together with GAP computations.}.

We make a further step to the proof of Conjecture \ref{Conj:2}. Namely, under some condition
on a nilpotent element $e\in\g$ we find a criterium for a finite dimensional module over $U(\g,e)$
to be one-dimensional in terms of the highest weight theory for W-algebras established in \cite{BGK} and further studied in \cite{LOCat}. The precise statement of this criterium, Theorem \ref{Thm:6.5}, will be given  in Subsection
\ref{SUBSECTION_highest_1dim}. The condition on $e$
is that the reductive part of the centralizer $\z_\g(e)$ is  semisimple. In particular, this is
the case for all rigid nilpotent elements. Moreover, our criterium
makes it possible to find explicitly  highest weights of some completely prime primitive
ideals in $U(\g)$, see Subsection \ref{SUBSECTION_toolkit}.

%

\subsection{Content of this paper}
The paper consists of sections that are divided into subsections. Definitions, theorems,
etc., are numbered within each subsection. Equations are numbered within sections.

Let us describe the content of this paper in more detail.
In Section \ref{SECTION_Notation} we gather some  standard notation
used in the paper. In Section \ref{SECTION_Ham} we recall different more or less standard
facts regarding deformation quantization, classical and quantum Hamiltonian actions,
and the Hamiltonian reduction. Essentially, it does not contain new results, with possible
exception of some technicalities. In Section \ref{SECTION_Walg} different facts about
W-algebras are gathered. In Subsection \ref{SUBSECTION_Walg_def} we recall two definitions
of W-algebras (from \cite{Wquant} and \cite{Premet1}). In Subsection \ref{SUBSECTION_decomp}
we recall a crucial technical result about W-algebras, the decomposition theorem from
\cite{Wquant}, as well as two maps between the sets of two-sided ideals
of $U(\g)$ and of $U(\g,e)$ also defined in \cite{Wquant}. Finally in Subsection \ref{SUBSECTION_Ocat}
we recall the notion of the categories $\Ocat$ for W-algebras introduced in \cite{BGK} and
a category equivalence theorem  from \cite{LOCat} (Theorem \ref{Thm:4.4}). This theorem
asserts that there is an equivalence between the category $\Ocat$ for $U(\g,e)$ and a certain category
of generalized Whittaker modules for $U(\g)$.

In Section \ref{SECTION_O_to_onedim} we apply Theorem \ref{Thm:4.4} to the study of irreducible finite dimensional
and, in particular, one-dimensional $U(\g,e)$-modules. Recall that, according to \cite{BGK}, one can consider an irreducible finite dimensional $U(\g,e)$-module $N$ as the irreducible highest weight module $L(N^0)$ (below we use the notation $L^\theta(N^0)$), where the ``highest weight'' $N^0$ is a module over
the W-algebra $U(\g_0,e)$,  $\g_0$ being a Levi subalgebra of $\g$ containing $e$.

In Subsection \ref{SUBSECTION_highest_findim} we prove a criterium, Theorem \ref{Thm:6.1}, for $L(N^0)$ to be finite dimensional. This criterium generalizes   \cite[Conjecture 5.2]{BGK} and is stated in terms of primitive ideals in $U(\g)$. In Subsection \ref{SUBSECTION_highest_1dim}, under the condition on $e$ mentioned
in the previous subsection, we prove a criterium for $L^\theta(N^0)$ to be one-dimensional,
Theorem \ref{Thm:6.5}. Roughly speaking, this criterium asserts  that a primitive ideal $J(\lambda)$ corresponds to a one-dimensional representation of $U(\g,e)$ if $\lambda$ satisfies certain four conditions, the most implicit (as well as the most difficult to check) one being that the associated variety of $J(\lambda)$ is $\overline{\Orb}$.
In Subsection \ref{SUBSECTION_toolkit} we provide some technical statements that allow to check whether the
last condition holds. Also we verify the  four conditions for an explicit highest weight in the case
of the nilpotent element $A_5+A_1$ in $E_8$.

Section \ref{SECTION_induction} is devoted to the study of a parabolic induction for
W-algebras. In Subsection \ref{SUBSECTION_LSind} we recall the definition and some
properties of the Lusztig-Spaltenstein induction for nilpotent orbits. Subsection
\ref{SECTION_classical} is very technical. There we prove some results about
certain classical Hamiltonian actions to be used in the next two  subsections.
The first of them, Subsection \ref{SUBSECTION_proof_main}, contains the proof of
Theorem \ref{Thm:1}. Subsection \ref{SUBSECTION_ideals} relates the parabolic induction
for W-algebras with the parabolic induction for ideals in universal enveloping algebras
(Corollary \ref{Cor:2.8.1}). In Subsection \ref{SUBSECTION_rep_schemes} we study the morphism
of representation schemes induced by the parabolic induction. Our main result is that
this morphism is finite. In Subsection \ref{SUBSECTION_adjoint} we show that the parabolic
induction functor has both left and right adjoint functors.

Finally, in Section \ref{SECTION_appendix} we establish two results, which are not
directly related to the main content of this paper but seem to be of interest. In Subsection \ref{SUBSECTION_Fedosov_vs_TDO}
we will relate the Fedosov quantization of the cotangent bundle to the quantization
 by twisted differential operator. The material of this subsection seems
 to be a folklore knowledge that was never written down explicitly.
 In Subsection \ref{SUBSECTION_Ass} we prove a general
result on filtered algebras that is used in Subsection \ref{SUBSECTION_rep_schemes}.

{\bf Acknowledgements.} This paper would never appear without numerous fruitful discussions
with Alexander Premet and David Vogan. Their precious help is gratefully acknowledged.
I also wish to thank Roman Bezrukavnikov, Jonathan Brundan, Alexander Elashvili, Pavel Etingof, Anthony Joseph, Alexander Kleshchev,  George Lusztig, and Roman Travkin for stimulating discussions. Finally, I wish to thank Simon Goodwin and Alexander Premet for pointing out several mistakes and gaps in
the previous versions of this text.

\section{Notation and conventions}\label{SECTION_Notation}
{\bf Algebraic groups and their Lie algebras}. If an algebraic group is denoted by a capital
Latin letter, e.g., $G$ or $N^-$, then its Lie algebra is denoted by the corresponding
small German letter, e.g., $\g,\n^-$.

{\bf Locally finite parts}. Let $\g$ be some Lie algebra and let $M$
be a module over $\g$. By the locally finite (shortly, l.f.) part of
$M$ we mean the sum of all finite dimensional $\g$-submodules of
$M$. The similar definition can be given for algebraic group
actions.

{\bf The main  algebras.} Below $G$ is a connected reductive algebraic group and
$\g$ is its Lie algebra. The universal enveloping algebra $U(\g)$ will be mostly denoted by
$\U$. Let $e$ be a nilpotent element in $\g$. The W-algebra $U(\g,e)$ will be denoted shortly by $\Walg$. In Section \ref{SECTION_O_to_onedim}
we will consider the situation when $e$ lies in a certain Levi subalgebra $\g_0\subset \g$. In this case
we set $\U^0:=U(\g_0), \Walg^0:=U(\g_0,e)$. In Section \ref{SECTION_induction} the element $e$
will be induced from a nilpotent element $\underline{e}$ of some Levi subalgebra
$\underline{\g}$ of $\g$. There we will write $\underline{\U}:=U(\underline{\g}), \underline{\Walg}:=
U(\underline{\g},\underline{e})$.

\begin{longtable}{p{3cm} p{12cm}}
$\W_V$& the Weyl algebra of a symplectic vector space $V$.\\
$\Ann_\A(M)$& the annihilator of an $\A$-module $M$ in an  algebra
$\A$.\\
$\gr \A$& the associated graded algebra of a filtered
algebra $\A$.\\
$G_x$& the stabilizer of a point $x$ under an action of a group $G$.\\
$(G,G)$& the derived subgroup of a group $G$.\\
$H^i_{DR}(X)$& $i$-th De Rham cohomology of a smooth algebraic variety (or of a formal
scheme) $X$.
\\
$\Id(\A)$& the set of all (two-sided) ideals of an algebra $\A$.\\
$\sqrt{\J}$& the radical of an ideal $\J$.\\
$\K[X]^\wedge_{Y}$& the completion of the algebra $\K[X]$ of regular
functions on $X$ with respect to a subvariety $Y\subset X$.\\
$\mult_{\Orb}(\M)$& the multiplicity of a Harish-Chandra $U(\g)$-bimodule $\M$
on an open orbit $\Orb\subset \VA(\M)$.\\
$N\leftthreetimes L$& the semidirect product of groups $N$ and $L$ ($N$ is normal).\\
$R_\hbar(\A)$& the Rees algebra of a filtered algebra $\A$.\\
$\Rad_u(H)$& the unipotent radical of an algebraic group $H$.\\
$\Span_A(X)$& the $A$-linear span of a subset $X$ in some $A$-module.
\\ $U^\skewperp$& the skew-orthogonal complement to a subspace $U$ in a symplectic
vector space.
\\$\VA(\M)$& the associated variety (in $\g^*$) of a Harish-Chandra $U(\g)$-bimodule $\M$.
\\ $X^G$& the fixed-point set for an action $G:X$.
\\
$\z_\g(x)$& the centralizer of $x$ in $\g$.\\
$Z_G(x)$& the centralizer of $x\in\g$ in $G$.
\\
$\Omega^i(X)$& the space of $i$-forms on a variety $X$.\\
$\widehat{\otimes}$& the completed tensor product (of complete topological vector spaces or
modules).
\end{longtable}

\section{Deformation quantization, Hamiltonian actions and Hamiltonian reduction}\label{SECTION_Ham}
\subsection{Deformation quantization}\label{SUBSECTION_defquant}
Let $A$ be a commutative associative algebra with unit equipped with
a Poisson bracket.

\begin{defi}\label{defi:1.1}
A $\K[[\hbar]]$-bilinear map $*:A[[\hbar]]\times A[[\hbar]]\rightarrow
A[[\hbar]]$ is called a {\it star-product} if it satisfies the
following conditions: \begin{itemize}
\item[(*1)] $*$ is associative, equivalently,  $(f*g)*h=f*(g*h)$ for all $f,g,h\in A$,
and $1\in A$ is a unit for $*$.
\item[(*2)] $f*g-fg\in \hbar^2 A[[\hbar]],f*g-g*f-\hbar^2\{f,g\}\in \hbar^4 A[[\hbar]]$ for all $f,g\in A$.
\item[(*3)] $A*A\subset A[[\hbar^2]]$
\end{itemize}
\end{defi}
Clearly, a star-product is uniquely determined by its restriction to
$A$. One may write $f*g=\sum_{i=0}^\infty D_i(f,g)\hbar^{2i},f,g\in A,
D_i:A\otimes A\rightarrow A$.  If all $D_i$ are
bidifferential operators, then the star-product $*$ is called {\it
differential}. In this case we can extend the star-product to $B[[\hbar]]$ for any
localization or completion $B$ of $A$. When we consider $A[[\hbar]]$ as an algebra with respect to
the star-product, we call it a {\it quantum algebra}.

Also we remark that usually in condition (*2) one has $f*g-fg\in \hbar A[[\hbar]],f*g-g*f-\hbar\{f,g\}\in \hbar^2 A[[\hbar]]$ (so we get our definition from the usual one replacing $\hbar$ with $\hbar^2$).

Let  $G$ be an algebraic group acting on $A$ by automorphisms. It
makes sense to speak about  $G$-invariant star-products ($\hbar$ is
supposed to be $G$-invariant). Now let $\K^\times$ act on $A,
(t,a)\mapsto t.a,$ by automorphisms. Consider the action
of $\K^\times$ on $A[[\hbar]]$ given by $$t.\sum_{i=0}^\infty
a_j\hbar^j=\sum_{j=0}^\infty t^{j}(t.a_j)\hbar^j.$$ If $\K^\times$
acts by automorphisms of $*$, then we say that $*$ is {\it homogeneous}.
Clearly, $*$ is homogeneous if and only if the map $D_l:A\otimes A\rightarrow
A$ is homogeneous of degree $-2l$.

Let $X$ be a smooth affine variety equipped with a symplectic form
$\omega$. Let $A:=\K[X]$ be its  algebra of regular functions.   There is a construction of a differential star-product due to Fedosov. According to this construction, see \cite[Section 5.3]{Fedosov2} one needs to fix a symplectic connection (=a torsion-free affine connection annihilating the symplectic form),
say $\nabla$,  and a $K[[\hbar^2]]$-valued 2-form $\Omega$ on $X$.
Then from $\omega,\nabla$
and $\Omega$ one canonically constructs differential operators $D_i$ defining a
star-product.

Now suppose $G\times\K^\times$ acts on $X$. If $\omega,\nabla,\Omega$ are $G$-invariant, then $*$
is $G$-invariant. If $\nabla,\Omega$ are $\K^\times$-invariant and $t.\omega=t^2\omega$ for any $t\in \K^\times$ (recall, that $t.\hbar=t\hbar$), then $*$ is homogeneous. Note also that
if $G$ is reductive, then there is always a $G\times\K^\times$-invariant symplectic connection on
$X$, \cite[Proposition 2.2.2]{Wquant}.

Let us now consider the question of the classification of  star-products on $\K[X][[\hbar]]$. We say that two
$G$-invariant homogeneous star-products $*,*'$ are {\it equivalent} if there is a $G\times\K^\times$-equivariant
map (equivalence) $T:=\id+\sum_{i=1}^\infty T_i\hbar^{2i}:\K[X][[\hbar]]\rightarrow \K[X][[\hbar]]$, where
$T_i:\K[X]\rightarrow \K[X]$ is a linear map, such that $(Tf)*'(Tg)=T(f*g)$. We say that the equivalence $T$ is {\it differential}
when all $T_i$ are differential operators.

\begin{Thm}\label{Thm:1.3}
\begin{enumerate}
\item Let $*,*'$ be Fedosov star-products constructed from the pairs $(\nabla,\Omega),$ $(\nabla',\Omega')$ such that
$\nabla,\nabla',\Omega,\Omega'$ are $G\times\K^\times$-invariant. Then $*,*'$ are equivalent if and only if
$\Omega'-\Omega$ is exact. In this case one can choose a differential equivalence.
\item Any $G$-invariant homogeneous star-product is $G\times\K^\times$-equivariantly equivalent to a Fedosov one.
\end{enumerate}
\end{Thm}

The first part of the theorem is, essentially, due to Fedosov, see \cite[Section 5.5]{Fedosov2}.
The second assertion is an easy special case of  \cite[Theorem 1.8]{Kaledin_Bezrukavnikov}.
The fact that an equivalence can be chosen $G\times \K^\times$-equivariant follows from the
proof of that theorem, compare with   \cite[Subsection 6.1]{Kaledin_Bezrukavnikov}.

\begin{Rem}
Although Fedosov worked in the $C^\infty$-category, his results remain valid for smooth affine
varieties too. Also  proofs of parts of Theorem \ref{Thm:1.3} in \cite{Fedosov2} and \cite{Kaledin_Bezrukavnikov}  of the results that we need can be  generalized to the case of smooth formal schemes  in a straightforward way.
\end{Rem}


\subsection{Classical Hamiltonian actions}\label{SUBSECTION_Ham_class}
In this subsection $G$ is an algebraic group and $X$ is a smooth affine
variety equipped with a regular symplectic form $\omega$ and an
action of $G$ by symplectomorphisms. Let $\{\cdot,\cdot\}$ denote
the Poisson bracket on $X$ induced by $\omega$.

To any element $\xi\in\g$ one assigns in a standard way the velocity vector field
$\xi_X$ on $X$. Suppose  there is a linear map $\g\rightarrow \K[X],
\xi\mapsto H_\xi,$ satisfying the following two conditions:
\begin{itemize}
\item[(H1)] The map $\xi\mapsto H_\xi$ is $G$-equivariant.
\item[(H2)] $\{H_\xi,f\}=\xi_Xf$ for any  $f\in \K[X]$.
\end{itemize}

\begin{defi}\label{Def:1.1.1}
The action  $G:X$ equipped with a linear map $\xi\mapsto H_\xi$
satisfying (H1),(H2) is said to be {\it Hamiltonian} and $X$ is
called a Hamiltonian $G$-variety. The functions $H_\xi$ are said to be the hamiltonians
of the action.
\end{defi}

For a Hamiltonian action $G:X$ we define a morphism
$\mu:X\rightarrow \g^*$ (called a {\it moment map}) by the formula
\begin{equation*}
\langle \mu(x),\xi\rangle= H_{\xi}(x),\xi\in\g,x\in X.
\end{equation*}
The map $\xi\mapsto H_\xi$ is often referred to as a {\it comoment map}.

Let $X_1,X_2$ be two $G$-varieties with Hamiltonian $G$-actions
and $\varphi$ be a $G$-equivariant morphism $X_1\rightarrow X_2$. We say that $\varphi$
is {\it Hamiltonian} if $\varphi$   intertwines the symplectic forms and the moment
maps.

Now let us recall a local description of Hamiltonian actions obtained in \cite{slice}.

Till the end of the subsection we assume that $G$ is reductive.
We identify
$\g$ with $\g^*$ using an invariant symmetric form $(\cdot,\cdot)$, whose restriction to the rational part of a Cartan subalgebra in $\g$ is positively definite. Then the restriction of $(\cdot,\cdot)$ to the Lie algebra of any reductive subgroup of $G$ is non-degenerate. So for any such subalgebra $\h\subset\g$ (in particular, for
$\h=\g$) we identify $\h$ with $\h^*$ using the form. For a point $x\in X$ let $\g_*x$ denote the tangent space
to the orbit $Gx$ in $x$.

Let $x\in X$ be a point with  closed $G$-orbit. Set $H=G_x, \eta=\mu_G(x)$ and $V:=\g_*x^\skewperp/(\g_*x\cap \g_*x^\skewperp)$ (in other words, $V$ is the symplectic part of the normal space to $Gx$ in $X$).
Then $V$ is a symplectic $H$-module.


\begin{Prop}[\cite{slice}]\label{Prop:2.41}
Let $X_1,X_2$ be two Hamiltonian $G$-varieties and $x_1\in X_1,x_2\in X_2$ be two points with
closed $G$-orbits. Suppose that $G_{x_1}=G_{x_2}, \mu_G(x_1)=\mu_G(x_2)$ and the $G_{x_i}$-modules
$V_i:=\g_*x_i^\skewperp/(\g_*x_i\cap \g_*x_i^\skewperp)$ are isomorphic. Then there is a $G$-equivariant Hamiltonian morphism $\varphi:(X_1)^\wedge_{Gx_1}\xrightarrow{\sim}(X_2)^\wedge_{Gx_2}$ mapping $x_1$ to $x_2$.

Suppose, in addition, that
we have $\K^\times$-actions on $X_1,X_2$ such that
\begin{itemize}
\item[(A)] they stabilize $Gx_i$ and the stabilizers of $x_1,x_2$ in $G\times \K^\times$ are the same.
\item[(B)] $t\in \K^\times$ multiplies
the symplectic forms and the functions $H_\xi$ by $t^2$.
\item[(C)] $V_1$ and $V_2$ are isomorphic as $(G\times \K^\times)_{x_i}$-modules.
\end{itemize}
Then an isomorphism $\varphi$ can be made,
in addition, $\K^\times$-equivariant.
\end{Prop}
\begin{proof}
We remark that if two symplectic modules are equivariantly isomorphic, then they are actually equivariantly symplectomorphic.
Now the part without a $\K^\times$-action basically follows from the uniqueness result in \cite{slice} (which was stated
in the complex-analytic category). See also the preprint \cite[Theorem 5.1]{Knop} for the proof
in the formal scheme setting.

The proof in \cite{slice} can  be adjusted to the $\K^\times$-equivariant situation too,
compare with the proof of   \cite[Theorem 3.1.3]{Wquant}.
\end{proof}

Let us construct explicitly a Hamiltonian $G$-variety corresponding to a triple $(H,\eta,V)$, a so called
{\it model variety $M_G(H,\eta,V)$} introduced in \cite{slice}. For simplicity we will assume that $\eta=e$
is nilpotent. In this case there is also a $\K^\times$-action satisfying the conditions (A),(B),(C)
of the previous proposition.

Let $V$ be a symplectic $H$-module with symplectic form $\omega_V$.
Define a  subspace $U\subset \g$ as follows. If $e=0$, then
$U:=\h^\perp$. Otherwise, consider an $\sl_2$-triple $(e,h,f)$ in $\g^H$
and set $U:=\z_\g(f)\cap \h^\perp$.

Consider the homogeneous vector bundle $X:=G*_H(U\oplus V)$ and embed $U\oplus V$ to $X$ by
$(u,v)\mapsto [1,(u,v)]$. Define the map $\mu_G:U\oplus V\rightarrow \g$ by
$$\mu_G([1,(u,v)])=e+u+\mu_{H}(v),$$
where $\mu_H:V\rightarrow \h$ is the moment map for the action of $H$ on $V$ given by $(\mu_H(v),\xi)=\frac{1}{2}\omega(\xi v,v)$. Since this map is $H$-equivariant, we can extend it
uniquely to a $G$-equivariant map $\mu_G:X\rightarrow \g$.

Also define  the 2-form $\omega_x\in \bigwedge^2 T^*_x X, x=[1,(u,v)],$ by
\begin{align*}&\omega_{x}(\xi_1+u_1+v_1,\xi_2+u_2+v_2)=\langle\mu_G(x),[\xi_1,\xi_2]\rangle+\langle\xi_1,u_2\rangle-
\langle\xi_2,u_1\rangle+\omega_V(v_1,v_2),\\
& \xi_1,\xi_2\in \h^\perp, u_1,u_2\in U, v_1,v_2\in V.\end{align*}

Note that the section  $\omega: x\mapsto \omega_x$ of $\bigwedge^2 T^*X|_{U\oplus V}$ is
$H$-invariant, so we can extend it to $X$ by $G$-invariance.
It turns out that the form $\omega$ is symplectic and $\mu_G:X\rightarrow \g$ is a moment
map for this form, see \cite{slice}.

Define the $\K^\times$-action on $X$ as follows. Let $\gamma:\K^\times\rightarrow G$ be
the composition of the homomorphism $\SL_2\rightarrow G$ induced by the $\sl_2$-triple
and of the embedding $\K^\times\rightarrow \SL_2$ given by $t\mapsto \operatorname{diag}(t,t^{-1})$.
Define a $\K^\times$-action
on $M_G(H,\eta,V)$ as follows:
$$t.(g,u,v)=(g\gamma(t)^{-1}, t^{-2}\gamma(t)u,t^{-1}\beta(t)v),$$
where $\beta$ is a group homomorphism $\K^\times\rightarrow \Sp(V)^H$.
The action of $t\in \K^\times$ multiplies $\omega$ and $\mu_G$ by $t^2$.

In the sequel we will need a technical result concerning Hamiltonian actions on formal schemes.

\begin{Lem}\label{Lem:1.1.2}
Let $G$ be a reductive group, and $X$ be an affine Hamiltonian $G$-variety with symplectic form $\omega$
and  moment map $\mu_G$. Let $x\in X$ be a point with closed $G$-orbit. Suppose that
$\mu_G(x)$ is nilpotent. Further, let  $\z$ be an algebraic Lie subalgebra in $\z(\g)$ satisfying $\z\oplus (\g_x+[\g,\g])=\g$. Finally, let $\zeta$ be a $G$-invariant symplectic vector field on $X^\wedge_{Gx}$ such that $\zeta H_\xi=0$ for all $\xi\in\z$. Then  $\zeta=v_f$ for a $G$-invariant element $f\in \K[X]^\wedge_{Gx}$,
where $v_f$ is the Hamiltonian vector field associated with $f$.
\end{Lem}
Actually, the proof generalizes directly to the case when $\eta$ is not necessarily nilpotent,
but we will not need this more general  result.
\begin{proof}
Note that the $v_f H_\xi=0$ for any $f\in (\K[X]^\wedge_{Gx})^G$.

Replacing $G$ with a covering we may assume that $G=Z\times G_0$, where $Z\subset Z(G)$ is the torus with Lie algebra
 $\z$ and $G_0$ is the product of
a torus and of a simply connected semisimple group with $G_0= H^\circ (G,G)$.
Thanks to Proposition \ref{Prop:2.41}, we can replace $X$ with a model variety
$M_G(H,\eta,V)$.

Let us consider the model variety $\widetilde{X}:= M_G(H^\circ,\eta,V)$.
From the construction of model varieties recalled above, we have the action of the finite group $\Gamma:=H/H^\circ$
on $\widetilde{X}$ by Hamiltonian automorphisms, and $X=\widetilde{X}/\Gamma$.
We have the decomposition $\widetilde{X}=T^*Z\times \widetilde{X}_0$, where
$\widetilde{X}_0:=M_{G_0}(H^\circ,\eta,V)$. The homogeneous space $G_0/H^\circ$
is simply connected. Note  also that the action of $\Gamma$ on $\widetilde{X}$ preserves the decomposition
$\widetilde{X}=T^*Z\times \widetilde{X}_0$. Pick a point $\widetilde{x}\in \widetilde{X}$ mapping
to $x=[1,(0,0)]\in X=M_G(H,\eta,V)$. Lift $\omega$ and $\zeta$ to $\widetilde{X}^\wedge_{G\widetilde{x}}$. We will denote the
liftings by the same letters.

Let us check that the contraction $\iota_\zeta\omega$ is an exact form.

The projection $\widetilde{X}^\wedge_{G\widetilde{x}}\twoheadrightarrow Z$ induces an isomorphism $H^1_{DR}(\widetilde{X}^\wedge_{Gx})
\rightarrow H^1_{DR}(Z)$.
Fix an identification
$Z\cong (\K^\times)^m$ and let $z_i, i=1,\ldots,m,$ be the coordinate on the $i$-th copy of $\K^\times$.
The forms $\frac{dz_i}{z_i}$ form a basis in $H^1_{DR}(Z)$. Choose the basis $p_1,\ldots,p_m$ in $\z^*$
corresponding to $z_1,\ldots,z_m$. Then the symplectic form on $T^*Z$ is written as
$\sum_{i=1}^m  dp_i\wedge dz_i$. Consider the vector field $\zeta_i:=\frac{1}{z_i}\frac{\partial}{\partial p_i}$,
then $\iota_{\zeta_i}\omega=\frac{dz_i}{z_i}$. The classes of the forms $\iota_{\zeta_i}\omega$ form a
basis in $H^1_{DR}(Z)$.
Finally, let $\xi_1,\ldots,\xi_m$ be the basis
of $\z$ corresponding to the choice of $z_1,\ldots,z_m$. Then $H_{\xi_i}=p_iz_i$ and
$\zeta_j H_{\xi_i}=\delta_{ij}$.

Now let $\zeta$ be an arbitrary $G$-invariant symplectic vector field on $\widetilde{X}^\wedge_{Gx}$. There are scalars
$a_1,\ldots,a_m$ such that $\iota_{\zeta}\omega-\sum_{i=1}^m a_i\iota_{\zeta_i}\omega=df$ for some
$f\in \K[\widetilde{X}]^\wedge_{Gx}$. Since the $G\times \Gamma$-module $\K[\widetilde{X}]^\wedge_{Gx}$ is pro-finite, we may assume that $f$ is $G\times \Gamma$-invariant. It follows that  $0=\zeta H_{\xi_i}=a_i$. Thus $\zeta=v_f$. Since $f\in (\K[\widetilde{X}]^\wedge_{Gx})^\Gamma= \K[X]^\wedge_{Gx}$, we are done.
\end{proof}

\subsection{Quantum Hamiltonian actions}\label{SUBSECTION_quantum_moment}

In this subsection we consider a quantum version of Hamiltonian actions. We preserve the notation of the previous subsection.

 Let $*$ be a star-product
on $\K[X][[\hbar]]$.   A {\it quantum comoment map} for the action $G:X$ is, by definition, a $G$-equivariant
linear map $\g\rightarrow \K[X][[\hbar]],\xi\mapsto \widehat{H}_\xi,$ satisfying the equality
\begin{equation}\label{eq:1.5}[\widehat{H}_\xi,f]=\hbar^2 \xi_X f, \forall \xi\in\g, f\in \K[X][[\hbar]],\end{equation}
where $[\widehat{H}_\xi,f]:=\widehat{H}_\xi*f-f*\widehat{H}_\xi$.
The elements $\widehat{H}_\xi$ are said to be the quantum hamiltonians of the action. A $G$-equivariant
homomorphism of quantum algebras is called {\it Hamiltonian} if it intertwines the quantum
comoment maps.

The following theorem gives a criterium for the existence of a quantum comoment  map in the case
when $G$ is reductive.

\begin{Thm}[\cite{GR_moment}, Theorem 6.2]\label{Thm:1.4}
Let $X$ be an affine symplectic Hamiltonian $G$-variety, $\xi\mapsto H_\xi$
being the comoment map, and $*$ be the
star-product  on $\K[X][[\hbar]]$  obtained by the Fedosov
construction with a $G$-invariant connection $\nabla$ and
$\Omega\in \Omega^2(X)[[\hbar^2]]^G$. Then the following conditions are
equivalent:
\begin{enumerate} \item The $G$-variety $X$ has a quantum comoment map $\xi\mapsto \widehat{H}_\xi$
with $\widehat{H}_\xi\equiv H_\xi\operatorname{mod}\hbar^2$.
\item The 1-form $i_{\xi_X}\Omega$ is exact for each $\xi$, where $i_{\xi_X}$ stands for the
contraction with $\xi_X$.
\end{enumerate}
Moreover, if (2) holds then for $\widehat{H}_\xi$ we can take $H_\xi+ \hbar^2 a_\xi$, where $\xi\mapsto a_{\xi}$
is a $G$-equivariant map $\g\rightarrow \K[X][[\hbar]]$ such that $d a_{\xi}=\iota_{\xi_X}\Omega$
for any $\xi\in\g$.
\end{Thm}

\begin{Rem}\label{Rem:1.5}
Although Gutt and Rawnsley worked with $C^\infty$-manifolds, their techniques can be carried
over to the algebraic setting without any noticeable modifications. Also Theorem \ref{Thm:1.4} can be
directly generalized to the case when $X$ is a smooth affine formal scheme (say, the completion
of a closed $G$-orbit in an affine variety).
\end{Rem}

\begin{Rem}\label{Rem:1.6}
We preserve the notation of Theorem \ref{Thm:1.4}. Suppose $\K^\times$ acts on $X$ such that
$\nabla,\Omega$ are $\K^\times$-invariant and $\omega, H_\xi$
 have degree 2 for all $\xi\in\g$. Then the last assertion of Theorem \ref{Thm:1.4}
insures that there are $\widehat{H}_\xi$  of degree 2.
\end{Rem}

In particular, Theorem \ref{Thm:1.4} implies that if $\Omega=0$, then we can take $0$ for $a_\xi$ and so  $H_\xi$
for $\widehat{H}_\xi$.

We will need a quantum version of Proposition \ref{Prop:2.41}.

\begin{Thm}\label{Thm:1.11}
Let $X_1,X_2,x_1,x_2$ be such as in Proposition \ref{Prop:2.41}. Suppose
there are actions of $\K^\times$ as in that proposition. Equip $X_1,X_2$ with homogeneous $G$-invariant
star-products corresponding
to $\Omega=0$ and let the quantum hamiltonians coincide with the classical ones.
Then there is a $G\times\K^\times$-equivariant differential Hamiltonian isomorphism $\K[X_1]^\wedge_{Gx_1}[[\hbar]]\rightarrow \K[X_2]^\wedge_{Gx_2}[[\hbar]]$   lifting
$(\varphi^{*})^{-1}$, where $\varphi$ is as in Proposition \ref{Prop:2.41}.
\end{Thm}
\begin{proof}
 Let $X,x,Z,G_0, \widetilde{X},\widetilde{x},\widetilde{X}_0,\Gamma$
 be such as in the proof of Lemma \ref{Lem:1.1.2}.
Let us equip $\K[\widetilde{X}][[\hbar]]$ with a star-product as follows.
Consider the Fedosov star-product on $T^*Z$ constructed from the {\it trivial} connection
and $\Omega=0$. This star-product is invariant w.r.t the action of $\Gamma$ on $T^*Z$.
Then consider some $G\times\K^\times\times \Gamma$-invariant symplectic
connection on $\widetilde{X}_0$ and construct the star-product from this connection and $\Omega=0$.
Taking the tensor product of the two star-products we get a star-product $*$ on
$\widetilde{X}$ (that can be extended to $\widetilde{X}^\wedge_{G\widetilde{x}}$), again, with
the quantum comoment map $\xi\mapsto H_\xi$. This star-product is $\Gamma$-invariant so it
descends to $\K[X][[\hbar]]$.

\begin{Lem}\label{Lem:1.4.1}
Let $\psi\in\z^*$. Then there is a derivation $\widetilde{\zeta}$ of the quantum algebra $\K[X][[\hbar]]$
such that $\widetilde{\zeta}(\hbar)=0,\widetilde{\zeta}(H_\xi)=\langle\psi,\xi\rangle$.
\end{Lem}
\begin{proof}[Proof of Lemma \ref{Lem:1.4.1}]
Again, as in the proof of Lemma \ref{Lem:1.1.2} choose coordinates $z_1,\ldots,z_m$ on $Z$
and let $p_i,\xi_i$ have the same meaning as in that proof. The algebra $\K[T^*Z][[\hbar]]$
is generated (as a $\K[[\hbar]]$-algebra) by the elements $z_i,z_i^{-1},p_i, i=1,\ldots,m,$
subject to the relations $[z_i,p_j]=\hbar^2\delta_{ij},[z_i,z_j]=0,[p_i,p_j]=0$.

Set $a_i=\langle\psi,\xi_i\rangle$.
It is clear that the map $\widetilde{\zeta}:z_i\mapsto 0, p_i\mapsto a_i/z_i$ can be uniquely extended to a
$\K[\hbar]$-linear derivation of the quantum algebra $\K[T^*Z][\hbar]$ also denoted by $\widetilde{\zeta}$.  Extend $\widetilde{\zeta}$ to the derivation of $\K[\widetilde{X}][\hbar]$ by making it act trivially on $\K[\widetilde{X}_0][\hbar]$. Note that $\widetilde{\zeta}(H_\xi)=\langle\psi,\xi\rangle$ for $\xi\in \z$.
By the construction, $\widetilde{\zeta}$ is $\Gamma$-invariant. So it descends to $\K[X][[\hbar]]$.
\end{proof}

By Proposition \ref{Prop:2.41}, we have $G\times\K^\times$-equivariant Hamiltonian isomorphisms $\K[X_1]^\wedge_{Gx_1}\cong \K[X]^\wedge_{Gx}\cong \K[X_2]^\wedge_{Gx_2}$.
Transfer the star-products and the quantum (=classical)
comoment maps from $(X_1)^\wedge_{Gx_1},(X_2)^\wedge_{Gx_2}$ to
$X^\wedge_{Gx}$.
So we get two star-products $*^1,*^2$ on $\K[X]^\wedge_{Gx}[[\hbar]]$
such that $\xi\mapsto H_\xi$ is a quantum comoment map for both of them.

It remains to check that there is a $G\times\K^\times$-equivariant differential equivalence
$T:=\id+\sum_{i=1}^\infty T_i\hbar^{2i}:\K[X]^\wedge_{Gx}\rightarrow \K[X]^\wedge_{Gx}$
such that $T(f*g)=(Tf)*^{1}(Tg)$ and $T(H_\xi)=H_{\xi}$. In general, $T(H_\xi)-H_{\xi}\in \K\hbar^2$
(compare with the proof of   \cite[Theorem 2.3.1]{HC}).

As we have seen in the proof of   \cite[Theorem 2.3.1]{HC}, $T_1$ is a Poisson derivation
of $\K[X]^\wedge_{Gx}$. Let $v$ denote the corresponding vector field.
Let $\zeta_i$ have the same meaning as in Lemma \ref{Lem:1.1.2}.
Then, as we have seen in the proof of that lemma, $v=v_f+\sum_{i=1}^m a_i\zeta_i$
for uniquely determined $a_i\in \K$
and some $f\in \K[X]^\wedge_{Gx}[[\hbar]]^{G}$.

Let $\widetilde{\zeta}$ be a derivation from Lemma \ref{Lem:1.4.1} constructed for the linear function
$\psi\in \z^*$ given by $\langle\psi,\xi_i\rangle=a_i$.
Set $T':=T\circ\exp(-\widetilde{\zeta})$. This is an isomorphism intertwining
the star-products $*,*^1$. Moreover, $T'=\id+\sum_{i=1}^\infty T_i'\hbar^{2i}$, where $T_1'$
now corresponds to a Hamiltonian
vector field.  So $T'(H_\xi)-H_{\xi}$ lies in $\hbar^4\K[[\hbar]]$ whence is zero.
\end{proof}

We will also need the following corollary of Theorems \ref{Thm:1.4},\ref{Thm:1.11}.
\begin{Cor}\label{Cor:1.6}
Let $X$ be an affine Hamiltonian $G$-variety with symplectic form $\omega$
and comoment map $\xi\mapsto H_\xi$.
Let $x\in X$ be a point with closed $G$-orbit and $G_x=\{1\}$. Next, suppose that $X$ is equipped with a $\K^\times$-action and also with an action of a reductive group $Q$ such that
\begin{itemize}
\item $Q$, $\K^\times$ and $G$ pairwise commute.
\item $Q$ acts by Hamiltonian automorphisms.
\item $t.\omega=t^2\omega, t.H_\xi=t^2 H_\xi$ for any $t\in \K^\times$.
\item $Q\times \K^\times$ preserves $Gx$.
\end{itemize}
Let $*$ be a star-product on $X^\wedge_{Gx}$ obtained by the Fedosov construction with
a $G\times Q\times\K^\times$-invariant connection and  $\Omega=0$.
Further, let $*'$ be some other $G\times Q$-equivariant homogeneous star-product
on $X^\wedge_{Gx}$ and $\g\rightarrow \K[X]^\wedge_{Gx}[[\hbar]], \xi\mapsto \widehat{H}'_\xi,$ be a quantum comoment  map for $*'$. We suppose that $\widehat{H}'_\xi\equiv H_\xi \mod \hbar^2$ and that $\widehat{H}'_\xi$ is $Q$-invariant
and of degree 2 with respect to $\K^\times$ for any $\xi\in\g$.
Then there exists a $G\times Q\times \K^\times$-equivariant linear map
$T:\K[X]^\wedge_{Gx}[[\hbar]]\rightarrow \K[X]^\wedge_{Gx}[[\hbar]]$ intertwining the star-products and
the quantum comoment maps.
\end{Cor}
\begin{proof}
We may assume that $*'$ is obtained by the Fedosov construction using a $Q\times\K^\times$-invariant two-form $\Omega'$. Let us show that $\Omega'$ is exact.

From Theorem \ref{Thm:1.4} (applied to formal schemes rather than to  varieties)
it follows that $\iota_{\xi_*}\Omega'$ is exact for any $\xi\in\z(\g)$. The inclusion $G=Gx\hookrightarrow
X^\wedge_{Gx}$ and the projection $G\twoheadrightarrow G/(G,G)$ induce isomorphisms of the second
cohomology groups. So we need to check that if $\theta$ is a closed (left) invariant 2-form on a torus $G/(G,G)$, then $\theta$ is exact whenever the forms $i_{\xi_*}\theta$ are exact for all $\xi\in\z(\g)$. Let $z_i, i=1,\ldots,m,$ have the same
meaning as before.
Then any second cohomology class of $(\K^\times)^m$ is uniquely represented in the form
$\sum a_{ij}\frac{dz_i\wedge dz_j}{z_iz_j}, a_{ij}=-a_{ji}$. The condition of the exactness of all
forms $\iota_{\xi_*}\theta$ means that the forms $\sum_{i}a_{ij}\frac{dz_i}{z_i}$ are exact
for all $j$. So $a_{ij}=0$ and $\Omega'$ is exact.

Applying Theorem \ref{Thm:1.11} (as the proof shows an equivalence there
can be chosen to be, in addition, $Q$-equivariant), we complete the proof of the corollary.
\end{proof}

\subsection{Hamiltonian reduction}\label{SUBSECTION_reduction}
In this subsection $\widetilde{G}$ is an algebraic group and $X$  an affine
Hamiltonian  $\widetilde{G}$-variety equipped with a symplectic form
$\omega$ and with a moment map $\mu_{\widetilde{G}}$.
We fix a normal subgroup $G\subset \widetilde{G}$
such that the $G$-action on $X$ is free. We suppose that there exists a
categorical quotient $X/G$ and that the quotient morphism $X\rightarrow X/G$
is affine.


Let $\mu_G$ be the moment map for the action $G:X$, in other words, $\mu_G$ is the composition
of $\mu_{\widetilde{G}}$ and the natural projection $\widetilde{\g}^*\twoheadrightarrow \g^*$.
By the standard properties of the moment map, compare, for example, with \cite[Section 26]{GS} $\mu_{G}^{-1}(0)$ is a smooth complete intersection in $X$. Further, $\mu_{G}^{-1}(0)/G\subset X/G$ is a smooth subvariety and there is
a unique symplectic form $\underline{\omega}$ on $\mu_{G}^{-1}(0)/G$ whose pullback to $\mu_{G}^{-1}(0)$
coincides with the restriction of $\omega$. So $\underline{\omega}$ is $\widetilde{G}/G$-invariant and the action
of $\widetilde{G}/G$ on   $\mu_{G}^{-1}(0)/G$ is Hamiltonian, the moment map sends an orbit $Gx, x\in \mu_{G}^{-1}(0),$ to $\mu_{\widetilde{G}}(x)$. The variety $\mu_{G}^{-1}(0)/G$ is said
to be the {\it Hamiltonian reduction} of $X$ under the action of $G$ and is denoted by $X\red G$. Note that the
algebra of functions on $X\red G$ is nothing else but $(\K[X]/I)^G$, where $I$
is the ideal in $\K[X]$ generated by $H_\xi,\xi\in\g$.

Suppose there is a $\K^\times$-action
on $X$ such that $t.\omega=t^2\omega$ and $t.\mu_{\widetilde{G}}=t^2\mu_{\widetilde{G}}$. Then this action
descends to $X\red G$ and has analogous properties.

Now let us consider a quantum version of this construction. Suppose that $X$ is equipped with
a $\widetilde{G}$-invariant homogeneous star-product $*$ and $\xi\mapsto \widehat{H}_\xi$ is a quantum
comoment map for this action. Set $\A_\hbar:=(\K[X][[\hbar]]/\I_\hbar)^G$, where $\I_\hbar$ is the left
ideal in $\K[X][[\hbar]]$ generated by $\widehat{H}_\xi,\xi\in\g$.

\begin{Prop}
$\A_\hbar$ is a complete flat $\K[[\hbar]]$-algebra, $\I_\hbar/(\I_\hbar\cap \hbar\K[X][[\hbar]])=I$, and the natural homomorphism $\A_\hbar/(\hbar)\rightarrow(\K[X]/I)^G$ is an isomorphism. Furthermore, if $\widetilde{G}/G$ is reductive, then there is a $(\widetilde{G}/G)\times\K^\times$-equivariant isomorphism $\A_\hbar\cong \K[X\red G][[\hbar]]$ of $\K[[\hbar]]$-modules, so we get a $\widetilde{G}/G$-invariant homogeneous star-product on $\K[X\red G][[\hbar]]$.
Finally, a map $ (\widetilde{\g}/\g)^*\rightarrow \A_\hbar$ sending $\xi $ to the image of $\widehat{H}_\xi$
 in  $\K[X][[\hbar]]/\I_\hbar$  (the image is easily seen to be $G$-invariant) is a quantum comoment map.
\end{Prop}
\begin{proof}Clearly, $\A_\hbar$ is complete and $\I_\hbar/(\I_\hbar\cap \hbar\K[X][[\hbar]])=I$. Flatness of $\A_\hbar$
stems from  the equality $\hbar\I_\hbar=\I_\hbar\cap \hbar\K[X][[\hbar]]$. This equality follows from the fact
that $H_{\xi_1},\ldots, H_{\xi_n}$ form a regular sequence in $\K[X]$ (here $\xi_1,\ldots,\xi_n$ is a basis
in $\g$), compare, for example, with the proof of   \cite[Lemma 3.6.1]{Wquant}.

Let us prove that $\A_\hbar/(\hbar)=(\K[X]/I)^G$. Set $\M_\hbar:=\K[X][[\hbar]]/\I_\hbar, \M_{\hbar,k}:=\M_\hbar/(\hbar^k)$.
We need to prove that the natural map $(\M_{\hbar,k})^G\rightarrow (\M_{\hbar,k-1})^G$ is surjective for any $k$.
We can rewrite $(\M_{\hbar,k})^G$ as $H^0(\n, M_{\hbar,k})^L$, where $G=N\leftthreetimes L$ is a Levi
decomposition. The group $L$ is reductive and all $L$-modules in consideration are locally finite.
From the exact sequence $0\rightarrow \K[X]/I\rightarrow \M_{\hbar,k}\rightarrow \M_{\hbar,k-1}\rightarrow 0$  we see that the following sequence
$$0\rightarrow H^0(\n,\K[X]/I)^L\rightarrow H^0(\n,\M_{\hbar,k})^L\rightarrow H^0(\n,\M_{\hbar,k-1})^L\rightarrow H^1(\n,\K[X]/I)^L$$
is exact.

Set, for brevity, $X_0:=\mu^{-1}(0)$, so that $\K[X]/I=\K[X_0]$.
We claim that $H^1(\n,\K[X_0])^L=0$.
Assume first that the quotient map $\pi:X_0\rightarrow X_0/G$ is a trivial principal $G$-bundle.
In particular, $\K[X_0]\cong \K[G]\otimes \K[X_0]^G$
as a $G$-module and hence $\K[X_0]\cong \K[N]\otimes \K[X_0]^N$ as an $N$-module. Since the group $N$ is unipotent,
the first cohomology $H^1(\n,\K[N])$ vanishes, compare, for example, with \cite[5.3]{GG}. So $H^1(\n,\K[X_0])=\{0\}$.

Proceed to the general case. There is a faithfully flat \'{e}tale morphism
$\varphi: Y\rightarrow X_0/G$ such that the natural morphism $Y\times_{X_0/G}X_0\rightarrow Y$
is a  trivial $G$-bundle. Then $H^1(\n,\K[Y\times_{X_0/G}X_0])^L=\K[Y]\otimes_{\K[X_0/G]} H^1(\n,\K[X_0])^L$.
As we have seen above, the left hand side of the previous equality vanishes.
Since the morphism $Y\rightarrow X_0/G$ is faithfully flat, we see that $H^1(\n,\K[X_0])^L=\{0\}$.

So we have proved that $\A_\hbar/(\hbar)=\K[X\red G]$.
Since the algebra $\A_\hbar$ is complete, from here it follows that $\A_\hbar\cong \K[X\red G][[\hbar]]$ as $\K[[\hbar]]$-modules and an isomorphism can be made $\widetilde{G}/G$-equivariant if $\widetilde{G}/G$
is reductive (because we can average an isomorphism over $\widetilde{G}/G$). The claim about a quantum comoment map
is straightforward to check.
\end{proof}

\section{W-algebras}\label{SECTION_Walg}
\subsection{Definition of W-algebras}\label{SUBSECTION_Walg_def}
In this subsection we are going to recall the definition of a W-algebra given
in \cite[Subsection 3.1]{Wquant} and also a variant of Premet's definition, \cite{Premet1}.

Recall that $G$ denotes a connected reductive algebraic group, and $\g$ is the Lie algebra
of $G$. Fix a nilpotent element $e\in \g$. Set $\Orb:=Ge$. Choose an
$\sl_2$-triple $(e,h,f)$ in $\g$ and set $Q:=Z_G(e,h,f)$. Let $T$ denote a maximal
torus in $Q$. Fix a $G$-invariant symmetric form $(\cdot,\cdot)$ on $\g$ and identify
$\g$ with $\g^*$ using it.

Define the Slodowy slice $S:=e+\z_\g(f)$. It will be convenient for
us to consider $S$ as a subvariety in $\g^*$.  Let $\gamma:\K^\times\rightarrow G$
have the same meaning as in Subsection \ref{SUBSECTION_Ham_class}.

Consider the cotangent bundle $T^*G$ of $G$.
Identify $T^*G$ with $G\times\g^*$ assuming that $\g^*$ consists of left-invariant 1-forms.
The variety $T^*G$ is equipped with a $G\times Q\times \K^\times$-action,
where $G$ acts by left translations: $g.(g_1,\alpha)=(gg_1,\alpha)$, $Q$ acts by right translations:
$q.(g_1,\alpha)=(g_1q^{-1}, q.\alpha)$, while $\K^\times$
acts by $$t.(g_1,\alpha)=(g_1\gamma^{-1}(t),t^{-2}\gamma(t)\alpha),$$ where $t\in \K^\times, q\in Q, g,g_1\in G,\alpha\in\g^*$.

 Recall the canonical symplectic form $\widetilde{\omega}$ on $T^*G$.
This form is $G\times Q$-invariant and $t.\widetilde{\omega}=t^2\widetilde{\omega}$.
Both $G$ and $Q$-actions are Hamiltonian with moment maps $\mu_G(g,\alpha)=g\alpha,
\mu_Q(g,\alpha)=\alpha|_{\q}$.

By the {\it equivariant Slodowy slice} we mean the subvariety $X:=G\times S\hookrightarrow G\times \g^*=T^*G$.
It turns out that $X$ is a $G\times Q\times \K^\times$-stable symplectic subvariety of $T^*G$, see \cite[Subsection 3.1]{Wquant}. Let $\omega$ denote the restriction of $\widetilde{\omega}$ to $X$. It is easy to see that the
symplectic variety  $X$ is nothing else but the model variety $ M_G(\{1\}, e,\{0\})$ introduced in \cite{slice}.

Choose a $G\times Q\times \K^\times$-invariant symplectic connection $\nabla$ on $X$
and construct a star-product on $\K[X][[\hbar]]$ from $\nabla$ and $\Omega=0$. It turns
out, see \cite[Subsection 3.1]{Wquant}, that $\widetilde{\Walg}_\hbar:=\K[X][\hbar]$ is a subalgebra in the quantum algebra $\K[X][[\hbar]]$. Set $\Walg_\hbar:=\widetilde{\Walg}_\hbar^G$. This $\K[\hbar]$-algebra is called
a {\it homogeneous  W-algebra}. We define a {\it W-algebra} $\Walg$ by $\Walg:=\Walg_\hbar/(\hbar-1)$.
This is a filtered associative algebra equipped with a $Q$-action and a quantum comoment
map $\q\rightarrow \Walg$. Also from the quantum comoment map $\g\rightarrow \widetilde{\Walg}_\hbar$
we get a homomorphism $\Centr\rightarrow \Walg$, where $\Centr$ stands for the center
of the universal enveloping algebra $\U:=U(\g)$.

There is another (earlier) definition of $\Walg$ due to Premet \cite{Premet1}.
Let us recall it. Introduce a grading on $\g$ by eigenvalues of $\ad h$:
$\g:=\bigoplus \g(i), \g(i):=\{\xi\in\g| [h,\xi]=i\xi\}$
so that $\gamma(t)\xi=t^i\xi$ for $\xi\in \g(i)$. Define the element $\chi\in \g^*$ by
$\chi=(e,\cdot)$ and the skew-symmetric form
$\omega_\chi$ on $\g(-1)$ by $\omega_\chi(\xi,\eta)=\langle\chi,[\xi,\eta]\rangle$.
It turns out that this form is symplectic. Pick a $\t$-stable lagrangian subspace $l\subset \g(-1)$
and define the subalgebra $\m:=l\oplus\bigoplus_{i\leqslant -2}\g(i)$. Then
$\chi$ is a character of $\m$. Define the shift $\m_{\chi}=\{\xi-\langle\chi,\xi\rangle,\xi\in \m\}\subset \g\oplus\K$.
Essentially, in \cite{Premet1} the W-algebra was defined as the quantum Hamiltonian reduction
$(\U/\U\m_\chi)^{\ad \m}$ (this variant of a Hamiltonian reduction is slightly
different from the one recalled above).

We checked in \cite{Wquant}, see also \cite[Theorem 2.2.1]{HC}, and the discussion after it,
that both definitions agree and also that the homomorphism
$\Centr\rightarrow \Walg$ coincides with one from \cite{Premet1}
and so is an isomorphism of $\Centr$ with the center of $\Walg$, \cite[footnote to Question 5.1]{Premet2}.

A useful feature of Premet's construction is that it allows to construct functors
between the categories of $\U$- and $\Walg$-modules. We say that a left $\U$-module $M$ is
a {\it Whittaker} module if $\m_\chi$ acts on $M$ by locally
nilpotent endomorphisms. In this case $M^{\m_\chi}=\{m\in M| \xi
m=\langle\chi,\xi\rangle m, \forall \xi\in\m\}$ is a nonzero $\Walg$-module. As
Skryabin proved in the appendix to \cite{Premet1}, the functor
$M\mapsto M^{\m_\chi}$ is an equivalence between the category of
Whittaker $\U$-modules and the category $\Walg$-$\Mod$ of
$\Walg$-modules. A quasiinverse equivalence
is given by $N\mapsto \Sk(N):=(\U/\U\m_\chi)\otimes_\Walg N$,
where $\U/\U\m_\chi$ is equipped with a natural structure of a
$\U$-$\Walg$-bimodule. In the sequel we will call $\Sk$ the {\it Skryabin functor}.

\subsection{Decomposition theorem and the correspondence between ideals}\label{SUBSECTION_decomp}
The decomposition theorem, roughly, says
that, up to suitably understood completions, the universal enveloping algebra is decomposed into the tensor product
of the W-algebra and of a Weyl algebra. We start with the equivariant version of this theorem, \cite[Subsection 3.3]{Wquant},
because it is similar in spirit to some constructions of the present paper.

First of all, let us note that we can apply the Fedosov construction (with the trivial 2-form
$\Omega$) to get a $G\times G$-equivariant homogeneous (with respect to the action of $\K^\times$
defined by $t.(g,\alpha)=(g,t^{-2}\alpha)$) star-product on $\K[T^*G][[\hbar]]$.
Automatically, $\K[T^*G][\hbar]$ is a subalgebra in $\K[T^*G][[\hbar]]$. The algebra
$\U_\hbar:=\K[T^*G][\hbar]^G$ of $G$-invariants (for the  action by left translations)
is said to be the homogeneous universal enveloping algebra of $\g$. The reason for the terminology
is that the quotient $\U_\hbar/(\hbar-1)$ is isomorphic to the universal
enveloping algebra $\U$, see \cite[Example 3.2.4]{Wquant}.

Set $V:=[\g,f]$. Equip $V$ with the symplectic form
$\omega(\xi,\eta)=\langle\chi,[\xi,\eta]\rangle$, the action of
$\K^\times, t.v:=\gamma(t)^{-1}v$, and the action of $Q$ restricted from $\g$.
Then we can equip the space $\W_{V,\hbar}:=S(V)[\hbar]$ with the Moyal-Weyl star-product,
compare with \cite[Example 3.2.3]{Wquant}.
The algebra $\W_{V,\hbar}$ is called the {\it homogeneous Weyl algebra}
of the vector space $V$. The quotient $\W_V:=\W_{V,\hbar}/(\hbar-1)$ is the usual
Weyl algebra of $V$.

Set $x:=(1,\chi)\in X$. Since the star-products on both $\K[X][[\hbar]]$ and $\K[T^*G][[\hbar]]$ are differential,
we can extend them to the completions $\K[X]^\wedge_{Gx}[[\hbar]], \K[T^*G]^\wedge_{Gx}[[\hbar]]$
along the $G$-orbit $Gx$. Next, we can form the completion $\W_{V,\hbar}^\wedge:=\K[V^*]^\wedge_0[[\hbar]]$
of the Weyl algebra. This space also has a natural star-product.
We remark that the algebras $\K[X]^\wedge_{Gx}[[\hbar]], \K[T^*G]^\wedge_{Gx}[[\hbar]],$ $ \W^\wedge_{V,\hbar}$
have natural topologies (of completions).

The following theorem was proved in \cite[Theorem 2.3.1, Remark 2.3.2]{HC}.

\begin{Thm}\label{Thm_decomp}
There is a $G\times Q\times \K^\times$-equivariant $\K[[\hbar]]$-linear Hamiltonian isomorphism $\Phi_\hbar:\K[T^*G]^\wedge_{Gx}[[\hbar]]
\rightarrow \W_{V,\hbar}^\wedge \widehat{\otimes}_{\K[[\hbar]]} \K[X]^\wedge_{Gx}[[\hbar]]$
of topological algebras.
\end{Thm}

Taking the $G$-invariants in the algebras from Theorem \ref{Thm_decomp}, we get a non-equivariant decomposition theorem. Namely, set $\Walg_\hbar^\wedge:=\K[S]^\wedge_\chi[[\hbar]],\U_\hbar^\wedge:=\K[\g^*]^\wedge_\chi[[\hbar]]$.

\begin{Cor}\label{Cor_decomp}
There is a $\K[[\hbar]]$-linear $Q\times \K^\times$-equivariant
Hamiltonian isomorphism $\Phi_\hbar:\U^\wedge_\hbar\rightarrow
\W_{V,\hbar}^\wedge\widehat{\otimes}_{\K[[\hbar]]}\Walg^\wedge_\hbar$
of topological algebras.
\end{Cor}

This proposition allows to define a map from the set $\Id(\Walg)$ of two-sided
ideals of $\Walg$ to the analogous set $\Id(\U)$ for $\U$. Namely, take a
two-sided ideal $\I\subset \Walg$. As we noted in \cite[Subsection 3.2]{Wquant},
 there is a unique ideal $\I_\hbar\subset
\Walg_\hbar$ with the following properties:
 \begin{itemize} \item $\I=\I_\hbar/(\hbar-1)$,
 \item $\I_\hbar$ is $\K^\times$-stable,
 \item $\I_\hbar$ is {\it $\hbar$-saturated}, i.e., $\I_\hbar\cap \hbar
\Walg_\hbar=\hbar \I_\hbar$.\end{itemize} Let $\I^\wedge_\hbar$ denote the closure of
$\I_\hbar$ in $\Walg^\wedge_\hbar$. Set $\J_\hbar:=\Phi_\hbar^{-1}(\W_{V,\hbar}^\wedge\widehat{\otimes}_{\K[[\hbar]]}\I^\wedge_\hbar)\cap\U_\hbar$. Finally, set
$\I^\dagger:=\J_\hbar/(\hbar-1)\subset \U$.

Reversing the procedure, we can construct a map $\J\mapsto \J_{\dagger}:\Id(\U)\rightarrow \Id(\Walg)$,
see \cite[Subsection 3.1]{HC}. More precisely, we first construct the homogeneous ideal $\J_\hbar\subset\U_\hbar$,
then complete it to get the ideal $\J^\wedge_\hbar\subset \U_\hbar^\wedge$. This ideal has the form
$\W^\wedge_\hbar(\I^\wedge_\hbar):=\W^\wedge_\hbar\widehat{\otimes}_{\K[[\hbar]]}\I^\wedge_\hbar$ for a unique
ideal $\I^\wedge_\hbar\subset \Walg_\hbar^\wedge$. The ideal $\I^\wedge_\hbar$ is the completion of a unique ideal $\I_\hbar\subset\Walg_\hbar$. Then for $\J_\dagger$ we take the image of $\I_\hbar$ in $\Walg$.
This map restricts to a surjection:
\begin{itemize}
\item from the set of all $\J\in \Id(\U)$ such that the   associated variety $\VA(\U/\J)$ equals $\overline{\Orb}$
\item to the set of all $Q$-stable ideals of finite codimension in $\Walg$.
\end{itemize}

The map  given by $\I\mapsto\I^\dagger$ is a section of $\J\mapsto \J_\dagger$, see \cite[Theorem 1.2.2]{HC}.
Moreover, $\codim_{\Walg}\J_\dagger$ coincides with the multiplicity
$\mult_{\Orb}\U/\J$. This follows from \cite[Proposition 3.4.2]{Wquant}.
In particular, $\J$ has multiplicity 1 on $\Orb$ if and only if $\J_\dagger$ is the annihilator
of a one-dimensional module. Moreover, if $V$ is a $\Walg$-module with
$\Ann_{\Walg}(V)^\dagger=\J$, then $\J_\dagger\subset \Ann_{\Walg}(V)$, see \cite[Theorem 1.2.2,(ii),  Proposition 3.4.4]{Wquant}.
So if $\mult_{\overline{\Orb}}\U/\J=1$, then $\dim V=1$ and $V$ is stable with respect to
the action of $Q$ on the set of irreducible $\Walg$-modules.


\subsection{Category $\OCat$}\label{SUBSECTION_Ocat}
Recall the Cartan subalgebra $\t\subset \q$.
Consider the centralizer $\g_0:=\z_\g(\t)$ of $\t$ in $\g$. This is a minimal Levi
subalgebra of $\g$   containing $e$. Fix an element $\theta$ lying in the cocharacter lattice
$\Hom(\K^\times,T)\hookrightarrow \t$  with $\z_\g(\theta)=\g_0$.
Let $\Walg:=\bigoplus_{\alpha\in \Z}\Walg_\alpha$ be the decomposition into eigenspaces of $\ad\theta$.
Set $$\Walg_{\geqslant 0}:=\bigoplus_{\alpha\geqslant 0}\Walg_\alpha,\quad \Walg_{>0}:=\bigoplus_{\alpha>0} \Walg_\alpha,\quad
\Walg_{\geqslant 0}^+:=\Walg_{\geqslant 0}\cap \Walg\Walg_{>0}.$$
Clearly, $\Walg_{\geqslant 0}$ is a subalgebra of $\Walg$, while $\Walg_{>0}$ and $\Walg_{\geqslant 0}^+$ are two-sided
ideals of $\Walg_{\geqslant 0}$. Note that we have an embedding $\t\hookrightarrow \q\hookrightarrow \Walg$. The image
of $\t$ lies in $\Walg_{\geqslant 0}$.  Hence we get  a homomorphism $\t\rightarrow \Walg_{\geqslant 0}/\Walg_{\geqslant 0}^+$ of Lie algebras.

Consider also the W-algebra $\Walg^0$ associated with the pair $(\g_0,e)$ (we remark that this notation is
different from \cite{LOCat}, where $\Walg^0$ denoted the quotient $\Walg_{\geqslant 0}/\Walg_{\geqslant 0}^+$,
while the W-algebra of $(\g_0,e)$ was denoted by $\underline{\Walg}$). Again, we have a natural
embedding $\t\hookrightarrow \Walg^0$. We are going to describe a relation between $\Walg_{\geqslant 0}/\Walg_{\geqslant 0}^+$ and $\Walg^0$. Fix a Cartan subalgebra $\h\subset \g_0$ containing $h$. Let $\Delta$ be the corresponding root system of $\g$. Fix a system $\Pi\subset \Delta$ of simple roots such that $\theta$ is dominant (in particular, $\Pi$ contains a system
of simple roots for $\g_0$). Set $\Delta^+:=\{\alpha\in \Delta| \langle\alpha,\theta\rangle>0\},\Delta^-=-\Delta^+$. Following \cite[Subsection 4.1]{BGK}, define an  element $\delta\in \h^*$ by
\begin{equation}\label{eq:4.1}
\delta:=\sum_{\alpha\in \Delta^-, \langle\alpha,h\rangle=-1}\alpha/2+\sum_{\alpha\in \Delta^-,\langle\alpha,h\rangle\leqslant -2}\alpha.
\end{equation}

Now we can state the following proposition, see \cite[Remark 5.4]{LOCat}.  A similar result was obtained by Brundan, Goodwin and Kleshchev in \cite[Theorem 4.3]{BGK}.

\begin{Prop}\label{Prop:4.2}
There is a $T$-equivariant isomorphism  $\Psi:\Walg_{\geqslant 0}/\Walg_{\geqslant 0}^+\rightarrow \Walg^0$ making the following diagram commutative.

\begin{picture}(60,30)
\put(2,2){$\t$}
\put(2,22){$\t$}
\put(46,22){$\Walg_{\geqslant 0}/\Walg_{\geqslant 0}^+$}
\put(50,2){$\Walg^0$}
\put(3,20){\vector(0,-1){14}} \put(4,12){\tiny $x\mapsto x-\langle\delta,x\rangle$}
\put(52,20){\vector(0,-1){14}}
 \put(53,12){\tiny $\Psi$}
\put(6,24){\vector(1,0){39}}
\put(6,4){\vector(1,0){42}}
\end{picture}
\end{Prop}


Below we consider $\Walg_{\geqslant 0}/\Walg_{\geqslant 0}^+$-modules as $\Walg^0$-modules and vice-versa
using the isomorphism $\Psi$.

By definition, the category $\widetilde{\OCat}(\theta)$ for the pair $(\Walg,\theta)$ consists of all $\Walg$-modules
$N$ satisfying the following conditions:

\begin{itemize}
\item $N$ is finitely generated.
\item $\t$ acts on $N$ by diagonalizable endomorphisms.
\item $\Walg_{> 0}$ acts on $N$ by locally nilpotent endomorphisms.
\end{itemize}

In \cite{LOCat} this category was denoted by $\widetilde{\Ocat}^{\t}(\theta)$.

There are analogs of Verma modules in $\widetilde{\Ocat}^{\t}(\theta)$.
Take a finitely generated $\Walg^0$-module $N^0$ with diagonalizable action of
$\t$. The module $\Verm^\theta(N^0):=\Walg\otimes_{\Walg_{\geqslant 0}}N^0$ lies in $\widetilde{\OCat}(\theta)$;
here we consider $N^0$ as a $\Walg_{\geqslant 0}$-module via the natural projection $\Walg_{\geqslant 0}\twoheadrightarrow
\Walg_{\geqslant 0}/\Walg_{\geqslant 0}^+\cong \Walg^0$. Also we have a functor between  $\Walg$-$\Mod$ and $\Walg^0$-$\Mod$ defined by $N\mapsto N^{\Walg_{>0}}:=\{v\in N| xv=0, \forall x\in \Walg_{>0}\}$.

Suppose now that $N^0$ is irreducible. Then $\Verm^\theta(N^0)$
has a unique irreducible quotient $L^\theta(N^0)$, see \cite[Theorem 4.5]{BGK}. Moreover, the modules
$L^\theta(N^0)$ form a complete set of pairwise non-isomorphic simple objects in  $\widetilde{\OCat}(\theta)$.
In particular,
any irreducible finite dimensional $\Walg$-module is of the form $L^\theta(N^0)$, where $N^0$ is  finite dimensional
and irreducible. However, $L^\theta(N^0)$ may be infinite dimensional even if $N^0$ is finite dimensional.
In  Subsection \ref{SUBSECTION_highest_findim} we will obtain a criterium for $L^\theta(N^0)$ to be finite dimensional. For this we will
need to recall the main result of \cite{LOCat}.

Set $\U^0:=U(\g_0)$.
Recall the subalgebra $\m\subset\g$ defined in Subsection \ref{SUBSECTION_Walg_def} and $\chi:=(e,\cdot)$. Let $\m_0:=\m\cap\g_0$. Then $\m_0$ plays the same role
for $(\g_0,e)$ as $\m$ played for $(\g,e)$. Let $\Sk_0:\Walg^0$-$\Mod\rightarrow \U^0$-$\Mod$ be the Skryabin
functor for $\g_0,e$.

Let $\g:=\bigoplus_{\alpha\in\Z} \g_\alpha$ be the grading
by eigenspaces of $\ad\theta$. Set \begin{align*}&\n_+:=\bigoplus_{\alpha>0}\g_\alpha,\quad\p:=\bigoplus_{\alpha\geqslant 0}\g_\alpha=\g_0\oplus\n_+,\\&\widetilde{\m}:=\n_+\oplus \m_0,\quad
\widetilde{\m}_\chi:=\{\xi-\langle\chi,\xi\rangle| \xi\in \widetilde{\m}\}.\end{align*}

Let $M$ be a $\U$-module. We say that $M$ is a {\it generalized Whittaker module} (for $e$ and $\theta$) if:
\begin{enumerate}
\item $M$ is finitely generated.
\item $\t$ acts on $M$ by diagonalizable endomorphisms.
\item $\widetilde{\m}_\chi$ acts by locally nilpotent endomorphisms.
\end{enumerate}

For example, let $N^0$ be a $\Walg^0$-module with diagonalizable action of $\t$. Set $\Verm^{e,\theta}(N^0):=
\U\otimes_{U(\p)}\Sk_0(N^0)$. Here $U(\p)$ acts on $\Sk_0(N^0)$ via the natural projection $U(\p)\twoheadrightarrow
\U^0$. Then $\Verm^{e,\theta}(N^0)$ is a generalized Whittaker module.
We denote the category of generalized Whittaker modules by $\widetilde{\operatorname{Wh}}(e,\theta)$.

We have a functor from $\widetilde{\operatorname{Wh}}(e,\theta)$ to $\Walg^0$-$\Mod$ constructed
as follows. For a $\g$-module $M$ the space $M^{\n_+}$ is a Whittaker $\g_0$-module. Now $\Sk_{0}^{-1}(M^{\n_+})=M^{\widetilde{\m}_\chi}$ is a $\Walg^0$-module.

The following theorem follows from  \cite[Theorem 4.1]{LOCat}.

\begin{Thm}\label{Thm:4.4}
There is  an equivalence of abelian categories $\mathcal{K}:\widetilde{\OCat}(\theta)\rightarrow \widetilde{\operatorname{Wh}}(e,\theta)$ having the following properties:
\begin{enumerate}
\item The functors $N\mapsto N^{\Walg>0}$ and $N\mapsto \mathcal{K}(N)^{\widetilde{\m}_\chi}$
from $\Walg$-$\Mod$ to $\Walg^0$-$\Mod$ are isomorphic.
\item The functors $N^0\mapsto\Verm^\theta(N^0)$ and $N^0\mapsto \mathcal{K}^{-1}(\Verm^{e,\theta}(N^0))$
from $\Walg^0$-$\Mod$ to $\Walg$-$\Mod$ are isomorphic.
\item For any $M\in \widetilde{\OCat}(\theta)$ we have $\Ann_{\U} \mathcal{K}(M)=(\Ann_{\Walg}(M))^{\dagger}$.
\end{enumerate}
\end{Thm}

\begin{Cor}\label{Cor:4.4.1}
The functor $\mathcal{K}$ induces a bijection between the following two sets:
\begin{itemize}
\item The set of finite dimensional irreducible $\Walg$-modules.
\item The set of irreducible modules $M$ in $\widetilde{\operatorname{Wh}}(e,\theta)$
with $\VA(\U/\Ann_\U(M))=\overline{\Orb}$.
\end{itemize}
This bijection sends $L^\theta(N^0)$ to a unique irreducible quotient
$L^{e,\theta}(N^0)$ of $\Verm^{e,\theta}(N^0)$ and, moreover, $\Ann_{\Walg}(L^\theta(N^0))^\dagger=
\Ann_{\U}(L^{e,\theta}(N^0))$.
\end{Cor}

\section{1-dimensional representations via category $\OCat$}\label{SECTION_O_to_onedim}
\subsection{Highest weights of finite dimensional representations}\label{SUBSECTION_highest_findim}
We use the notation of Subsection \ref{SUBSECTION_Ocat}. In particular, let
$\g_0:=\z_\g(\t)$ and let $\Walg^0$ be the W-algebra
constructed from the pair $(\g_0,e)$.

Recall the element
$\theta\in \Hom(\K^\times,T)\hookrightarrow \t$, the subalgebras $\n_+,\p\subset \g$, the category
$\widetilde{\OCat}(\theta)$ and also the Verma modules $\Verm^\theta(N^0)$ and  the irreducible modules
$L^\theta(N^0)$ in this category introduced in Subsection \ref{SUBSECTION_Ocat}.

Let $\I\mapsto \I^{\dagger_0}$ denote the map  $\Id(\Walg^0)\rightarrow \Id(\U^0)$
defined analogously to $\bullet^\dagger$,  where  $\U^0:=U(\g_0)$.
Let $L(\lambda)$ (resp., $L_{0}(\lambda)$) denote the irreducible highest weight module for
$\g$ (resp., $\g_0$) with highest weight $\lambda$. Set $J(\lambda):=\Ann_{\U}(L(\lambda)), J_0(\lambda)=\Ann_{\U^0}(L_0(\lambda))$.

\begin{Thm}\label{Thm:6.1}
Let $N^0$ be an irreducible $\Walg^0$-module. Suppose that $(\Ann_{\Walg^0}(N^0))^{\dagger_0}=J_0(\lambda)$
for some $\lambda\in\h^*$.
Then $\left(\Ann_{\Walg}(L^\theta(N^0))\right)^{\dagger}=J(\lambda)$.
\end{Thm}

This theorem has the following immediate corollary.

\begin{Cor}\label{Cor:6.2}
In the notation of Theorem \ref{Thm:6.1} the following conditions are equivalent:
\begin{enumerate}
\item $\dim L^\theta(N^0)<\infty$.
\item $\dim N^0<\infty$ and $\VA(\U/J(\lambda))=\overline{\Orb}$.
\end{enumerate}
\end{Cor}

In the proof of the theorem a crucial role is played by the following lemma.
\begin{Lem}\label{Lem:6.3}
Let $M^0$ be an irreducible $\U^0$-module. Then the induced $\U$-module
$\U\otimes_{U(\p)}M^0$ has a unique irreducible quotient $M$. The annihilator of $M$ depends only on the
annihilator of $M^0$.
\end{Lem}
The first claim of this lemma is completely standard. The idea of the proof of the second one was communicated
to me by David Vogan.
\begin{proof}
Let $\alpha$ be the eigenvalue of $\theta$ on $M^0$. For any other eigenvalue   $\beta$ of $\theta$ on $M$
the difference $\alpha-\beta$ is a positive integer. There is the largest submodule of $\U\otimes_{U(\p)}M^0$
contained in $\sum_{\beta<\alpha}M_\beta$, say $R(M^0)$. Since $\U\otimes_{U(\p)}M^0$ is generated  by $M^0$,
we see that $R(M^0)$ is the largest proper submodule in $\U\otimes_{U(\p)}M^0$, hence the first claim.

To prove the second claim let us define a certain map  $\pi:\U\twoheadrightarrow \U^0$. Namely, let
$\pi$ be a unique $T$-equivariant linear map such that $\pi$ is the identity on $\U^0\subset \U$ and $\ker\pi\cap \U_0=\U_0\cap\U\p$
(recall that $\U_\beta$ denotes the  eigenspace of $\ad\theta$ in $\U$
with eigenvalue $\beta$). We claim that $\Ann_\U M$ coincides with
\begin{equation}\label{eq:6.4}\I:=\{u\in \U| \pi(aub)\in \Ann_{\U^0}M^0,\forall a,b\in \U\}\end{equation}
Let us note  that $u\in \U_0$ acts on $M^0$ by $\pi(u)$. Let $u\in \Ann_\U M\cap \U_\beta$. Then for $a\in \U_\gamma,
b\in \U_{-\beta-\gamma}$ the element $aub$  lies in $\Ann_\U M\cap\U_0$ and acts trivially on $M$ and, in particular,
on $M^0$. So $aub\in \I$. It follows that $\Ann_\U M\subset \I$.

Let us show that $\I\subset \Ann_\U(M)$, i.e., any element  $u\in\I$ acts trivially on $M$.
We may assume that $u\in U_\beta$ for some $\beta$. Then $aub$ acts trivially on $M^0$ for any
$a\in \U_\gamma, b\in \U_{-\beta-\gamma}$. It follows
that $\U u \U M^0$ has zero intersection with $M^0=M_{\alpha}$. Since $M$ is irreducible, we see that $\U u \U M^0=\{0\}$ whence $u\in \Ann_{\U}(M)$.
\end{proof}

\begin{proof}[Proof of Theorem \ref{Thm:6.1}]
By Theorem \ref{Thm:4.4}, $\left(\Ann_{\Walg}(L^\theta(N^0))\right)^{\dagger}=\Ann_\U L^{e,\theta}(N^0)$.
Note  that $L^{e,\theta}(N^0)$ is obtained from $M^0:=\Sk_0(N^0)$ as described in Lemma \ref{Lem:6.3}.
The  $\U$-module $L(\lambda)$  is obtained from the  $\U^0$-module
$L_0(\lambda)$  by the same construction. Applying Lemma \ref{Lem:6.3}, we complete the proof.
\end{proof}

\subsection{Highest weights for 1-dimensional representations}\label{SUBSECTION_highest_1dim}
Recall that the group $Q$ acts on $\Walg$ by automorphisms and so acts also on the set of
isomorphism classes of  irreducible $\Walg$-modules.
The action of  $Q^\circ$ on the latter is trivial.

Recall that in Subsection \ref{SUBSECTION_Ocat} we chose the Cartan subalgebra $\h\subset \g$ and the set of simple roots $\Pi\subset\h^*$ and then defined the element $\delta\in\h^*$ by
(\ref{eq:4.1}).

\begin{Thm}\label{Thm:6.5}
Let $N^0$ be an irreducible $\Walg^0$-module such that $\dim L^\theta(N^0)<\infty$.
Suppose that the stabilizer of $L^\theta(N^0)$ in $N_Q(\t)$ acts on $\t$ without
nonzero fixed points.
Then the following conditions are equivalent:
\begin{enumerate}
\item $\dim N^0=1$ and $\t$ acts on $N^0$ by $\delta|_{\t}$.
\item $\dim L^\theta(N^0)=1$.
\end{enumerate}
\end{Thm}
\begin{proof}
Consider the action of $\t$ on $L^\theta(N^0)$ via the embedding $\t\hookrightarrow \q\hookrightarrow \Walg$.
By the construction of $L^\theta(N^0)$,
$\t$ acts on $N^0\subset L^\theta(N^0)$  by a single character, say $\alpha$. For any other $\t$-weight $\beta$ of $L^\theta(N^0)$
we have $\langle\beta,\theta\rangle<\langle\alpha,\theta\rangle$. By Proposition \ref{Prop:4.2}, the second condition
of (1) is equivalent to $\alpha=0$.

 Let $\Gamma$ denote the stabilizer of $L^\theta(N^0)$ in $N_Q(\t)$.
Some central extension of $\Gamma$ acts on $L^\theta(N^0)$. So the set of $\t$-weights
in $L^\theta(N^0)$ is $\Gamma$-stable.

Suppose (1) holds.  So $\langle \beta,\theta\rangle<0$
for all other weights $\beta$ of $L^\theta(N^0)$. On the other hand, $\sum_{\gamma\in \Gamma/T}\gamma.\theta=0$
for $\Gamma $ has no fixed points in $\t$. So $0$ is the only weight of $\t$ in $L^\theta(N^0)$ whence $L^\theta(N^0)=N^0$.

Now let $\dim L^\theta(N^0)=1$. Then $\alpha$ is $\Gamma$-stable.
So $\alpha$ is zero.
\end{proof}

If $\g=\sl_n$, then $Q/Z(G)$ is a connected group with a nontrivial center (unless $e$ is principal).
So the theorem does not work in this situation. However, in this case the classification
of one-dimensional $\Walg$-modules should follow from results of Brundan and Kleshchev,
\cite{BK1}, where an explicit presentation of $\Walg$ in terms of generators and relations
is found, see \cite[Subsection 3.8]{Premet4}.

For the  types $B$ and $C$ we have the following corollary.

\begin{Cor}\label{Cor:6.6}
Let $\g$ be of type $B$ or $C$ and let $\lambda\in \h^*$ be such that $\VA(\U/J(\lambda))=\overline{\Orb}$
and $\VA(\U^0/J_0(\lambda))=\overline{\Orb_0}$, where $\Orb_0$ is the orbit of $e$ in $\g_0$.
Then the following conditions (1) and (2) are equivalent.
\begin{enumerate}
\item $\mult_\Orb \U/J(\lambda)=1$.
\item
\begin{itemize}
\item[(A)]
$\lambda-\delta$ vanishes on $\t$,
\item[(B)]$\mult_{\Orb_0}\U^0/J_0(\lambda)=1$,
\item[(C)]
and for any $\lambda'\in \h^*$ the conditions
$\VA(\U^0/J_0(\lambda'))=\overline{\Orb_0}$ and $J(\lambda)=J(\lambda')$ imply
$J_0(\lambda)=J_0(\lambda')$.
\end{itemize}
\end{enumerate}
\end{Cor}
\begin{proof}
First of all, let us check that $N_Q(\t)$ acts on $\t$ without nonzero fixed points. Let $e$ correspond to the partition $(n_1^{r_1},\ldots,n_k^{r_k}),n_1>n_2>\ldots>n_k$ (the superscripts $r_i$ denote the multiplicities).
We may assume  that $G$ is $\SO_n$ ($n$ is odd) or $\Sp_n$ ($n$ is even). Then, as a linear group, $Q$ is the subgroup in $G_1\times\ldots \times G_k$ consisting of all matrices with determinant 1. Here $G_i=\operatorname{O}_{r_i}$
if either $n_i$ is even and $G=\Sp_n$ or $n_i$ is odd and $G=\SO_n$. Otherwise, $G_i=\Sp_{r_i}$.
This description implies that $N_Q(\t)$ acts on $\t$ without nonzero fixed vectors.

Let us prove the implication (1)$\Rightarrow$(2). Let $N^0$ be a finite dimensional irreducible $\Walg^0$-module with
$\Ann_{\Walg^0}(N^0)^{\dagger_0}=J_0(\lambda)$. Then, by Theorem \ref{Thm:6.1}, $\Ann_{\Walg}(L^\theta(N^0))^\dagger=J(\lambda)$. Since $\mult_{\Orb}\U/J(\lambda)=1$, the discussion at the end of
Subsection \ref{SUBSECTION_decomp} implies that
$\dim L^\theta(N^0)=1$ and $L^\theta(N^0)$ is $Q$-stable. By the choice of $N^0$, $\t$ acts on
$N^0$ by $\lambda$, hence $\lambda-\delta$ vanishes on $\t$. Since $L^\theta(N^0)$ is $Q$-stable,
we see that $N^0$ is $Q\cap G_0$-stable (where $G_0$ stands for the Levi subgroup of $G$ with Lie algebra $\g_0$). Also $\dim N^0=1$. This implies $\mult_{\Orb_0}(\U^0/J_0(\lambda))=1$. Now let $\lambda'$ be such that $J(\lambda)=J(\lambda')$,
$\VA(\U^0/J_0(\lambda'))=\overline{\Orb}_0$. Let $N'^0$ be an irreducible $\Walg^0$-module with
$\Ann_{\Walg^0}(N'^0)^{\dagger_0}=J_0(\lambda')$. Then $\Ann_\Walg(L^\theta(N'^0))^\dagger=J(\lambda')=J(\lambda)=
\Ann_{\Walg}(L^\theta(N^0))^\dagger$. Therefore $\Ann_\Walg(L^\theta(N'^0))$ is $Q$-conjugate to
$\Ann_\Walg(L^\theta(N^0))$. But the last ideal coincides with $J(\lambda)_\dagger$ and hence is $Q$-stable.
So $L^\theta(N'^0)\cong L^\theta(N^0)$ hence $N'^0\cong N^0$. We conclude that $J_0(\lambda)=\Ann_{\Walg^0}(N^0)^\dagger= \Ann_{\Walg^0}(N'^0)^\dagger=J_0(\lambda')$.

To prove (2)$\Rightarrow$(1) one reverses the argument of the previous paragraph. Namely, pick
$N^0$ with $\Ann_{\Walg^0}(N^0)=J_0(\lambda)$. Using (B) one sees that $\dim N^0=1$ and $\Ann_{\Walg^0}(N^0)$
is $G_0\cap Q$-stable. From (C) one deduces that $\Ann_{\Walg}L^\theta(N^0)$ is $Q$-stable.
Now one uses (A) and Theorem 6.5 to show that $\dim L^\theta(N^0)=1$. (1) follows.
\end{proof}

For many orbits in type $D$ the group $N_Q(\t)$ still acts on $\t$ without nonzero fixed points (and for these orbits
both claims of the corollary hold). However, this is not always the case (for instance, recall that $\mathfrak{so}_6\cong \sl_4$). The partitions for which $N_Q(\t)$ has a nonzero fixed vector in $\t$ are precisely those
where there are exactly two odd parts and they are equal.

 On the other hand Theorem \ref{Thm:6.5} works (for any representation $L^\theta(N^0)$)  whenever
 $\q$ is semisimple.  This is the case when $e$ is a rigid nilpotent element. For classical $\g$ this follows from the
 combinatorial description of rigid elements, see, for example, \cite[Corollary 7.3.5]{CM}. For the case of an exceptional Lie algebra $\g$ see the following table.

\setlongtables
\begin{longtable}{|c|c|c|c|c|}
\caption{Rigid elements in exceptional algebras}\label{Tbl:1.4}\\\hline
 N&$\g$&$e$&$\q$&$\dim \z_\g(e)$\\\endfirsthead\hline
 N&$\g$&$e$&$\q$&$\dim\z_\g(e)$\\\endhead\hline
 1&$G_2$&$A_1$&$A_1$&8\\\hline
 2&$G_2$&$\widetilde{A}_1$&$A_1$&6\\\hline
 3&$F_4$&$A_1$&$C_3$&36\\\hline
 4&$F_4$&$\widetilde{A}_1$&$A_3$&30\\\hline
 5&$F_4$&$A_1+\widetilde{A}_1$&$A_1+A_1$&24\\\hline
 6&$F_4$&$A_2+\widetilde{A}_1$&$A_1$&18\\\hline
 7&$F_4$&$\widetilde{A}_2+A_1$&$A_1$&16\\\hline
 8&$E_6$&$A_1$&$A_5$&56\\\hline
 9&$E_6$&$3A_1$&$A_2+A_1$&38\\\hline
 10&$E_6$&$2A_2+A_1$&$A_1$&24\\\hline
 11&$E_7$&$A_1$&$D_6$&99\\\hline
 12&$E_7$&$2A_1$&$B_4+A_1$&81\\\hline
 13&$E_7$&$(3A_1)'$&$C_3+A_1$&69\\\hline
 14&$E_7$&$4A_1$&$C_3$&63\\\hline
 15&$E_7$&$A_2+2A_1$&$3A_1$&51\\\hline
 16&$E_7$&$2A_2+A_1$&$2A_1$&43\\\hline
 17&$E_7$&$(A_3+A_1)'$&$3A_1$&41\\\hline
 18&$E_8$&$A_1$&$E_7$&190\\\hline
 19&$E_8$&$2A_1$&$B_6$&156\\\hline
 20&$E_8$&$3A_1$&$F_4+A_1$&136\\\hline
 21&$E_8$&$4A_1$&$C_4$&120\\\hline
 22&$E_8$&$A_2+A_1$&$A_5$&112\\\hline
 23&$E_8$&$A_2+2A_1$&$B_3+A_1$&102\\\hline
 24&$E_8$&$A_2+3A_1$&$G_2+A_1$&94\\\hline
 25&$E_8$&$2A_2+A_1$&$G_2+A_1$&86\\\hline
 26&$E_8$&$A_3+A_1$&$B_3+A_1$&84\\\hline
 27&$E_8$&$2A_2+2A_1$&$B_2$&80\\\hline
 28&$E_8$&$A_3+2A_1$&$B_2+A_1$&76\\\hline
 29&$E_8$&$D_4(a_1)+A_1$&$3A_1$&72\\\hline
 30&$E_8$&$A_3+A_2+A_1$&$2A_1$&66\\\hline
 31&$E_8$&$2A_3$&$B_2$&60\\\hline
 32&$E_8$&$A_4+A_3$&$A_1$&48\\\hline
 33&$E_8$&$A_5+A_1$&$2A_1$&46\\\hline
 34&$E_8$&$D_5(a_1)+A_2$&$A_1$&46\\\hline
 \end{longtable}

The information in this table is taken from \cite[Subsection 5.7]{McG}, and \cite[Subsection 13.1]{Carter}.
A nilpotent element is given by its Bala-Carter label. This label also indicates a minimal Levi subalgebra
containing the element. In all cases but NN29,34 the nilpotent element is regular in the Levi subalgebra. In the remaining
two cases its $D_l$-component, $l=4,5$, is subregular.

\subsection{Toolkit}\label{SUBSECTION_toolkit}
Recall that we have fixed a Cartan subalgebra $\h\subset\g$ and a system of simple roots $\Pi$.
We assume that $\g$ is an exceptional Lie algebra and $\Orb$ is its rigid nilpotent orbit. To prove that
the corresponding algebra $\Walg$ has a one-dimensional representation we need to find
\begin{itemize}
\item A Levi subalgebra $\g_0$, whose simple root system $\Pi_0$ is contained in $\Pi$ and  a nilpotent element $e\in \Orb\cap\g_0$ that is not contained in a proper Levi subalgebra of
$\g_0$. These data are needed to construct the category $\Ocat$, see Subsection \ref{SUBSECTION_Ocat}.
Let $\Orb_0$ denote the orbit of $e$ in $\g_0$.
\item  $\lambda\in\h^*$ such that
\begin{itemize}
\item[(A)] $\VA(\U^0/J_0(\lambda))=\overline{\Orb}_0$ so that there is an irreducible finite dimensional
$\Walg^0$-module $N^0$ with $\Ann_{\Walg^0}(N^0)=J_0(\lambda)$. This yields $\Orb\subset
\VA(\U/J(\lambda))$.
\item[(B)] $\VA(\U/J(\lambda))=\overline{\Orb}$ or, equivalently, modulo (A), $\dim \VA(\U/J(\lambda))\leqslant \dim \Orb$. By Corollary \ref{Cor:6.2}, this means $\dim L^\theta(N^0)<\infty$.
\item[(C)] $\lambda-\delta\in \operatorname{Span}_{\K}\Pi_0$, where $\delta$ is defined by (\ref{eq:4.1}).
By the discussion in Subsection \ref{SUBSECTION_highest_1dim}, this means that $L(N^0)=N^0$.
\item[(D)] $J_0(\lambda)=\I_0^{\dagger_0}$ for some ideal $\I_0$ of $\Walg^0$ of codimension 1
(this condition always holds when $e$ is regular in $\g_0$, because in this case $\Walg^0$ is commutative,
and so all irreducible $\Walg^0$-modules are 1-dimensional). This means that $\dim N^0=1$.
\end{itemize}
\end{itemize}

Recall the notion of a {\it special} nilpotent orbit due to Lusztig. One of their characterizations
is that an orbit $\Orb$ is special if and only if there is an integral  weight $\lambda\in \h^*$  such that
$\overline{\Orb}=\VA(\U/J(\lambda))$.

 Consider the Langlands
dual algebra $\g^\vee$ with the Cartan subalgebra $\h^\vee:=\h^*$ and the root system $\Delta^\vee$, the dual
root system of $\g$. Let $\Pi^\vee\subset \Delta^\vee$ be the simple root system corresponding to $\Pi$. Note that $\g\cong\g^\vee$ provided $\g$ is exceptional. Spaltenstein constructed an order-reversing bijection $\Orb\mapsto\Orb^\vee$ between
the sets of special nilpotent orbits.

Suppose $e$ is a special element. Take $e^\vee\in \Orb^\vee$ and let $(e^\vee, h^\vee, f^\vee)$ be the corresponding
$\sl_2$-triple. We may assume that $h^\vee$ lies in $\h^*$ and is dominant. This specifies $h^\vee$ uniquely.
As Barbash and Vogan checked in \cite[Proposition 5.1]{BV}, the element $e^\vee$ is even (i.e., all eigenvalues
of $\ad h^\vee$ are even) provided $e$ is rigid.

Let, as usual, $\rho$ denote  half the sum of all positive roots of $\g$.

\begin{Prop}[\cite{BV}, Proposition 5.10]\label{Prop:6.6}
One has the equality $\VA(\U/J(h^\vee-\rho))=\overline{\Orb}$.
\end{Prop}

The list of special
nilpotent orbits in exceptional algebras and the description of the Spaltenstein duality are given in \cite[Subsection 13.4]{Carter}. In Table \ref{Tbl:1.5} we list the special {\it rigid} elements $e$ together with their duals
$e^\vee$. Here numbers in the first column are those from Table \ref{Tbl:1.4}.

\begin{longtable}{|c|c|c|c|}
\caption{Non-minimal special rigid elements in exceptional algebras}\label{Tbl:1.5}\\\hline
 N&$\g$&$e$&$e^\vee$\\\endfirsthead\hline
 N&$\g$&$e$&$e^\vee$\\\endhead\hline
 4&$F_4$&$\widetilde{A}_1$&$F_4(a_1)$\\\hline
 5&$F_4$&$A_1+\widetilde{A}_1$&$F_4(a_2)$\\\hline
 12&$E_7$&$2A_1$&$E_7(a_2)$\\\hline
 15&$E_7$&$A_2+2A_1$&$E_7(a_4)$\\\hline
 19&$E_8$&$2A_1$&$E_8(a_2)$\\\hline
 22&$E_8$&$A_2+A_1$&$E_8(a_4)$\\\hline
 23&$E_8$&$A_2+2A_1$&$E_8(b_4)$\\\hline
 29&$E_8$&$D_4(a_1)+A_1$&$E_8(a_6)$\\\hline
 \end{longtable}

It turns out that for 5 special elements (with exception of $(F_4,\widetilde{A}_1),(E_8,A_2+A_1),(E_8,D_4(a_1)+A_1)$) the weight $\lambda=h^\vee-\rho$ satisfies conditions (A) and (C)
(under a suitable choice of $\g_0$). (B) follows from Proposition \ref{Prop:6.6}, while (D)
 holds because the five elements are principal in $\g_0$. Hence in these 5 cases $J(h^\vee-\rho)$ has the form
$\I^\dagger$ with $\dim \Walg/\I=1$.

Now suppose $\Orb$ is not special. Pick some  $\lambda\in \h^*$.
Let $\Delta(\lambda)^\vee$ be the subset of all coroots $\alpha^\vee\in \Delta^\vee$
such that $\langle\alpha^\vee,\lambda\rangle\in \Z$. Let $\Delta(\lambda)$ be the dual root system of $\Delta^\vee(\lambda)$. Note that for $\g=E_6,E_7,E_8$ the root system $\Delta(\lambda)$ is a subsystem in $\Delta$.

Let $\Pi^\vee(\lambda)$ be a unique system of simple roots in $\Delta^\vee(\lambda)$ such that the corresponding system of positive coroots in $\Delta^\vee(\lambda)$ consists of positive coroots for $\Delta^\vee$. Let $W(\lambda)$ denote the Weyl group of $\Delta(\lambda)$ (or of $\Delta^\vee(\lambda)$). Let $w\in W(\lambda)$ be a unique element of minimal length such that $w(\lambda+\rho)$ is antidominant, i.e. $\langle w(\lambda+\rho),\alpha^\vee\rangle\leqslant 0$ for all $\alpha^\vee\in \Pi^\vee(\lambda)$.

Recall that the group $W(\lambda)$ is decomposed into equivalence classes called {\it two-sided} cells and there is
a bijection between the set of two-sided cells in $W(\lambda)$ and the set of special nilpotent orbits
in the Lie algebra $\g(\lambda)$ associated with $\h,\Delta(\lambda)$. We will denote the nilpotent orbit
corresponding to $w\in W(\lambda)$ by $\Orb_{w}$.

\begin{Prop}\label{Prop:6.7}
We have $\dim \VA(\U/J(\lambda))=\dim\g-\dim \g(\lambda)+\dim \Orb_w$.
\end{Prop}
This result is well known to  specialists. However, as far as we know, the proof was never written
explicitly, so it is written  down below. I wish to thank A. Joseph and D. Vogan for explaining the
details.
\begin{proof}
Using an appropriate  translation functor, see, for example, \cite[4.12,4.13]{Jantzen},
we may assume that $\lambda$ is regular.
Let $w_0$ denote the longest element in $W(\lambda)$. Then $\lambda=w_0w\cdot \mu$, where $\mu$ is dominant.
Here for $x\in W(\lambda)$ we set $x\cdot \mu:=x(\mu+\rho(\lambda))-\rho(\lambda)$, where $\rho(\lambda)$ is
half of the sum of all positive roots in $\Delta(\lambda)$.

For $x,y\in W(\lambda)$ let $a(x,y)$ denote the coefficient of the character of the Verma module $\Verm(y\cdot \mu)$
in the character of the irreducible module $L(x\cdot \mu)$. Set $n=\frac{1}{2}(\# \Delta-\dim \VA(\U/J(\lambda)))$.  Joseph proved in \cite{JosephII} that $n$ coincides with the smallest non-negative integer $m$ such that
$$\sum_{y\in W(\lambda)}a(x,y)y^{-1}\rho(\lambda)^m\neq 0.$$
But the Kazhdan-Lusztig conjecture implies that $a(x,y)$ coincides with the coefficient of the character of the Verma module $\Verm_{(\lambda)}(y\cdot \mu)$
in the character of the irreducible module $L_{(\lambda)}(x\cdot \mu)$ for the algebra $\g(\lambda)$. It follows
that $n=\#\Delta(\lambda)-\dim \Orb_w$. Since $\#\Delta-\#\Delta(\lambda)=\frac{1}{2}(\dim\g-\dim\g(\lambda))$, we are done.
\end{proof}

The easiest case here is when $w$ is the longest element in $W(\lambda)$. Then $\Orb_{w}$ is just the zero orbit.
So we have the following corollary originally suggested to us by A. Premet.

\begin{Cor}[\cite{Joseph}, Corollary 3.5]\label{Cor:6.8}
If $\langle \lambda+\rho,\alpha^\vee\rangle>0$ for all $\alpha^\vee\in \Pi(\lambda)^{\vee}$, then $\dim \VA(\U/J(\lambda))=\dim\g-\dim \g(\lambda)$.
\end{Cor}

To finish the section let us provide an example of computation.

As we have seen, it is more convenient to work with the element $\lambda+\rho$ rather than with $\lambda$.
Also from the point of view of computations it is better to replace $\delta$ with  the element
$$\delta':=\frac{1}{2}\sum_{\alpha\in\Delta^+,\langle\alpha,h\rangle=0,1}\alpha.$$
We claim that $\delta'-\delta-\rho\in \Span_{\Q} \Pi_0$. Indeed, since $e,h,f\in \g_0$, we see that the $\sl_2$-triple
$(e,h,f)$ preserves all weight subspaces of $\t$. In particular, by the representation
theory of $\sl_2$, for any $k,l$ we have
$$\sum_{\alpha,\langle\alpha,\theta\rangle=k,\langle\alpha,h\rangle=l}\alpha|_{
\t}=\sum_{\alpha,\langle\alpha,\theta\rangle=k,\langle\alpha,h\rangle=-l}\alpha|_{
\t}$$ or equivalently,
\begin{equation}\label{eq:5.1.1}
\sum_{\alpha,\langle\alpha,\theta\rangle=k,\langle\alpha,h\rangle=l}\alpha\equiv\sum_{\alpha,\langle\alpha,\theta\rangle=k,\langle\alpha,h\rangle=-l}\alpha
\mod\Span_\Q(\Pi_0).
\end{equation}
Applying (\ref{eq:5.1.1}), we get $$\delta\equiv \frac{1}{2}\sum_{\langle\alpha,\theta\rangle<0, \langle\alpha,h\rangle\neq 0,1}\alpha \mod \Span_\Q(\Pi_0).$$
It follows that $$\delta+\rho\equiv -\frac{1}{2}\sum_{\langle\alpha,\theta\rangle<0, \langle\alpha,h\rangle=0,1}\alpha \equiv \frac{1}{2}\sum_{\langle\alpha,\theta\rangle>0, \langle\alpha,h\rangle=0,-1}\alpha\mod \operatorname{Span}_{\Q}(\Pi_0).$$
Applying (\ref{eq:5.1.1}) once more, we have $$\sum_{\langle\alpha,\theta\rangle>0, \langle\alpha,h\rangle=-1}\alpha\equiv
\sum_{\langle\alpha,\theta\rangle>0, \langle\alpha,h\rangle=1}\alpha \mod\operatorname{Span}_{\Q}(\Pi_0). $$
It follows that $\delta'\equiv \delta+\rho \mod \operatorname{Span}_{\Q}(\Pi_0)$.

We will check the existence of  $\lambda\in \h^*$
with properties (A)-(D) for  the non-special orbit $A_5+A_1$
in $E_8$. We use the notation for
roots from \cite{VO}. Namely, fix an orthonormal basis $\varepsilon_i,i=1,\ldots,9,$ in the Euclidean
space $\mathbb{R}^9$. Then we represent the real form $\h(\mathbb{R})^*$ of $\h^*$ as the quotient of $\mathbb{R}^9$ by the diagonal.
The image of $\varepsilon_i$ in $\h(\mathbb{R})^*$ is again denoted by $\varepsilon_i$. Now simple roots
are $\alpha_i:=\varepsilon_i-\varepsilon_{i+1}, i=1,\ldots,7, \alpha_8=\varepsilon_6+\varepsilon_7+\varepsilon_8$.
The element $\rho$ equals $\sum_{i=1}^7 (8-i)\varepsilon_i-22\varepsilon_9$. Let $\pi_i, i=1,\ldots,8,$ denote the
fundamental weights.  They are given by $\pi_i=\sum_{i=1}^7 \varepsilon_i-\min(i,15-2i)\varepsilon_9, i=1,2,\ldots,7, \pi_8=3\varepsilon_9$. We identify $\h(\mathbb{R})^*$ with $\h(\mathbb{R})$ by means of the scalar product
on $\h(\mathbb{R})^*$ (induced from $\mathbb{R}^9$).

First of all, we choose $\g_0$ such that $\Pi_0=\{\alpha_1,\alpha_2,\alpha_3,\alpha_4,\alpha_5,\alpha_7\}$. Then for $h$ we can take
$5\varepsilon_1+3\varepsilon_2+\varepsilon_3-\varepsilon_4-3\varepsilon_5-5\varepsilon_6+\varepsilon_7-\varepsilon_8$.
We get
$$\delta'=\frac{3}{2}\varepsilon_1+\varepsilon_2+2\varepsilon_3+\varepsilon_4+\frac{3}{2}\varepsilon_5+
\frac{1}{2}\varepsilon_6+\varepsilon_7-4\varepsilon_9.$$
Set $$\lambda':=\varepsilon_1+\frac{7}{6}\varepsilon_2+\frac{1}{3}\varepsilon_3+\frac{1}{2}\varepsilon_4
+\frac{2}{3}\varepsilon_5+\frac{5}{6}\varepsilon_6+\frac{1}{6}\varepsilon_7-\frac{1}{6}\varepsilon_8-\frac{9}{2}\varepsilon_9,
\lambda:=\lambda'-\rho.$$

One can check directly that $\lambda'-\delta'\in \Q\Pi_0$, which is equivalent to (C). The restriction of $\lambda$ to $\g_0$ is antidominant, so (A) holds.  The system $\Delta(\lambda)$ is isomorphic to $A_5+A_2+A_1$ and the set of simple roots in $\Delta(\lambda)$ is $\varepsilon_7+\varepsilon_4+\varepsilon_3,
\varepsilon_2-\varepsilon_7,\varepsilon_8+\varepsilon_7+\varepsilon_1,\varepsilon_6-\varepsilon_8,
\varepsilon_8+\varepsilon_5+\varepsilon_4,\varepsilon_8+\varepsilon_6+\varepsilon_3,
\varepsilon_7+\varepsilon_5+\varepsilon_2,\varepsilon_1+\varepsilon_3+\varepsilon_5$.
The element $\lambda'$ is positive on all these roots. By Corollary \ref{Cor:6.8}, $\dim \VA(\U/J(\lambda))=248-46=202=\dim \Orb$, which is (B).
So we have checked   (A),(B),(C). The condition (D) is vacuous because $e$ is of principal
Levi type.

\section{Parabolic induction for W-algebras}\label{SECTION_induction}
\subsection{Lusztig-Spaltenstein induction}\label{SUBSECTION_LSind}
Let $\underline{G}$ be a Levi subgroup of $G$, $\underline{\g}$ its Levi subalgebra, and
$\underline{e}\in \underline{\g}$ be a nilpotent element. Let $\underline{\Orb}\subset \underline{\g}\cong \underline{\g}^*$ be the
orbit of $\underline{e}$.

\begin{defi}[\cite{LS}, Introduction, Theorem 2.2]\label{defi:2.1}
Let $P$ be a parabolic subgroup of $G$ with Levi subgroup $\underline{G}$ and
let $\p$ be the Lie algebra of $P$. There is a unique nilpotent
orbit $\Orb$ in $\g$ such that
$(\underline{\Orb}+\Rad_u(\p))\cap\Orb$ is dense in $\underline{\Orb}+\Rad_u(\p)$.
The orbit $\Orb$ does not depend on the choice of $P$ and
is said to be induced from $\underline{\Orb}$ and is denoted by
$\Ind_{\underline{\g}}^\g(\underline{\Orb})$.
\end{defi}

We remark that, although Lusztig and Spaltenstein considered
unipotent orbits in $G,\underline{G}$, their results are transferred to the case
of nilpotent orbits in a straightforward way, see \cite{McG} for an
exposition.

Here are some properties of $\Ind_{\underline{\g}}^\g(\underline{\Orb})$ due to Lusztig and
Spaltenstein.

\begin{Prop}[\cite{LS}, Theorem 1.3]\label{Prop:2.2}
Let $e\in (\underline{\Orb}+\Rad_u(\p))\cap\Ind_{\underline{\g}}^\g(\underline{\Orb})$. Then the following
assertions hold:
\begin{enumerate}
\item $\dim Z_G(e)=\dim Z_{\underline{G}}(\underline{e})$.
\item $Ge\cap (\underline{\Orb}+\Rad_u(\p))=Pe$.
\item $Z_G(e)^\circ\subset P$.
\item If $g\in G$ is such that $\Ad(g)e\in \underline{\Orb}+\Rad_u(\p)$, then
there is $z\in Z_G(e)$ such that $z^{-1}g\in P$.
\end{enumerate}
\end{Prop}

In the sequel we will also need the following lemma.

\begin{Lem}\label{Prop:2.3}
Let $h_0\in\g$ be such that $[h_0,e]=2e$. Then $h_0\in\p$.
\end{Lem}
\begin{proof}
Since $\z_\g(e)\subset \p$, we may assume that $h_0=h$ is the
semisimple element of some $\sl_2$-triple $(e,h,f)$. Let
$\gamma:\K^\times\rightarrow G$ be the one-parameter subgroup corresponding to
$h$, i.e. such that $\frac{d}{dt}|_{t=0}\gamma=h$. By Proposition
\ref{Prop:2.2}, (4), for any $t\in \K^\times$ there is $z_t\in
Z_G(e)$ such that $z_t^{-1}\gamma(t)\in P$. Since
$Z_G(e)^\circ\subset P$, we see that $\{\gamma(t)P, t\in \K^\times\}$ is a finite subset
of $G/P$. Therefore $\gamma(t)\in P$ for all $t\in \K^\times$,
equivalently, $h\in\p$.
\end{proof}

\subsection{Construction on the classical level}\label{SECTION_classical}
Let $\underline{G},\underline{e}, P,\Orb$ be such as in the previous subsection. Set $N:=\Rad_u(P)$.
Choose $e\in \Orb\cap (\underline{\Orb}+\n)$. Let $(e,h,f)$ be an $\sl_2$-triple.
Also recall the group $Q:=Z_G(e,h,f)$.
By Proposition \ref{Prop:2.2} and Lemma \ref{Prop:2.3}, we have $h\in \p,\q\subset \p$. The subalgebra
$\K h\oplus \q\subset\p$ is reductive, so replacing $\underline{G}$ by a $P$-conjugate subgroup,
we may assume $h\in \underline{\g},\q \subset \underline{\g}$. Hence $Q^\circ\subset \underline{G}$.

Our goal in this subsection is to study the Hamiltonian reductions
of some  open subsets in $X,T^*G$.

Let $P^-$ denote the opposite parabolic to $P$ containing $\underline{G}$ and $N^-:=\Rad_u(P^-)$. Set $Z=T^*(PN^-)=
PN^-\times\g^*\subset G\times \g^*=T^*G$.
This is an affine $P$ and $\K^\times$-stable open subvariety in $T^*G$ containing $(1,\chi)$.
Note that $Z= P\times N^-\times\g^*$. Set $\widetilde{Y}:=Z\red N, Y:=(Z\cap X)\red N=\widetilde{Y}\cap (X/N)$.

Since $\underline{G}$ normalizes $N$,   we have Hamiltonian actions of $\underline{G}$  on $Y,\widetilde{Y}$.
The  Kazhdan actions of $\K^\times$ on $X,T^*G$ descend to
$\K^\times$-actions on $Y,\widetilde{Y},$ because
$Z$ is $\K^\times$-stable and $\mu_{G}(t.x)=t^2\mu_{G}(x)$ for all $x\in X$.
 The reduction $\widetilde{Y}$ can be naturally identified with
the Hamiltonian $\underline{G}$-variety $T^*P^-$. The latter, in its turn, is naturally
isomorphic to $T^*\underline{G}\times T^*N^-$, where $\underline{G}$ acts only on the first
factor (however, $\K^\times$  acts non-trivially on both factors).  Set $\underline{Q}:=Q\cap \underline{G}$. Clearly, $Z\subset T^*G$ is stable with respect to the right $\underline{Q}$-action.

Let $y$ denote the image of $(1,\chi)$ in $Y\subset \widetilde{Y}$ (the inclusion $(1,\chi)\in \mu^{-1}(\p)$
holds because $\chi$ vanishes on $\n$). Note that the orbit $\underline{G}y$ is stable under the $\K^\times$-action because $\gamma(t)\in \underline{G}$ for any $t\in\K^\times$.

According to Proposition \ref{Prop:2.41}, there is a $G\times Q\times\K^\times$-equivariant Hamiltonian isomorphism $\varphi:(T^*G)^\wedge_{Gx}\rightarrow
(X\times V^*)^\wedge_{Gx}$, which is identical on $Gx$.
Completing further, we get a $P\times \underline{Q}\times \K^\times$-equivariant
symplectomorphism $(T^*G)^\wedge_{Px}\rightarrow (X\times V^*)^\wedge_{Px}$ and, since $Px\subset Z$,
also a symplectomorphism $Z^\wedge_{Px}\rightarrow [(X\cap Z)\times V^*]^\wedge_{Px}$, which we also
denote by $\varphi$. Taking the reduction, we get a $\underline{G}\times \underline{Q}\times \K^\times$-equivariant Hamiltonian isomorphism $\underline{\varphi}:\widetilde{Y}^\wedge_{\underline{G}y}\rightarrow
(Y\times V^*)^\wedge_{\underline{G}y}=Y^\wedge_{\underline{G}y}\times (V^*)^{\wedge}_{0}$.



Pick an $\sl_2$-triple $(\underline{e},\underline{h},\underline{f})$ in $\underline{\g}^{\underline{Q}}$ such that
$[h,\underline{h}]=0, [h,\underline{f}]=-2\underline{f}$ (if $\underline{e}=0$, we set $\underline{h}=\underline{f}=0$).

Now let $\underline{X}$ denote the equivariant Slodowy slice constructed for $\underline{G}$ and $\underline{e}$. In other words, as a $\underline{G}$-variety $\underline{X}$ is nothing else but $\underline{G}\times \underline{S}$, where $\underline{S}$ is the Slodowy slice for $\underline{e}$ in $\underline{\g}$. Since $\dim \z_\g(e)=\dim \z_{\underline{\g}}(\underline{e})$,
we see that $\dim \underline{S}=\dim S$. Let $\underline{\chi}\in \underline{\g}^*$ be such that $\langle\underline{\chi},\xi\rangle=(\underline{e},\xi)$
for $\xi\in\underline{\g}$.  The group $\underline{Q}$ acts
on $\underline{X}$ and its action commutes with the actions of $\underline{G}$ and $\K^\times$.

Set $\underline{x}=(1,\underline{\chi})$. We have the projection $\widetilde{Y}=T^*P^-=T^*\underline{G}\times T^*N^-
\twoheadrightarrow T^*\underline{G}$ mapping $y$ to $\underline{x}$. So it gives rise to the projection $\widetilde{Y}^\wedge_{\underline{G}y}\twoheadrightarrow
(T^*\underline{G})^\wedge_{\underline{G}\underline{x}}$.

\begin{Lem}\label{Lem:2.3.1}
There is a $\underline{G}\times \underline{Q}\times \K^\times$-equivariant Hamiltonian isomorphism  $\underline{X}^\wedge_{\underline{G}\underline{x}}\rightarrow Y^\wedge_{\underline{G}y}$.
\end{Lem}
\begin{proof}
According to Proposition \ref{Prop:2.41}, we need to check that
\begin{enumerate}
\item The stabilizers of $\underline{x},y$ in $\underline{G}\times \underline{Q}\times \K^\times$ coincide.
\item  $\mu_{\underline{G}}(\underline{x})=\mu_{\underline{G}}(y), \mu_{\underline{Q}}(\underline{x})=\mu_{\underline{Q}}(y)$, where $\mu_{\underline{G}},\mu_{\underline{Q}}$
are the moment maps  for the $\underline{G}$- and $\underline{Q}$-actions.
\item The normal spaces to $\underline{G}\underline{x}$ in $\underline{X}$ and to
$\underline{G}y$ in $Y$ are isomorphic as $(\underline{G}\times \underline{Q}\times \K^\times)_{\underline{x}}$-modules.
\end{enumerate}

The stabilizer of $x$ in $G\times \underline{Q}\times \K^\times$ is $\{(q\gamma(t),q,t), q\in \underline{Q}, t\in \K^\times\}$. Since
the latter is contained in $\underline{G}\times \underline{Q}\times\K^\times$, we see that it coincides with
the stabilizer of $y$ in $\underline{G}\times \underline{Q}\times \K^\times$. Also, this subgroup coincides
with the stabilizer of $\underline{x}$. Hence (1). Moreover, it is easy to see from the definition of the moment map
for $Y$ that $\mu_{\underline{G}}(y)=\mu_{\underline{G}}(\underline{x})=\underline{e}, \mu_{\underline{Q}}(y)=\mu_{\underline{Q}}(\underline{x})=0$. Hence (2).
Since the orbit $\underline{G}\underline{x}$ is a coisotropic subvariety in $\underline{X}$, (3) reduces  to checking that $\dim Y=\dim \underline{X}$.

From the definition of $Y$, we have $\dim Y=\dim X-2\dim N$. But $\dim X=\dim G+\dim S=\dim \underline{G}+2\dim N+\dim S$ and hence $\dim Y=\dim S+\dim\underline{G}$. On the other hand, $\dim \underline{X}=\dim \underline{S}+\dim \underline{G}$. As we have seen above, $\dim S=\dim \underline{S}$. Hence $\dim Y=\dim S+\dim \underline{G}=\dim \underline{X}$.
\end{proof}

Now we will discuss some constructions to be used in the proof of Theorem \ref{Prop:2.8}.
Set $\underline{V}:=[\underline{\g},\underline{f}] $. Consider the action of $\K^\times$ on $\underline{V}$ given by $t\mapsto \gamma(t)^{-1}$. Also $\underline{V}$ has a natural $\underline{Q}$-action.
Similarly to $\varphi$, we have a $\underline{G}\times \underline{Q}\times \K^\times$-equivariant
Hamiltonian morphism $\varphi_0:(T^*\underline{G})^\wedge_{\underline{G}\underline{x}}\rightarrow
\underline{X}^\wedge_{\underline{G}\underline{x}}\times (\underline{V}^*)^\wedge_0$.


\begin{Lem}\label{Lem:2.3}
There is a $\underline{G}\times \underline{Q}\times\K^\times$-equivariant Hamiltonian isomorphism
$$\underline{\varphi}':\widetilde{Y}^\wedge_{\underline{G}y}\rightarrow Y^\wedge_{\underline{G}y}\times (V^*)^\wedge_0$$
making the following diagram commutative.

\begin{picture}(90,50)
\put(6,2){$Y^\wedge_{\underline{G}y}$}
\put(40,2){$\underline{X}^\wedge_{\underline{G}\underline{x}}$}
\put(30,22){$\underline{X}^\wedge_{\underline{G}\underline{x}}\times (\underline{V}^*)^\wedge_0$}
\put(73,22){$(T^*\underline{G})^\wedge_{\underline{G}\underline{x}}$}
\put(75,42){$\widetilde{Y}^\wedge_{\underline{G}y}$}
\put(1,42){$Y^\wedge_{\underline{G}y}\times (V^*)^\wedge_0$}
\put(36,3){\vector(-1,0){20}}
\put(8,40){\vector(0,-1){33}}
\put(45,20){\vector(0,-1){13}}
\put(71,23){\vector(-1,0){15}}
\put(79,40){\vector(0,-1){13}}
\put(73,43){\vector(-1,0){48}}
\put(45,45){\tiny $\underline{\varphi}'$}
\put(64,24){\tiny $\varphi_0$}
\end{picture}

Here all vertical arrows are projections and the arrow $\underline{X}^\wedge_{\underline{G}\underline{x}}\rightarrow Y^\wedge_{\underline{G}y}$
is an isomorphism constructed in Lemma \ref{Lem:2.3.1}.
\end{Lem}
\begin{proof}
Set $\widetilde{V}:=T_y(\underline{G}y)^{\skewperp}/[T_y(\underline{G}y)^\skewperp\cap T_y(\underline{G}y)] $. This is a symplectic vector space acted on by $(\underline{G}\times \underline{Q}\times \K^\times)_y$.
Recall that the stabilizer $(\underline{G}\times \underline{Q}\times \K^\times)_y$
consists of all elements of the form $(q\gamma(t),q,t), q\in \underline{Q},t\in \K^\times,$
so we can identify it with $\underline{Q}\times \K^\times$. The action of $\underline{Q}\times \K^\times$ on $\widetilde{Y}=T^*(P^-)$ is given by
\begin{align*}
&q.(p,\alpha):=(qpq^{-1}, q\alpha), \, t.(p,\alpha):=(\gamma(t)p\gamma(t)^{-1}, t^{-2}\gamma(t)\alpha),\\
&q\in \underline{Q}, t\in \K^\times, p\in P^-,\alpha\in \p^{-*}.\end{align*}
So we see that as a $\underline{Q}\times \K^\times$-module, $\widetilde{V}$ is isomorphic to $\underline{V}\oplus \n^-\oplus\n^{-*}$, where the action of $\underline{Q}\times \K^\times$ on $\underline{V}$ was described previously, the action on $\n^-$ factors through the adjoint action of $\underline{G}$ via the homomorphism
$\underline{Q}\times \K^\times\rightarrow \underline{G}, (q,t)\mapsto q\gamma(t)$, and, finally,
the action on $\n^{-*}\cong \n$ is given by $(q,t)\mapsto t^{-2}\Ad(q\gamma(t))$.

Proposition \ref{Prop:2.41} implies that $$\widetilde{Y}^\wedge_{\underline{G}y}\cong \underline{X}^\wedge_{\underline{G}\underline{x}}\times (\widetilde{V}^*)^\wedge_0.$$
So to complete the proof we need to check that $V$ is isomorphic to $\widetilde{V}$ as a $\underline{Q}\times \K^\times$-module (and then an isomorphism can be chosen to be symplectic). But we have already constructed  a $\underline{G}\times \underline{Q}\times \K^\times$-equivariant isomorphism $\underline{\varphi}:\widetilde{Y}^\wedge_{\underline{G}y}\rightarrow Y^\wedge_{\underline{G}y}\times (V^*)^\wedge_{0}$ mapping $y$ to $y$. The existence of $\underline{\varphi}$ implies $V\cong \widetilde{V}$.
\end{proof}

Unfortunately, we do not know whether $\underline{\varphi}=\underline{\varphi}'$. So let us consider
$\lambda:=\underline{\varphi'}\circ\underline{\varphi}^{-1}$. This is an $\underline{G}\times \underline{Q}\times \K^\times$-equivariant Hamiltonian automorphism of $Y^\wedge_{\underline{G}y}\times (V^*)^\wedge_0$, which is trivial on $\underline{G}y$.

We are going to describe the structure of $\lambda$, more precisely, to prove that $\lambda$ is, in a sense,
an inner automorphism.

\begin{Lem}\label{Lem:2.3.2}
The isomorphism $\lambda$ can be written as a composition $A\exp(v)$, where $A\in \Sp(V^*)$
(when we consider $A$ as an automorphism of $Y^\wedge_{\underline{G}y}\times (V^*)^\wedge_0$ we mean that
$A$ acts trivially on the first factor)
 and $v$ is the Hamiltonian
vector field $v_f$ of a $\underline{G}\times \underline{Q}$-invariant element $f\in \K[Y\times V^*]^\wedge_{\underline{G}y}$ of degree 2 with respect to the $\K^\times$-action.
\end{Lem}
We will see in the proof that  $v$ will be constructed in such a way that $\exp(v)$ converges.
\begin{proof}
Since $\lambda$ is the identity on $\underline{G}y$, the subspace $(T_{y}\underline{G}y)^\skewperp$ is stable w.r.t
$d_{y}\lambda$. Let $A$ be the operator induced by $d_{y}\lambda$ on
$V^*=(T_{y}\underline{G}y)^\skewperp/([T_{y}\underline{G}y)^\skewperp\cap T_{y}\underline{G}y]$. Set $\lambda_0:=A^{-1}\circ\lambda$. Note that $d_{y}\lambda_0$ is a unipotent operator on $T_{y}(Y\times V^*)$. Therefore $\lambda_0^*$ induces a unipotent operator on all quotients $\K[Y\times V^*]/I^k$, where $I$ stands for the ideal of
$\underline{G}y$. This means that the operator $\ln \lambda_0^*:\K[Y\times V^*]^\wedge_{\underline{G}y}\rightarrow
\K[Y\times V^*]^\wedge_{\underline{G}y}$ is well-defined. This operator is a derivation, and the corresponding vector field $v$ is symplectic and annihilates all hamiltonians $H_\xi, \xi\in \underline{\g}$.
It remains to apply Lemma \ref{Lem:1.1.2} to conclude that $v$ has the required form.
\end{proof}

\subsection{Proof of the main theorem}\label{SUBSECTION_proof_main}
Consider the formal scheme $X^\wedge_{Px}$.
 Its algebra of functions is equipped with the differential star-product induced from
 $\K[X]$.  Since $Y^\wedge_{\underline{G}y}$ is naturally identified with $X^\wedge_{Px}\red N$, we can apply the construction of Subsection \ref{SUBSECTION_reduction} to get a star-product on $\K[Y]^\wedge_{\underline{G}y}$.
 Note that we have a natural map
 \begin{equation}\label{eq:2.21}\Walg^\wedge_\hbar=\K[X]^\wedge_{Gx}[[\hbar]]^G\hookrightarrow
 \K[X]^\wedge_{Px}[[\hbar]]^P\rightarrow \K[Y]^\wedge_{\underline{G}y}[[\hbar]]^{\underline{G}}.\end{equation}

The following proposition is a quantum analog of Lemma \ref{Lem:2.3.1}.

\begin{Prop}\label{Prop:2.4.1}
There is a $\underline{G}\times\underline{Q}\times\K^\times$-equivariant Hamiltonian isomorphism of quantum algebras $$\K[Y]^\wedge_{\underline{G}y}[[\hbar]]\rightarrow \K[\underline{X}]^\wedge_{\underline{G}\underline{x}}[[\hbar]].$$
\end{Prop}
\begin{proof}
By Lemma \ref{Lem:2.3.1}, we have a $\underline{G}\times\underline{Q}\times \K^\times$-equivariant Hamiltonian isomorphism $Y^\wedge_{\underline{G}y}\rightarrow \underline{X}^\wedge_{\underline{G}\underline{x}}$. It follows
from Corollary \ref{Cor:1.6} 
that the corresponding isomorphism
$\K[Y]^\wedge_{\underline{G}y}\rightarrow \K[\underline{X}]^\wedge_{\underline{G}\underline{x}}$ can be lifted  to an $\underline{G}\times\underline{Q}\times \K^\times$-equivariant Hamiltonian isomorphism $\K[Y]^\wedge_{\underline{G}y}[[\hbar]]\rightarrow \K[\underline{X}]^\wedge_{\underline{G}\underline{x}}[[\hbar]]$.
\end{proof}

Set $h_0:=h-\underline{h}$, then $h_0\in \z_{\underline{\g}}(\underline{e},\underline{h},\underline{f})$.
There is an embedding $\z_{\underline{\g}}(\underline{e},\underline{h},
\underline{f})\hookrightarrow \underline{\Walg}$ coming from the quantum  comoment map. So consider $ h_0$
as an element in $\underline{\Walg}$. Form the completion $\underline{\Walg}'$ of $\underline{\Walg}$ consisting
of all infinite sums $\sum_{i\leqslant k}f_i$, where $[h_0,f_i]=if_i$. The algebra $\underline{\Walg}'$ has a topology,
the sets $O_k:=\{\sum_{i\leqslant k}f_i\}$ are fundamental neighborhoods of 0. Clearly, $\underline{\Walg}'$
is complete and separated with respect to this topology.

\begin{Thm}\label{Thm:2.5}
There is an embedding $\Xi:\Walg\hookrightarrow \underline{\Walg}'$.
\end{Thm}
\begin{proof}
The isomorphism $\K[Y]^\wedge_{\underline{G}y}[[\hbar]]\rightarrow \K[\underline{X}]^\wedge_{\underline{G}\underline{x}}[[\hbar]]$ restricts  to a
$\K^\times$-equivariant isomorphism $\K[Y]^\wedge_{\underline{G}y}[[\hbar]]^{\underline{G}}\rightarrow \K[\underline{X}]^\wedge_{\underline{G}\underline{x}}[[\hbar]]^{\underline{G}}$. The latter isomorphism intertwines the
embeddings of $\q$. Let us note that $\K[\underline{X}]^\wedge_{\underline{G}\underline{x}}[[\hbar]]^{\underline{G}}$ is nothing else
but $\K[\underline{S}]^\wedge_{\underline{\chi}}[[\hbar]]=\underline{\Walg}^\wedge_\hbar$.
On the other hand, from the construction of the star-product on $\K[Y]^\wedge_{\underline{G}y}[[\hbar]]$ we have a $\K^\times\times \underline{Q}$-equivariant embedding
$\Xi_\hbar:\Walg_\hbar=\K[X][\hbar]^G\hookrightarrow\K[Y]^\wedge_{\underline{G}y}[[\hbar]]^{\underline{G}}$. Therefore we get an embedding $\Xi_\hbar:\Walg_\hbar\rightarrow (\underline{\Walg}^\wedge_\hbar)_{\K^\times-l.f}$,
where the subscript ``$\K^\times-l.f$'' means the subalgebra of locally finite vectors, and also homomorphisms
\begin{equation}\label{eq:2.5.1}
\Xi_1:\Walg=\Walg_\hbar/(\hbar-1)\rightarrow \underline{\Walg}^\heartsuit:=(\underline{\Walg}^\wedge_\hbar)_{\K^\times-l.f}/(\hbar-1).
\end{equation}
\begin{equation}\label{eq:2.5.2}
\Xi_0:\K[S]=\Walg_\hbar/(\hbar)\rightarrow (\underline{\Walg}^\wedge_\hbar)_{\K^\times-l.f}/(\hbar)=(\K[\underline{S}]^\wedge_{\underline{\chi}})_{\K^\times-l.f}.
\end{equation}
Let us note that both algebras in (\ref{eq:2.5.1}) are filtered and $\Xi_1$ preserves
the filtrations. The algebras and the homomorphism in (\ref{eq:2.5.2}) are associated graded of
those in (\ref{eq:2.5.1}). Also let us note that the homomorphism in (\ref{eq:2.5.2}) is injective.
To show this it is enough to check that the morphism $Y/\underline{G}\rightarrow X/G$ is dominant, that is, that
any $G$-orbit in $X$  intersects $\mu_{G}^{-1}(\p)$. The latter is clear, because any adjoint orbit
intersects $\p$.

So the homomorphism in (\ref{eq:2.5.1}) is also injective. It remains to embed the algebra $\underline{\Walg}^\heartsuit$ into $\underline{\Walg}'$. The algebra $\K[\underline{S}]$ is the symmetric algebra $A:=S(\vf)$, where
$\vf$ is $\underline{S}$ considered as a vector space (with $\underline{\chi}$ as an origin).
The space $\vf$ is equipped with two  $\K^\times$-actions. The first one is
the Kazhdan action: $t.\underline{s}=t^{-2}\underline{\gamma}(t)\underline{s}$, where $\underline{\gamma}$ is the
one-parameter subgroup of $\underline{G}$ corresponding to $\underline{h}$. The second action is given
by the one-parameter subgroup $\gamma_0(t)=\gamma(t)\underline{\gamma}(t)^{-1}$, whose differential
 at $t=0$ is $h_0$. So we have a bi-grading
$A=\bigoplus_{i,j\in \Z} A(i,j)$, where $A(i,j)=\{a\in A| t.a=t^ia,\gamma_0(t).a=t^ja\}$.
It follows that $A(0,0)=\K$, $A(i,j)=0$ for  $i<0$ and all $j$, and $A(0,k)=0$ provided $k\neq 0$.
Analogously to \cite[Subsection 3.2]{Wquant},  $\underline{\Walg}^\heartsuit$ consists of all infinite sums $a=\sum_{i,j\in \Z}a_{ij}, a_{ij}\in A(i,j),$ such that there is
$n\in \Z$ (depending on $a$) with $a_{ij}=0$ provided $i+j\geqslant n$. On the other hand,
$\underline{\Walg}'$ consists of all sums $\sum_{i,j\in \Z}a_{ij}$, such that
\begin{enumerate}
\item  there is $n$ with $a_{ij}=0$ for $j>n$,
\item and for any $j$ only finitely many elements $a_{ij}$ are nonzero.
\end{enumerate}

The products on both algebras can be written as
\begin{equation}\label{eq:2.5.3}(\sum_{i,j}a_{ij})(\sum_{i',j'}a'_{i'j'})=\sum_{i,j,i',j'}a_{ij}*_1a'_{i'j'}.\end{equation}
Here for $f,g\in \K[\underline{S}]$ we write $f*_1g=\sum_{i=0}^\infty D_i(f,g)$, where $f*g=\sum_{i=0}^\infty D_i(f,g)\hbar^{2i}$
is the star-product on $\K[\underline{S}][\hbar]$.

Since $A(i,j)=0$ for $i<0$ and $\dim\bigoplus_j A(i,j)<\infty$ for any $i$, we see  that $\underline{\Walg}^{\heartsuit}\subset \underline{\Walg}'$. By (\ref{eq:2.5.3}), this inclusion is a homomorphism of algebras.
\end{proof}

\begin{proof}[Proof of Theorem \ref{Thm:1}]
Any finite dimensional representation of $\underline{\Walg}$ has only finitely many eigenvalues of $h_0$.
 So it uniquely extends to a continuous (with respect to the discrete topology) representation of $\underline{\Walg}'$. Now the functor $\rho$ we need is just the pull-back from $\underline{\Walg}'$ to $\Walg$.
\end{proof}

We would like to deduce a corollary from the proof of Theorem \ref{Thm:2.5} concerning the behavior of the centers.
Let $\Centr,\underline{\Centr}$ denote the centers of $\U,\underline{\U}$, respectively. We have identifications of $\Centr$ (resp., $\underline{\Centr}$) with the center of $\Walg$ (resp., $\underline{\Walg}$).

Note that the choice of the parabolic subgroup $P$ with Levi subgroup $\underline{G}$ gives rise to an embedding $\iota:\Centr\rightarrow \underline{\Centr}$. More precisely, $\iota$ is the restriction of the map $\U\mapsto \underline{\U}$ defined analogously to the map from the proof of Lemma \ref{Lem:6.3}.

\begin{Cor}\label{Cor:2.9}
The image of $\Centr\subset \Walg$ under $\Xi$ lies in $\underline{\Centr}\subset
\underline{\Walg}\hookrightarrow \underline{\Walg}'$. The corresponding homomorphism
$\Centr\rightarrow \underline{\Centr}$ coincides with the one described in the previous paragraph.
\end{Cor}
\begin{proof}
The action of $G$ on $X$ is free, so the homomorphism $\U_\hbar\hookrightarrow \widetilde{\Walg}_\hbar$
induced by the quantum comoment map is an embedding.
 Set $\Centr_\hbar:=\U_\hbar^G\subset \widetilde{\Walg}_\hbar^G=\Walg_\hbar$.
It is clear that $\Centr_\hbar/(\hbar-1)=\Centr$.

The homomorphism $\Walg_\hbar\rightarrow \K[Y]^\wedge_{\underline{G}y}[[\hbar]]$ factors
through $\Walg_\hbar\rightarrow (\widetilde{\Walg}_\hbar/\widetilde{\Walg}_\hbar\n)^{\ad \n} \rightarrow\K[Y]^\wedge_{\underline{G}y}[[\hbar]]$.
Therefore the homomorphism $\Centr_\hbar=\U_\hbar^G\hookrightarrow \Walg_\hbar\rightarrow \K[Y]^\wedge_{\underline{G}y}[[\hbar]]$ factors
through $(\U_\hbar/\U_\hbar\n)^{\n}$. For any element $u\in \Centr_\hbar$ there
are unique $\underline{u}\in \underline{\Centr}_\hbar:=\underline{\U}_\hbar^{\underline{G}}$ and $u_+\in \U_\hbar\n$ such that $u=\underline{u}+u_+$.
The image of $u$ in $\K[Y]^\wedge_{\underline{G}y}[[\hbar]]$ coincides with the image of $\underline{u}$
under the quantum comoment map. This image lies in the center of $\K[Y]^\wedge_{\underline{G}y}[[\hbar]]^{\underline{G}}$.
Now the claim of the corollary follows from the construction of $\Xi$.
\end{proof}

\begin{Rem}\label{Rem:2.7}
In some special cases we can embed $\Walg$ into $\underline{\Walg}$ itself not just into a completion
of $\underline{\Walg}$. Namely, suppose that the element $h$
is even. In this case we can take $\underline{\g}=\z_\g(h)$ and $\underline{e}=0$. Then $\underline{\Walg}=\underline{\U}$, and $h_0=h$. The algebra $\underline{\Walg}'$ coincides with $\underline{\Walg}$.
So we get an embedding $\Walg\hookrightarrow \underline{\Walg}$.
Note that, at least, some embedding $\Walg\hookrightarrow \underline{\Walg}$ in this case was known
previously: this is a so called
{\it generalized Miura transform} first discovered in \cite{Lynch}. Using techniques developed
in \cite[Subsections 3.2,3.3]{Wquant} it should not be very difficult to show that these two embeddings coincide
(at least, up to an automorphism of $\underline{\U}$).
However, to save space we are not going to study this question.

We do not know if it is always possible to embed $\Walg$ into $\underline{\Walg}$. Neither it is clear
for us how to define the map $\Theta$ using Premet's construction of W-algebras (the generalized Miura transform
can be easily seen from that construction).
\end{Rem}

\subsection{Behavior of ideals}\label{SUBSECTION_ideals}
Recall the procedure of the parabolic induction for ideals in universal enveloping algebras. To a two-sided ideal $\underline{\J}\subset \underline{\U}$
we can assign a two-sided ideal $\J$ in $\U$
as follows. Let $\underline{\J}_\p$ denote the inverse image of $\underline{\J}$ in $U(\p)$ under the projection $U(\p)\twoheadrightarrow \underline{\U}$.
Then for $\J$ we take the largest (with respect to inclusion) two-sided ideal of $\U$ contained
in $\U\underline{\J}_\p$. It is well known that if $\underline{M}$ is a $\underline{\U}$-module with $\Ann_{\underline{\U}}\underline{M}= \underline{\J}$, then
$\J:=\Ann_{\U}(\U\otimes_{U(\p)}\underline{M})$. We say that $\J$ is obtained from
$\underline{\J}$ by parabolic induction.

The goal of this subsection is to relate the inclusion $\Walg\hookrightarrow \underline{\Walg}'$ to
parabolic induction for ideals.

\begin{Thm}\label{Prop:2.8}
Let $\underline{\I}$ be a two-sided ideal  in $\underline{\Walg}$ and let $\underline{\I}'$ denote the closure
of $\underline{\I}$ in $\underline{\Walg}'$. Set $\I:=\Xi'^{-1}(\underline{\I}')\subset \Walg,
\underline{\J}:=\underline{\I}^\dagger\subset \underline{\U}$. Finally, let $\J\subset \U$ be the ideal
obtained from $\underline{\J}$ by parabolic induction. Then $\I^\dagger=\J$.
\end{Thm}

Here is an immediate corollary of this theorem.

\begin{Cor}\label{Cor:2.8.1}
Let $\underline{N}$ be a finite dimensional $\underline{\Walg}$-module. Then the ideal $\Ann_{\Walg}(\rho(\underline{N}))^\dagger\subset\U$
is obtained from $\Ann_{\underline{\Walg}}(\underline{N})^\dagger\subset \underline{\U}$ by the parabolic induction.
\end{Cor}

\begin{proof}[Proof of Theorem \ref{Prop:2.8}]
Let us explain the scheme of the proof. In Step 1 we will give an alternative construction
of the parabolic induction map $\underline{\J}\mapsto \J$. The main result of Step 3 is  commutative
diagram (\ref{eq:diag1}), whose vertices are various sets of two-sided ideals.
The sink of this diagram is $\Id(\underline{\Walg})$, and the source is $\Id(\U)$.
Diagram (\ref{eq:diag1}) is obtained from diagram (\ref{eq:diag_alg}) of algebra homomorphisms
to be established on Step 2.
It is more or less easy to see that (\ref{eq:diag1}) implements the composition $\underline{\I}\mapsto \I\mapsto \I^\dagger$.
The analogous claim for  the composition $\underline{\I}\mapsto \underline{\J}\mapsto \J$
will be verified on Step 4 using the description of Step 1.

{\it Step 1.}
We can extend the star-product from $\K[T^*G]^\wedge_{Gx}$ to $\K[T^*G]^\wedge_{Px}$. Applying the
quantum Hamiltonian reduction described in Subsection \ref{SUBSECTION_reduction} to the action
of $N$ on $(T^*G)^\wedge_{Px}$, we get
the structure of a quantum algebra on $\K[\widetilde{Y}]^\wedge_{\underline{G}y}[[\hbar]]$.
 Similarly to (\ref{eq:2.21}), we have a homomorphism
$\U^\wedge_\hbar\hookrightarrow \K[\widetilde{Y}]^\wedge_{\underline{G}y}[[\hbar]]^{\underline{G}}$.

A remarkable feature here, however, is that, unlike in the case of W-algebras, we can get a non-completed version
of this homomorphism. The quantum algebra $\K[T^*G][\hbar]$ is $G\times G\times \K^\times$-equivariantly isomorphic
to the  ``homogeneous version'' of the  algebra of differential operators $\D_\hbar(G)$, for a rigorous proof see, e.g., Subsection
\ref{SUBSECTION_Fedosov_vs_TDO}. This isomorphism gives rise to a $P\times P^-\times \K^\times$-equivariant
 identification $\K[Z][\hbar]\xrightarrow{\sim} \D_\hbar(NP^-)$ (recall that $Z=P^-\times\g^*=T^*(NP^-)$). The quantum
 Hamiltonian reduction of $\D_\hbar(NP^-)=\D_\hbar(N)\otimes_{\K[\hbar]}\D_\hbar(P^-)$ under the action of $N$ is naturally identified with
 $\D_\hbar(P^-)$. So we get a $\underline{G}\times P^-\times\K^\times$-equivariant identification
 $\K[\widetilde{Y}][\hbar]\xrightarrow{\sim} \D_\hbar(P^-)$.


Embed $P^-$ into $G/N$ via $p\mapsto Np$.
Consider a  homomorphism $\U_\hbar\rightarrow \D_\hbar(P^-)$ sending $\xi\in\g$ to the velocity
vector $\xi_{G/N}$ on $P^-\hookrightarrow G/N$.
The  image of this homomorphism consists of $\underline{G}$-invariants, so we have a homomorphism
$\U_\hbar\rightarrow \K[\widetilde{Y}][\hbar]^G\xrightarrow{\sim} \D_\hbar(P^-)^{\underline{G}}=
\underline{\U}_\hbar\otimes_{\K[\hbar]} \D_\hbar(N^-)$.
%
%
The latter homomorphism is $\underline{G}$-equivariant, where we consider the adjoint action
of $\underline{G}$ on $\U$, and the $\underline{G}$-action on $\D_\hbar(P^-)^{\underline{G}}$  coming from the $\underline{G}$-action on $P^-$ by right translations.  We can also take the quotients by $\hbar-1$ and get the $\underline{G}$-equivariant homomorphism $\widetilde{\Xi}:\U\rightarrow \underline{\U}\otimes \D(N^-)$.

The reason why we need this construction is that we can give an alternative definition of the
parabolic induction for ideals in universal enveloping algebras.

\begin{Lem}\label{Lem:3.31}
Let $\underline{\J}$ be a two-sided ideal in $\underline{\U}$ and $\J$ be the ideal in $\U$
obtained from $\underline{\J}$ by parabolic induction. Then $\J=\widetilde{\Xi}^{-1}(\underline{\J}\otimes \D(N^-))$.
\end{Lem}
\begin{proof}[Proof of Lemma \ref{Lem:3.31}]
Set, for brevity, $\B:=\underline{\U}\otimes \D(N^-)$. As a vector space, $\D(N^-)=\B_+\otimes
\B_-$, where $\B_+=\K[N^-],\B_-=U(\n^-)$. Note that the construction of the homomorphism
$\widetilde{\Xi}:\U\rightarrow\B$ implies that the restrictions of $\widetilde{\Xi}$ to
$\underline{\U},\B_-\subset\U$ are the identity maps.

Pick a rational element $\vartheta\in \z(\underline{\g})$ such that all eigenvalues of $\ad\vartheta$ on $\n$
are positive (and then the eigenvalues of $\ad\vartheta$ on $\n^-$ are  automatically negative).
Recall that $\underline{G}$ acts on $\B$. Let $\vartheta_*$ denote the image of $\vartheta$ in $\operatorname{Der}(\B)$. All eigenvalues of $\vartheta_*$ on $\underline{\U}$, (resp., $\B_-,\B_+$) are zero (resp., nonpositive, nonnegative). Moreover, $\B_+=\K\oplus\B_{++}$, where all eigenvalues of $\vartheta_*$ on $\B_{++}$ are strictly positive. Clearly, $\B_{++}$ is an ideal in $\B_+$.

Pick a $\underline{\U}$-module $\underline{M}$ with $\Ann_{\underline{\U}}(\underline{M})=\underline{\J}$. Consider the induced module $M=\B\otimes_{\underline{\U}\otimes\B_+}\underline{M}$. Here $\B_+$ acts on $\underline{M}$ via the projection $\B_+\twoheadrightarrow
\B_+/\B_{++}=\K$.

We claim that, as a $\U$-module, $M=\U\otimes_{U(\p)}\underline{M}$.
Since $\B=\B_-\otimes\underline{\U}\otimes\B_+$ we see that $M=U(\n^-)\otimes \underline{M}$
as a $U(\p^-)$-module, where we consider $U(\p_-)$ as a subalgebra of $\B$.
All eigenvectors of $\vartheta_*$ with positive eigenvalues act by zero
on $\underline{M}\subset M$. Since the map $\widetilde{\Xi}$ is $\underline{G}$-equivariant, we have $\widetilde{\Xi}([\vartheta,u])=\vartheta_*u$ for all $u\in U$. Therefore $\n$
acts trivially on $\underline{M}$. The identity map $\underline{M}\rightarrow \underline{M}$
extends to a homomorphism $\U\otimes_{U(\p)}\underline{M}\rightarrow
M$ of $\U$-modules. Since both $\U\otimes_{U(\p)}\underline{M}$ and $M$ are naturally identified
with $U(\n^-)\otimes \underline{M}$, we see that this homomorphism is, in fact, an isomorphism.

It is easy to see the annihilator of $M$ in $\mathcal{B}$ coincides with $\underline{\J}\otimes \D(N^-)$.
So $\Ann_\U(M)=\widetilde{\Xi}^{-1}(\underline{\J}\otimes \D(N^-))$. On the other hand, from the previous
paragraph it follows that $\Ann_\U(M)=\J$.
\end{proof}

{\it Step 2.}
Recall the isomorphism $$\Phi_\hbar:\K[T^*G]^\wedge_{Gx}[[\hbar]]\xrightarrow{\sim} \W_{V,\hbar}^\wedge(\K[X]^\wedge_{Gx}[[\hbar]])(:=\W_{V,\hbar}^\wedge\widehat{\otimes}_{\K[[\hbar]]}\K[X]^\wedge_{Gx}[[\hbar]])$$
from Theorem \ref{Thm_decomp}. We have a similar isomorphism $\Phi_{0\hbar}:\K[T^*\underline{G}]^\wedge_{\underline{G}\underline{x}}\xrightarrow{\sim}
\W_{\underline{V},\hbar}^\wedge(\K[\underline{X}]^\wedge_{\underline{G}\underline{x}}[[\hbar]])$.

We extend $\Phi_\hbar$ to the isomorphism
$\K[T^*G]^\wedge_{Px}[[\hbar]]\xrightarrow{\sim} \W_{V,\hbar}^\wedge(\K[X]^\wedge_{Px}[[\hbar]])$ and then, performing the quantum
Hamiltonian reduction, we get
an isomorphism $\underline{\Phi}_\hbar: \K[\widetilde{Y}]^\wedge_{\underline{G}y}[[\hbar]]\xrightarrow{\sim} \W_{V,\hbar}^\wedge(\K[Y]^\wedge_{\underline{G}y}[[\hbar]])$.

The following lemma is a quantum analog of Lemmas \ref{Lem:2.3.1},\ref{Lem:2.3},\ref{Lem:2.3.2}.

\begin{Lem}\label{Prop:2.4}
There is a commutative diagram

\begin{equation}\label{eq:diag3}
\begin{picture}(150,50)
\put(6,2){$\K[Y^\wedge_{\underline{G}y}][[\hbar]]$}
\put(72,2){$\K[\underline{X}]^\wedge_{\underline{G}\underline{x}}[[\hbar]]$}
\put(70,22){$\W_{\underline{V},\hbar}^\wedge(\K[\underline{X}]^\wedge_{\underline{G}\underline{x}}[[\hbar]])$}
\put(121,22){$\K[T^*\underline{G}]^\wedge_{\underline{G}\underline{x}}[[\hbar]]$}
\put(124,42){$\K[\widetilde{Y}]^\wedge_{\underline{G}y}[[\hbar]]$}
\put(1,42){$\W_{V,\hbar}^\wedge(\K[Y]^\wedge_{\underline{G}y}[[\hbar]])$}
\put(45,42){$\W_{V,\hbar}(\K[Y]^\wedge_{\underline{G}y}[[\hbar]])$}
\put(27,3){\vector(1,0){42}}
\put(12,7){\vector(0,1){32}}
\put(80,7){\vector(0,1){13}}
\put(103,23){\vector(1,0){17}}
\put(108,25){\tiny $\Phi_{0\hbar}^{-1}$}
\put(130,27){\vector(0,1){13}}
\put(77,43){\vector(1,0){45}}
\put(90,44.5){\tiny $\underline{\Phi}_\hbar^{-1}$}
\put(33,43){\vector(1,0){11}}
\put(37,44){\tiny $\Lambda_\hbar$}
\end{picture}
\end{equation}
Here all vertical maps are natural embeddings, while $\Lambda_\hbar$ is an automorphism of the form $A\exp(\frac{1}{\hbar^2}\ad(\hat{f}))$ for some $\underline{G}\times \underline{Q}$-invariant element $\hat{f}\in \W_{V,\hbar}^\wedge(\K[Y]^\wedge_{\underline{G}y}[[\hbar]])$
of  degree 2 with respect to $\K^\times$. The bottom horizontal arrow is an isomorphism from Proposition \ref{Prop:2.4.1}.
\end{Lem}
\begin{proof}
Similarly to Lemma \ref{Lem:2.3}, we get a $\underline{G}\times \underline{Q}\times\K^\times$-equivariant automorphism $\Lambda_\hbar$ of $\W_{V,\hbar}^\wedge(\K[Y]^\wedge_{\underline{G}y}[[\hbar]])$ making the diagram commutative
 and preserving the quantum comoment map. In addition, we may assume that modulo $\hbar^2$ the automorphism $\Lambda_\hbar$ equals $A\exp(\frac{1}{\hbar^2}\ad(f))$, where $f,A$ are  as in Lemma \ref{Lem:2.3.2}. In particular, $\exp(\frac{1}{\hbar^2}\ad(f))$ converges.
So $\Lambda_{\hbar, 1}:=\Lambda_\hbar \exp(-\frac{1}{\hbar^2}\ad(f))A^{-1}$ has the form $\operatorname{id}+\sum_{i=1}^\infty
T_i\hbar^{2i}$. Again, $\Lambda_{\hbar,1}$ is a $\underline{G}\times \underline{Q}\times \K^\times$-equivariant Hamiltonian automorphism. As in the proof of   \cite[Theorem 2.3.1]{HC}, we see that $T_1$ is a $\underline{G}\times \underline{Q}$-equivariant Poisson derivation of $\K[Y\times V^*]^\wedge_{\underline{G}y}$. Also let us note that $T_1$ annihilates
the image of the classical comoment map and therefore, by Lemma \ref{Lem:1.1.2}, $T_1=v_{f_1}$ for some
$\underline{G}\times \underline{Q}$-invariant element $f_1\in \K[Y\times V^*]^\wedge_{\underline{G}y}$ that has degree 0 w.r.t $\K^\times$.
Replace $\Lambda_{\hbar,1}$ with  $\Lambda_{\hbar,2}:=\Lambda_{\hbar,1}\exp(-\ad(f_1))$. We get $\Lambda_{\hbar,2}=\operatorname{id}+\sum_{i=2}^\infty T_i' \hbar^{2i}$. Repeating the procedure we obtain
a presentation $\Lambda_{\hbar}=A\exp(\frac{1}{\hbar^2}\ad(f))\exp(\ad(f_1))\exp(\hbar^2\ad(f_2))\ldots$.
Using the Campbell-Hausdorff formula, we get a presentation of $\Lambda_\hbar$ in the  form required
in the statement of the proposition.
\end{proof}

By taking the $\underline{G}$-invariants, we get the following commutative diagram (the embeddings $\Walg_\hbar\hookrightarrow\K[Y]^\wedge_{\underline{G}y}[[\hbar]], \U_\hbar\hookrightarrow \K[\widetilde{Y}]^\wedge_{\underline{G}y}[[\hbar]]$ come from the quantum
Hamiltonian reduction).

\begin{equation}\label{eq:diag_alg}
\begin{picture}(160,70)
\put(22,2){$\K[Y]^\wedge_{\underline{G}y}[[\hbar]]^{\underline{G}}$}
\put(70,2){$\underline{\Walg}^\wedge_\hbar$}
\put(65,22){$\W_{\underline{V},\hbar}^\wedge(\underline{\Walg}^\wedge_\hbar)$}
\put(127,22){$\underline{\U}^\wedge_\hbar$}
\put(125,42){$\K[\widetilde{Y}]^\wedge_{\underline{G}y}[[\hbar]]^{\underline{G}}$}
\put(19,42){$\W_{V,\hbar}^\wedge(\K[Y]^\wedge_{\underline{G}y}[[\hbar]]^{\underline{G}})$}
\put(70,42){$\W_{V,\hbar}^\wedge(\K[Y]^\wedge_{\underline{G}y}[[\hbar]]^{\underline{G}})$}
\put(43,2.5){\vector(1,0){25}}
\put(34,7){\vector(0,1){33}}
\put(73,7){\vector(0,1){13}}
\put(86,23){\vector(1,0){40}}
\put(130,27){\vector(0,1){13}}
\put(53,43){\vector(1,0){15}}
\put(60,44){\tiny $\Lambda_\hbar$}
\put(104,43){\vector(1,0){20}}
\put(4,15){$\Walg_\hbar^\wedge$}
\put(1,62){$\W_{V,\hbar}^\wedge(\Walg^\wedge_\hbar)$}
\put(5,20){\vector(0,1){40}}
\put(10,14){\vector(1,-1){10}}
\put(12,60){\vector(1,-1){14}}
\put(110,62){$\U^\wedge_\hbar$}
\put(22,63){\vector(1,0){87}}
\put(114,60){\vector(1,-1){14}}
\end{picture}
\end{equation}

{\it Step 3}.
Recall from \cite{Wquant} that for a flat $\K[\hbar]$-algebra $\A_\hbar$ equipped with a $\K^\times$-action one defines
the set $\Id_\hbar(\A_\hbar)$ of two-sided $\hbar$-saturated $\K^\times$-stable ideals.

Now let $\A$ be an algebra equipped with an increasing exhaustive separated filtration $\F_i\A$.
Then to $\A$ we can assign the Rees algebra $R_\hbar(\A)=\bigoplus_{i\in \Z}(\F_i\A)\hbar^{i}\subset \A[\hbar^{-1},\hbar]$. The group $\K^\times$ naturally acts on $R_\hbar(\A)$ by algebra
automorphisms. Also to any two-sided ideal $\I\subset \A$ we can assign the two-sided
ideal $R_\hbar(\I)=\bigoplus_{i\in \Z}(\F_i\A\cap\I)\hbar^{i}$. As we have seen in \cite[Subsection 3.2]{Wquant},
the sets $\Id(\A)$ and $\Id_\hbar(R_\hbar(\A))$ are identified via the maps $\I\mapsto R_\hbar(\I), \I_\hbar\mapsto
\I_\hbar/(\hbar-1)\I_{\hbar}$. In particular, we have identifications $\Id(\Walg)\xrightarrow{\sim} \Id_\hbar(\Walg_\hbar),
\Id(\U)\xrightarrow{\sim} \Id_\hbar(\U_\hbar),\Id(\underline{\Walg})\xrightarrow{\sim} \Id(\underline{\Walg}_\hbar)$.

Also we have the following result whose proof is straightforward.

\begin{Lem}\label{Lem:4.11}
Let $\A,\B$ be filtered algebras and $\Phi:\A\rightarrow\B$ be a homomorphism preserving the filtrations.
Let $\Phi_\hbar:R_\hbar(\A)\rightarrow R_\hbar(\B)$ be the induced homomorphism of the Rees algebras.
Then $\Phi_\hbar^{-1}(R_\hbar(\I))=R_\hbar(\Phi^{-1}(\I))$ for any two-sided ideal $\I\subset\B$.
\end{Lem}

In \cite[Subsection 3.4]{Wquant},  we established a natural identification
$\Id_\hbar(\Walg_\hbar)\xrightarrow{\sim} \Id_\hbar(\Walg_\hbar^\wedge)$.  We have a similar identification
for $\underline{\Walg}$. Also we have a natural identification
$\Id_\hbar(\A_\hbar)\xrightarrow{\sim} \Id_\hbar(\W_{V,\hbar}^\wedge(\A_\hbar))$ provided $\A_\hbar$
is complete in the $\hbar$-adic topology. More precisely, the maps $\I_\hbar\mapsto \W_{V,\hbar}^\wedge\widehat{\otimes}_{\K[[\hbar]]}\I_\hbar, \J_\hbar\mapsto \J_\hbar\cap \A_\hbar$
are mutually inverse bijections between $\Id_\hbar(\A_\hbar)$ and $\Id_\hbar(\W_{V,\hbar}^\wedge(\A_\hbar))$.

\begin{Lem}\label{Lem:ideal_ident}
There is an identification
$\Id_\hbar(\underline{\U}^\wedge_\hbar)\xrightarrow{\sim} \Id_\hbar(\K[\widetilde{Y}]^\wedge_{\underline{G}y}[[\hbar]]^{\underline{G}})$.
\end{Lem}
\begin{proof}
Recall the identification $\K[\widetilde{Y}][\hbar]^{\underline{G}}\xrightarrow{\sim} \underline{\U}_\hbar\otimes_{\K[\hbar]}\D_\hbar(N^-)$. Taking the completions
at the point $(0,1,0)\in \underline{\g}^*\times N^-\times \n^{-*}=\underline{\g}^*\times T^*N^-$,
we get an isomorphism $\K[\widetilde{Y}]^\wedge_{\underline{G}y}[[\hbar]]^{\underline{G}}\xrightarrow{\sim}
\underline{\U}_\hbar^\wedge\widehat{\otimes}_{\K[[\hbar]]}\D_\hbar(N^-)^\wedge$, where the last factor
is the completion of $\D_\hbar(N^-)$ at $(1,0)$. This isomorphism is $\K^\times$-equivariant. So we get a natural map \begin{equation}\label{eq:ideals}\Id_\hbar(\K[\widetilde{Y}]^\wedge_{\underline{G}y}[[\hbar]]^{\underline{G}})\rightarrow \Id_\hbar(\underline{\U}^\wedge_\hbar)\end{equation} given by taking the intersection with $\underline{\U}^\wedge_\hbar$.

The algebra $\D_\hbar(N^-)^\wedge$ is isomorphic
to $\W_{\n\oplus \n^*, \hbar}^\wedge$.  
Hence the map $\underline{\J}'_\hbar\mapsto \D_\hbar(N^-)^\wedge\widehat{\otimes}_{\K[[\hbar]]}\underline{\J}'_\hbar$
is inverse to (\ref{eq:ideals}).
\end{proof}

Note that, thanks to Lemma \ref{Prop:2.4}, $\Lambda_\hbar$ acts trivially on $\Id_\hbar\left(\W_{V,\hbar}^\wedge(\K[Y]^\wedge_{\underline{G}y}[[\hbar]]^{\underline{G}})\right)$.

Summarizing, we have the following commutative diagram.
\begin{equation}\label{eq:diag1}
\begin{picture}(160,70)
\put(28,2){$\Id_\hbar(\K[Y]^\wedge_{\underline{G}y}[[\hbar]]^{\underline{G}})$}
\put(80,2){$\Id_\hbar(\underline{\Walg}^\wedge_\hbar)$}
\put(73,22){$\Id_\hbar(\W_{\underline{V},\hbar}^\wedge(\underline{\Walg}^\wedge_\hbar))$}
\put(127,22){$\Id_\hbar(\underline{\U}^\wedge_\hbar)$}
\put(115,42){$\Id_\hbar(\K[\widetilde{Y}]^\wedge_{\underline{G}y}[[\hbar]]^{\underline{G}})$}
\put(19,42){$\Id_\hbar(\W_{V,\hbar}^\wedge(\K[Y]^\wedge_{\underline{G}y}[[\hbar]]^{\underline{G}}))$}
\put(80,3){\vector(-1,0){20}}\put(70,4){\tiny $\cong$}
\put(34,7){\vector(0,1){33}}\put(35,20){\tiny $\cong$}
\put(90,7){\vector(0,1){13}}\put(91,13){\tiny $\cong$}
\put(105,23){\vector(1,0){20}}\put(115,24){\tiny $\cong$}
\put(131,27){\vector(0,1){13}}\put(132,32){\tiny $\cong$}
\put(63,43){\vector(1,0){50}}\put(80,44){\tiny $\cong$}
\put(4,15){$\Id_\hbar(\Walg_\hbar^\wedge)$}
\put(1,62){$\Id_\hbar(\W_{V,\hbar}^\wedge(\Walg^\wedge_\hbar))$}
\put(10,20){\vector(0,1){40}}\put(11,36){\tiny $\cong$}
\put(27,4){\vector(-1,1){9}}
\put(32,46){\vector(-1,1){14}}
\put(101,62){$\Id_\hbar(\U^\wedge_\hbar)$}
\put(30,63){\vector(1,0){70}}\put(60,64){\tiny $\cong$}
\put(124,46){\vector(-1,1){14}}
\put(5,2){$\Id(\Walg)$}
\put(128,2){$\Id(\underline{\Walg})$}
\put(128,62){$\Id(\U)$}
\put(10,7){\vector(0,1){6}}\put(11,10){\tiny $\cong$}
\put(131,7){\vector(0,1){12}}\put(132,13){\tiny $\cong$}
\put(131,47){\vector(0,1){12}}
\put(27,3){\vector(-1,0){10}}
\put(127,3){\vector(-1,0){30}}\put(110,4){\tiny $\cong$}
\put(117,63){\vector(1,0){11}}
\end{picture}
\end{equation}

Here almost all arrows are either natural identifications or are obtained from isomorphisms
of algebras. The maps $$\Id_\hbar(\K[Y]^\wedge_{\underline{G}y}[[\hbar]]^{\underline{G}})\rightarrow \Id_\hbar(\Walg^\wedge_\hbar),\quad \Id_\hbar(\W_{V,\hbar}^\wedge(\K[Y]^\wedge_{\underline{G}y}[[\hbar]]^{\underline{G}}))
\rightarrow \Id_\hbar(\W_{V,\hbar}^\wedge(\Walg^\wedge_\hbar)),$$
$$\Id_\hbar(\K[\widetilde{Y}]^\wedge_{\underline{G}y}[[\hbar]]^{\underline{G}})\rightarrow\Id_\hbar(\U^\wedge_\hbar),
\quad\Id_\hbar(\U^\wedge_\hbar)\rightarrow \Id(\U)(=\Id_\hbar(\U_\hbar))$$
are obtained by pull-backs with respect to the corresponding inclusions. Finally, the maps
$$\Id_\hbar(\K[Y]^\wedge_{\underline{G}y}[[\hbar]]^{\underline{G}})\rightarrow \Id(\Walg),\quad
\Id_\hbar(\K[\widetilde{Y}]^\wedge_{\underline{G}y}[[\hbar]]^{\underline{G}})\rightarrow\Id(\U)$$
are completely determined by the condition that the diagram is commutative.

{\it Step 4.}
Let $\underline{\I}_\hbar$ be the ideal in $\underline{\Walg}_\hbar$ corresponding to
$\underline{\I}$. Let $\underline{\I}^\wedge_\hbar$ denote the closure of $\underline{\I}_\hbar$ in $\underline{\Walg}^\wedge_\hbar$.
One has the equality
$$\underline{\I}'\cap \underline{\Walg}^\heartsuit=
(\underline{\I}^\wedge_{\hbar})_{\K^\times-l.f.}/(\hbar-1)(\underline{\I}^\wedge_{\hbar})_{\K^\times-l.f.},$$
where $\underline{\I}'$ is the closure of $\underline{\I}$ in $\underline{\Walg}'$.
Indeed, any ideal in $\underline{\Walg}^\heartsuit$ is generated by its intersection with $\underline{\Walg}$,
compare with the proof of   \cite[Lemma 5.3]{LOCat},
and both sides of the previous equality intersect $\underline{\Walg}$ in $\underline{\I}$.

Since the embedding $\Walg\hookrightarrow \underline{\Walg}'$ factors through $\Walg\hookrightarrow \underline{\Walg}^{\heartsuit}$,
we see that $\I$ is the image of $\underline{\I}$ under the maps of the bottommost row of the commutative diagram
(\ref{eq:diag1}).
From the construction, $\I^\dagger$ is the image of $\I$ under the maps of the leftmost
column and the topmost row of the diagram. So it remains to check that
\begin{itemize}
\item[(*)]
$\J$ is the image of
$\underline{\I}$ under the maps of the rightmost column.
\end{itemize}

We have the following commutative diagram, where all vertical arrows are natural embeddings.

\begin{equation}\label{eq:diagr2}
\begin{picture}(120,30)
\put(2,2){$\underline{\U}_\hbar$}
\put(2,20){$\underline{\U}^\wedge_\hbar$}
\put(20,2){$\D_\hbar(N^-)\otimes_{\K[\hbar]}\underline{\U}_\hbar$}
\put(27,20){$\W_{\n\oplus \n^*,\hbar}^\wedge(\underline{\U}_\hbar^\wedge)$}
\put(72,2){$\K[\widetilde{Y}][\hbar]^{\underline{G}}$}
\put(70,20){$\K[\widetilde{Y}]^\wedge_{\underline{G}y}[[\hbar]]^{\underline{G}}$}
\put(106,2){$\U_\hbar$}
\put(106,20){$\U^\wedge_\hbar$}
\put(4,7){\vector(0,1){11}}
\put(34,7){\vector(0,1){11}}
\put(75,7){\vector(0,1){11}}
\put(108,7){\vector(0,1){11}}
\put(7,3){\vector(1,0){12}}
\put(7,21){\vector(1,0){18}}
\put(52,3){\vector(1,0){18}}
\put(50,21){\vector(1,0){19}}
\put(105,3){\vector(-1,0){16}}
\put(105,21){\vector(-1,0){12}}
\put(61,4){\tiny $\cong$}
\put(58,22){\tiny $\cong$}
\end{picture}
\end{equation}


The diagram (\ref{eq:diagr2}) gives rise to the following commutative diagram of maps
between the sets of ideals.

\begin{equation}\label{eq:diagr11}
\begin{picture}(150,50)
\put(2,22){$\Id_\hbar(\underline{\U}_\hbar)$}
\put(2,40){$\Id_\hbar(\underline{\U}^\wedge_\hbar)$}
\put(30,22){$\Id_\hbar(\underline{\U}_\hbar\otimes_{\K[\hbar]}\D_\hbar(N^-))$}
\put(37,40){$\Id_\hbar(\W_{\n\oplus\n^*,\hbar}^\wedge(\underline{\U}_\hbar^\wedge))$}
\put(92,22){$\Id_\hbar(\K[\widetilde{Y}][\hbar])$}
\put(90,40){$\Id_\hbar(\K[\widetilde{Y}]^\wedge_{\underline{G}y}[[\hbar]])$}
\put(136,22){$\Id_\hbar(\U_\hbar)$}
\put(136,40){$\Id_\hbar(\U^\wedge_\hbar)$}
\put(6,38){\vector(0,-1){11}}
\put(46,38){\vector(0,-1){11}}
\put(97,38){\vector(0,-1){11}}
\put(144,38){\vector(0,-1){11}}
\put(18,23){\vector(1,0){12}}
\put(18,41){\vector(1,0){18}}
\put(70,23){\vector(1,0){20}}
\put(65,41){\vector(1,0){25}}
\put(117,23){\vector(1,0){18}}
\put(121,41){\vector(1,0){14}}
\put(25,42){\tiny $\cong$}
\put(23,24){\tiny $\cong$}
\put(75,42){\tiny $\cong$}
\put(81,24){\tiny $\cong$}
\put(3,2){$\Id(\underline{\U})$}
\put(32,2){$\Id(\underline{\U})\otimes \D(N^-)$}
\put(138,2){$\Id(\U)$}
\put(6,20){\vector(0,-1){14}}
\put(46,20){\vector(0,-1){14}}
\put(144,20){\vector(0,-1){14}}
\put(15,3){\vector(1,0){15}}
\put(62,3){\vector(1,0){76}}
\put(7,13){\tiny $\cong$}
\put(22,4){\tiny $\cong$}
\put(47,13){\tiny $\cong$}
\put(145,13){\tiny $\cong$}
\end{picture}
\end{equation}

The map in the rightmost column of diagram (\ref{eq:diag1}) is the composition of the identification
$\Id(\underline{\Walg})\xrightarrow{\sim} \Id_\hbar(\underline{\U}^\wedge_\hbar)$
and the map $\Id_\hbar(\underline{\U}_\hbar^\wedge)\rightarrow \Id_\hbar(\U_\hbar)$ from  diagram
(\ref{eq:diagr11}). Now the claim  (*) follows from
Lemma \ref{Lem:3.31}.
\end{proof}

\subsection{Parabolic induction and representation schemes}\label{SUBSECTION_rep_schemes}
In this subsection we study the morphism of representation schemes
induced by the parabolic induction functor.

Let us recall some generalities on representation schemes.

Let $\A$ be a finitely generated associative algebra with generators $x_1,\ldots,x_n$.
Consider the ideal $\J$ of relations for $\A$ in the free algebra $\K\langle x_1,\ldots,x_n\rangle$  so that $\A\cong \K\langle x_1,\ldots,x_n\rangle/\J$. Fix some positive
integer $d$ and consider the subscheme $X\subset \Mat_d(\K)^n=\{(X_1,\ldots,X_n), X_i\in \Mat_d(\K)\}$ defined by the
equations $f(X_1,\ldots,X_n)=0$ for $f\in \J$. By definition the {\it representation
scheme} $\Rep(\A,d)$ is the categorical quotient $X\quo \GL_d$.
The points of $\Rep(\A,d)$ are in bijection with isomorphism classes
of semisimple $\A$-modules of dimension $d$. A homomorphism
$\A_1\rightarrow \A_2$ induces a morphism $\Rep(\A_2,d)\rightarrow \Rep(\A_1,d)$.

For elements $x_1,\ldots,x_{2n}$ in an associative algebra we put
$$s_{2n}(x_1,\ldots,x_{2n}):=\sum_{\sigma\in \mathfrak{S}_{2n}}\operatorname{sgn}(\sigma)x_{\sigma(1)}x_{\sigma(2)}\ldots x_{\sigma(2n)}.$$
An algebra is said to satisfy the identity $s_{2n}$ if $s_{2n}(x_1,\ldots,x_{2n})=0$ for all elements
$x_1,\ldots,x_{2n}$ of this algebra.
According to the Amitsur-Levitzki Theorem (see e.g. \cite[Theorem 13.3.3]{MR}), the algebra of
$\Mat_{d}(\K)$ satisfies $s_{2n}$ provided $d\leqslant n$.

 For an algebra $\A$ let $\A^{(n)}$ denote the quotient
of $\A$ by the two-sided ideal generated by all elements $s_{2n}(a_1,\ldots,a_{2n}), a_i\in \A$.
Clearly, the schemes $\Rep(\A,d)$ and $\Rep(\A^{(d)},d)$ are canonically isomorphic.

\begin{Prop}\label{Lem:8.1}
\begin{enumerate}
\item If $\underline{e}$ is rigid, then the algebra $U([\underline{\g},\underline{\g}],\underline{e})^{(d)}$ is finite dimensional for any $d$.
\item We have an isomorphism $\underline{\Walg}^{(d)}\xrightarrow{\sim} U(\z(\underline{\g}))\otimes U([\underline{\g},\underline{\g}],\underline{e})^{(d)}$.
\item
The inclusion $\underline{\Walg}\hookrightarrow \underline{\Walg}'$ induces an
isomorphism $\underline{\Walg}^{(d)}\rightarrow \underline{\Walg}'^{(d)}$.
\end{enumerate}
\end{Prop}
\begin{proof}
We  derive  assertion (1) from  Corollary \ref{Cor:3.11} proved below
in the Appendix.  To do this we need to verify
that the Poisson variety $\underline{S}_0:=\underline{S}\cap [\underline{\g},\underline{\g}]=\Spec(\gr U([\underline{\g},\underline{\g}],\underline{e}))$ has only one 0-dimensional symplectic leaf.
Symplectic leaves in $\underline{S}_0$
are exactly the irreducible components of the
intersections of $\underline{S}_0$ with coadjoint orbits of $\underline{G}$, see, for example, \cite[3.1]{GG}.
So we need to check that the only orbit in $[\underline{\g},\underline{\g}]$ of dimension $\dim\underline{\Orb}$
that intersects $\underline{S}_0$ is $\underline{\Orb}$ itself.

Assume the converse. Thanks to the contraction property of the Kazhdan action, an orbit $\underline{\Orb}'$
intersects $\underline{S}'$ if and only if $\underline{\Orb}\subset \overline{\K^\times \underline{\Orb}'}$.
 Recall that by a sheet in $[\underline{\g},\underline{\g}]^*$ one means an irreducible component of the locally
closed subvariety $\{\alpha\in [\underline{\g},\underline{\g}]^*| \dim \underline{G}\alpha=k\}$ for some fixed
$k$.  A sheet in a semisimple algebra  containing a rigid nilpotent orbit consists only of this orbit,
see, for example, \cite[$\S$5.5]{McG}. So if $\underline{\Orb}\subset \overline{\K^\times \underline{\Orb}'}$
and $\underline{\Orb}$ is rigid,
then $\dim \underline{\Orb}'>\dim \underline{\Orb}$.

It follows that the only zero-dimensional symplectic leaf in $\underline{S}_0$ is $\underline{e}$.
So Corollary \ref{Cor:3.11} does apply in the present situation.

Assertion (2) follows from the decomposition $\underline{\Walg}=U(\z(\underline{\g}))\otimes
U([\underline{\g},\underline{\g}],\underline{e})$.

Let us proceed to assertion (3). Let $\I$, resp. $\I'$, denote the ideal in $\underline{\Walg}$
(resp., $\underline{\Walg}'$) generated by $s_{2d}(x_1,\ldots, x_{2d})$ with $x_i\in \underline{\Walg}$,
(resp., $x_i\in \underline{\Walg}'$). Clearly, $\I'$ is the closure of $\I$
in $\underline{\Walg}'$. We need to check that $\I'\cap \underline{\Walg}=\I$ and
$\underline{\Walg}+\I'=\underline{\Walg}'$. Let $\underline{\Walg}=\bigoplus_{\alpha\in \Z}\underline{\Walg}_\alpha$
denote the eigenspace decomposition with respect to $\ad h_0$. We claim that $\Walg_\alpha\subset \I$
for all  $\alpha$ less than some number $\alpha_0$ depending on $d$.

Let $Y$ denote the union of all sheets containing the orbit $\underline{\Orb}$.
By Katsylo's results, \cite{Katsylo}, the group $Q^\circ$ acts trivially on
$Y\cap S$. By Theorem \ref{Thm:3}, the maximal spectrum $\operatorname{Specm} \gr\Walg^{(d)}$ is contained in $Y\cap S$ so the $\K^\times$-action
on $\operatorname{Specm}(\gr\Walg^{(d)})$ (via $\gamma_0$) is trivial. It follows that
$\ad h_0$ has finitely many eigenvalues on $\gr\underline{\Walg}^{(d)}$ and hence on $\underline{\Walg}^{(d)}$.

Therefore $\I'=\I+ \prod_{\alpha\leqslant \alpha_0} \underline{\Walg}_\alpha$.
So $\I=\I'\cap\underline{\Walg}$ and $\I'+\underline{\Walg}=\underline{\Walg}'$.
\end{proof}

So we have a homomorphism $\Walg^{(d)}\rightarrow \underline{\Walg}^{(d)}$. It gives
rise to a morphism $\Rep(\underline{\Walg},d)\rightarrow \Rep(\Walg,d)$. Below we assume
that the element $\underline{e}$ is rigid.

\begin{Thm}\label{Thm:2.7}
The morphism $\Rep(\underline{\Walg},d)\rightarrow \Rep(\Walg,d)$ is finite.
\end{Thm}
\begin{proof}
Recall that the centers $\underline{\Centr},\Centr$ of $\underline{\U},\U$ are identified with
the centers of $\underline{\Walg},\Walg$.
 By Corollary \ref{Cor:2.9}, the following diagram is commutative

\begin{picture}(60,30)
\put(2,2){$\underline{\Centr}$}\put(2,22){$\Centr$}
\put(40,2){$\underline{\Walg}^{(d)}$}
\put(41,22){$\Walg^{(d)}$}
\put(4,19){\vector(0,-1){13}}
\put(42,19){\vector(0,-1){13}}
\put(7,3){\vector(1,0){32}}
\put(7,23){\vector(1,0){32}}
\end{picture}

By Proposition \ref{Lem:8.1}, $\underline{\Walg}^{(d)}=U(\z(\underline{\g}))\otimes\A$, where
$\A$ is a finite dimensional algebra. Let us check that the morphism $\Rep(\underline{\Walg},d)=\Rep(\underline{\Walg}^{(d)},d)\rightarrow
\Rep(\underline{\Centr},d)$ is finite. For this we will describe the varieties in interest.
Clearly, a morphism of schemes is finite if the induced morphism  of the underlying varieties
is so.

As a variety,  $\Rep(\underline{\Centr},d)$ is just
$\Spec(\underline{\Centr})^d/S_d$.  Indeed, this variety parameterizes  semisimple
$\underline{\Centr}$-modules of dimension $d$. Such a module is determined  up to an isomorphism
by an unordered $d$-tuple of characters of $\underline{\Centr}$.

Proceed to $\Rep(\underline{\Walg}^{(d)},d)$. As we have seen above, $\underline{\Walg}^{(d)}=U(\z(\underline{\g}))\otimes \A$, where $\A$ is a finite dimensional algebra over $\Centr([\underline{\g},\underline{\g}])$.
An irreducible $\underline{\Walg}^{(d)}$-module has the form $\K_\chi\otimes V$, where $V$
is an irreducible $\A$-module and $\K_\chi$ is the one-dimensional $U(\z(\underline{\g}))$-module
corresponding to a character $\chi$. So $\Rep(\underline{\Walg}^{(d)})$ is a disjoint union of
irreducible components $\Rep(\underline{\Walg}^{(d)})_{v}$, where $v$ is an unordered collection
of irreducible $\A$-modules with sum of dimensions equal $d$. Let $|v|$ denote the total number of representations
in $v$. The component is naturally identified with the quotient of $(\z(\underline{\g})^*)^{|v|}$ by
an action of  an appropriate
product of symmetric groups. The natural morphism $\Rep(U(\z(\underline{\g}))\otimes \A,d)\rightarrow
\Rep(\underline{\Centr},d)=\Rep(U(\z(\underline{\g}))\otimes \Centr([\underline{\g},\underline{\g}]),d)$ can be read from the characters
appearing in the restrictions of irreducible $\A$-modules to $\Centr([\underline{\g},\underline{\g}])$.
It is easy to see that this morphism is finite.

From the description of $\Rep(\underline{\Centr},d)$ (the variety $\Rep(\Centr,d)$
can be described similarly) it is clear that the morphism $\Rep(\underline{\Centr},d)\rightarrow
\Rep(\Centr,d)$ is also finite. So the composition $\Rep(\underline{\Walg}^{(d)},d)\rightarrow \Rep(\Centr,d)$
is finite.
\end{proof}

\subsection{Adjoint functors}\label{SUBSECTION_adjoint}
The goal of this subsection is to construct the left and right adjoint functors
of the parabolic induction functor $\rho$.

Let $M$ be a finite dimensional $\Walg$-module. Pick an ideal $\J\subset\Centr$ of finite codimension annihilating
$M$. So $M$ is a module over $\Walg_\J:=\Walg/\Walg \J$.

\begin{Prop}\label{Lem:3.11}
There is a minimal two-sided ideal $\I_0$ of finite codimension in $\Walg_\J$.
\end{Prop}
\begin{proof}
It is enough to check that the $\Walg$-bimodule $\Walg_\J$ has finite length.
This will follow for an arbitrary $\J$ if we prove it  for a {\it maximal} one.
The claim for maximal $\J$ is  \cite[Theorem 4.2.2,(i)]{Ginzburg}.
\end{proof}

In particular, the sequence $\ker(\Walg_\J\rightarrow\Walg_\J^{(d)})$ stabilizes
for any given $\J$. Similarly, for $\underline{\Walg}_\J:=\underline{\Walg}/\underline{\Walg}\J$
(where we consider $\J$ as a subspace in $\underline{\Centr}$ via the embedding $\Centr\hookrightarrow
\underline{\Centr}$) the sequence
$\underline{\Walg}_\J^{(d)}$ stabilizes. So there is $N\in \N$ such that
$\Walg_\J^{(d)}=\Walg_\J^{(N)}$ and $\underline{\Walg}_\J^{(d)}=\underline{\Walg}_\J^{(N)}$
for all $d>N$.

Let $\Walg_\J$-$\Mod^{fin},\underline{\Walg}_\J$-$\Mod^{fin}$ denote the categories
of finite dimensional modules for the corresponding algebras. Let $\rho_\J:\underline{\Walg}_\J$-$\Mod^{fin}\rightarrow
\Walg_\J$-$\Mod^{fin}$ be the corresponding pullback functor.

Set $\kappa_\J^l(\bullet):=\underline{\Walg}_\J^{(N)}\otimes_{\Walg_\J^{(N)}}\bullet, \kappa_\J^r(\bullet):=\Hom_{\Walg_\J^{(N)}}(\underline{\Walg}_\J^{(N)},\bullet)$. These are the
left and right adjoint functors of $\rho_\J$.

Now let $\J_1\subset\J_2$ be two ideals of $\Centr$ of finite codimension. Then we have
the natural embedding $\iota_{12}:\Walg_{\J_2}$-$\Mod^{fin}\rightarrow \Walg_{\J_1}$-$\Mod^{fin}$
and the similar embedding $\underline{\iota}_{12}$ for $\underline{\Walg}$.

It is clear that $\kappa_{\J_1}^l\circ \iota_{12}=\underline{\iota}_{12}\circ \kappa_{\J_2}^l$
and $\kappa_{\J_1}^{r}\circ \iota_{12}=\underline{\iota}_{12}\circ \kappa_{\J_2}^{r}$. Since the
category $\Walg$-$\Mod^{fin}$ is the direct limit of its full subcategories $\Walg_\J$-$\Mod^{fin}$
(and  similarly for $\underline{\Walg}$) we get the well-defined direct limit
functors $\kappa^l,\kappa^r: \Walg$-$\Mod^{fin}\rightarrow \underline{\Walg}$-$\Mod^{fin}$
that are left and right adjoint to $\rho$.

\section{Appendices}\label{SECTION_appendix}
\subsection{Fedosov quantization vs twisted differential operators}
\label{SUBSECTION_Fedosov_vs_TDO}
Let $X_0$ be a smooth affine variety. Consider the cotangent bundle $X:=T^*X_0$.
Then $X$ has a canonical symplectic form $\omega$.
Our first goal in this section is to relate the Fedosov quantization corresponding to the zero
curvature form to a certain algebra of twisted differential operators on $X_0$.

For computations we will need a precise description of the Poisson structure on $\K[X]$.
Recall that the symplectic form $\omega$ on $X$ equals $-d\lambda$, where
$\lambda$ is the canonical 1-form on $X$ defined as follows. A typical point of $X$
is $x=(x_0,\alpha)$, where $x_0\in X_0,\alpha\in T_x^*X_0$. Equip $X$ with the $\K^\times$-action
by setting $t.(x_0,\alpha)=(x_0,t^{-2}\alpha)$.
Let $\pi$ denote the natural projection
$X\twoheadrightarrow X_0$.
Then for $v\in T_xX$ we set $\langle\lambda_x,v\rangle:=\langle\alpha,d_x\pi(v)\rangle$.

Identify $\K[X_0]$ with a subalgebra in $\K[X]$ via $\pi^*$ and also embed $\Vect(X_0)$
into $\K[X]$ as the space of functions of degree 2 with respect to the $\K^\times$-action. Then
the Poisson bracket on $\K[X]$ is determined by the following equalities:
\begin{equation}\label{eq:3}
\begin{split}
&\{f_1,f_2\}=0,\\
&\{f,v\}=\partial_v f=-L_vf,\\
&\{v_1,v_2\}=[v_1,v_2],\\
&f,f_1,f_2\in \K[X_0], v,v_1,v_2\in \Vect(X),
\end{split}
\end{equation}
where $L_v$ denotes the Lie derivative w.r.t $v$.
With this sign convention the Cartan magic formula has the form
$$L_v\lambda=-d\iota_v\lambda-\iota_vd\lambda, v\in \Vect(X_0), \lambda\in \Omega^{\bullet}(X_0).$$

Consider the algebra
$\D(X_0)$ of linear differential operators on $X_0$ and equip it with the
filtration $\D^{\leqslant i}(X_0)$, where $\D^{\leqslant i}(X_0)$ consists
of differential operators of order  $\leqslant\frac{i}{2}$.
Form the Rees algebra $D_\hbar(X_0):=\bigoplus_{i=0}^\infty \D^{\leqslant i}(X_0)\hbar^i$.
This algebra has a natural
$\K^\times$-action.
Note also that $\D_\hbar(X_0)$ can be considered as the algebra of global section of
the {\it sheaf} of homogeneous differential operators on $X_0$.

Consider the bundle $\Omega^{top}(=\Omega^{top}_{X_0})$ of top differential forms on $X_0$, and let
$\widetilde{X}_0$ denote the total space of this bundle. This is a smooth algebraic variety
acted on by the torus $\K^\times$. The algebra $\K[\widetilde{X}_0]^{\K^\times}$ is naturally identified
with $\K[X_0]$, while the space of functions of degree 1 is nothing else but $\Gamma(X_0,\Omega^{top})$.
The $\K^\times$-action gives rise to the Euler vector field $\eu$ on $\widetilde{X}_0$. Consider the
algebra $\D^{\frac{1}{2}\Omega^{top}}(X_0)$ of twisted differential operators on $\frac{1}{2}\Omega^{top}$, i.e. the algebra $\D(\widetilde{X}_0)^{\K^\times}/(\eu-\frac{1}{2})$. This algebra has a  filtration
$\F_i\D^{\frac{1}{2}\Omega^{top}}(X_0)$ similar to the one above and we can form the Rees algebra
$$\D_\hbar^{\frac{1}{2}\Omega^{top}}(X_0):=\bigoplus_{i=0}^\infty \F_i\D^{\frac{1}{2}\Omega^{top}}(X_0)\hbar^{i}=
\D_\hbar(\widetilde{X}_0)^{\K^\times}/(\eu-\frac{1}{2}\hbar^2).$$

Let $\varrho$ denote the projection $\D(\widetilde{X}_0)^{\K^\times}\twoheadrightarrow \D^{\frac{1}{2}\Omega^{top}}(X_0)$.
We have the natural embedding $\iota:\K[X_0]\hookrightarrow \D(\widetilde{X}_0)$. For
$v\in \Vect(X_0)$ define a $\K^\times$-invariant vector field $\iota(v)$ on $\widetilde{X}_0$ by
$$L_{\iota(v)}f=L_vf, L_{\iota(v)}\sigma=L_{v}\sigma, f\in \K[X_0], \sigma\in \Gamma(X_0,\Omega^{top}).$$
For brevity, put $\theta=\varrho\circ\iota$.

\begin{Prop}\label{Prop:3.1}
Equip $\K[X][\hbar]$ with a homogeneous
Fedosov star-product corresponding to the zero curvature form. Then there is a unique homomorphism
$\K[X][\hbar]\rightarrow \D_\hbar^{\frac{1}{2}\Omega^{top}}(X_0)$ mapping
$f\in \K[X_0]$ to $\theta(f)$ and $v\in \Vect(X_0)$ to $\hbar^2\theta(v)$. This homomorphism
is an isomorphism.
\end{Prop}
\begin{proof}
Let $f*g=\sum_{i=0}^\infty D_i(f,g)\hbar^{2i}$ be the star-product on $\K[X][\hbar]$. Then it is well known, see, for example, \cite[Lemma 3.3]{BW}, that
$D_i(f,g)=(-1)^i D_i(g,f)$ for all $f,g\in \K[X]$. In particular, $D_1(f,g)=\frac{1}{2}\{f,g\}$ and $D_2$ is symmetric.
So we have
\begin{equation}\label{eq:3.2}
\begin{split}
&f_1*f_2=f_1f_2.\\
&f*v=fv-\frac{1}{2}L_vf\hbar^2.\\
&v*f=fv+\frac{1}{2}L_vf\hbar^2.\\
&v_1*v_2-v_2*v_1=[v_1,v_2]\hbar^2,\\
&f,f_1,f_2\in \K[X_0], v,v_1,v_2\in \Vect(X_0).
\end{split}
\end{equation}

Consider the $\K[\hbar]$-algebra $\widehat{\D}_\hbar$ generated by $\K[X_0]$ and $\Vect(X_0)$ subject to the relations
(\ref{eq:3.2}). Equip $\widehat{\D}_\hbar$ with the filtration ``by the order of a differential
operator'': $\F_k\widehat{\D}_\hbar=\K[X_0][\hbar]\Vect(X_0)^k$.
We have a natural epimorphism $\widehat{\D}_\hbar\rightarrow \K[X][\hbar]$. Passing to
the associated graded algebras we get a homogeneous homomorphism $\gr\widehat{\D}_\hbar\rightarrow \K[X][\hbar]$
(where the last space is considered as an algebra with respect to the commutative product).
This homomorphism has an inverse, a natural epimorphism
$\K[X][\hbar]\twoheadrightarrow \gr \widehat{\D}_\hbar$, because the algebra $\K[X][\hbar]$ is naturally
identified with $S_{\K[X_0][\hbar]}(\Vect(X_0)\otimes \K[\hbar])$. It follows that the natural epimorphism $\widehat{\D}_\hbar\rightarrow \K[X][\hbar]$ (here $\K[X][\hbar]$ is a quantum algebra) is an isomorphism.

We are going to check that in $\D^{\frac{1}{2}\Omega^{top}}(X_0)$
the relations analogous to (\ref{eq:3.2}) hold. Clearly, $\iota(f_1)\iota(f_2)=\iota(f_1f_2)$. Let us compute now $\iota(f)\iota(v)$. We have
\begin{equation}\label{eq:3.3.1}
\iota(f)\circ\iota(v)g=fL_vg, g\in \K[X_0],
\end{equation}
\begin{equation}\label{eq:3.3.2}
\iota(f)\circ \iota(v)\sigma=f L_v\sigma=L_{fv}\sigma+df\wedge \iota_v\sigma=L_{vf}\sigma-(L_vf)\sigma.
\end{equation}
So we see that $\iota(f)\iota(v)=\iota(fv)-\iota(L_vf)\eu$. Therefore in $\D^{\frac{1}{2}\Omega^{top}}(X_0)$
we have $$\theta(f)\theta(v)=\theta(fv)-\frac{1}{2}\theta(L_vf).$$
Similarly, we have
$$\theta(v)\theta(f)=\theta(fv)+\frac{1}{2}\theta(L_vf).$$
Finally, from the definition of $\theta$ we have $\theta([v_1,v_2])=[\theta(v_1),\theta(v_2)]$.

Since $\K[X][\hbar]$ is identified with $\widehat{\D}_\hbar$ as above,  we have a unique homomorphism $\K[X][\hbar]\rightarrow \D_{\hbar}^{\frac{1}{2}\Omega^{top}}(X_0)$ sending $f\in \K[X_0], v\in \Vect(X_0)$ to $\theta(f),\hbar^2\theta(v)$.
Define a filtration on $\D_\hbar^{\frac{1}{2}\Omega^{top}}(X_0)$ by
$$\operatorname{G}_i\D_\hbar^{\frac{1}{2}\Omega^{top}}(X_0)=\theta(\K[X_0])[\hbar](\hbar^2\theta(\Vect(X_0)))^i$$
(i.e., again by the order of a differential operator). Then the associated graded algebra
is again naturally identified with the commutative algebra $\K[X][\hbar]$. Moreover, under this identification,
the associated graded of the homomorphism under consideration is the identity.
So we see that the homomorphism is actually an isomorphism.
\end{proof}

Suppose that we have an action of $G$ on $X_0$.
 Then the map
$\xi\mapsto \widehat{H}_\xi:=\xi_{X_0}\hbar^2\in \D^{\leqslant 2}(X_0)\hbar^2$, where $\xi_{X_0}$
is the velocity vector field associated to $\xi$, is a quantum comoment map.

The $G$-action on $X_0$ lifts naturally to $\widetilde{X}_0$.
The quantum moment map for the corresponding $G$-action on $\D_\hbar(\widetilde{X}_0)$ is
$\xi\mapsto \iota(\xi_{X_0})\hbar^2$. The $G$-action descends to $\D^{\frac{1}{2}\Omega^{top}}_\hbar(X_0)$
with a quantum comoment map $\xi\mapsto  \theta(\xi_{X_0})\hbar^2$.

In particular, if $\Omega^{top}$ is $G$-equivariantly trivial, then there is a $G\times\K^\times$-equivariant
Hamiltonian isomorphism $\K[X][\hbar]\rightarrow \D_\hbar(X_0)$.

\subsection{A general result about filtered algebras}\label{SUBSECTION_Ass}
 Let $\A$ be an associative algebra
with unit equipped with an increasing exhaustive filtration $\F_i\A, i\geqslant 0$.
Let $A:=\sum_i A_i$, where $A_i:=\F_i\A/\F_{i-1}\A$, be the corresponding associated graded algebra. We assume that $A$ is finitely generated over
$\K$. Suppose that there is $d>0$ with $[\F_i\A,\F_j\A]\subset \F_{i+j-d}\A$ for all $i,j$.
Then $A$ has a canonical Poisson bracket induced from $\A$ with $\{A_i,A_j\}\subset A_{i+j-d}$.
Consider the Poisson ideal $I:=A\{A,A\}$. Then $I$
is the minimal ideal of $A$ such that the algebra $A/I$ is Poisson commutative.

Fix a positive integer $n$. Recall the quotient $\A^{(n)}$ of $\A$ defined
in Subsection \ref{SUBSECTION_rep_schemes}. Let $\I_n$ be the kernel of the natural epimorphism
$\A\twoheadrightarrow\A^{(n)}$.

\begin{Thm}\label{Thm:3}
In the above notation,
$\sqrt{\gr \I_n}\supset I$.
\end{Thm}
\begin{proof}
The filtration on $\A$ induces a filtration $\F_i\A^{(n)}$ on $\A^{(n)}$. Set $B:=\gr\A^{(n)}$. We need
to show that all symplectic leaves of the Poisson subscheme $\Spec B\subset \Spec A$ are 0-dimensional.
Assume the converse, pick a point  $x\in \Spec B$ whose symplectic leaf has positive dimension.
Without loss of generality, we may assume that $x$ lies in the smooth locus of the reduced scheme associated with $B$.

  The algebra  $\A^{(n)}$ satisfies the identities $s_{2n}$.
Consider the Rees algebra $R_\hbar(\A^{(n)})=\bigoplus_{i\geqslant 0}\hbar^i \F_i\A^{(n)}$. Then
$R_\hbar(\A^{(n)})/(\hbar-a)\cong\A^{(n)}$ for $a\in \K, a\neq 0,$ and $R_\hbar(\A^{(n)})/(\hbar)\cong B$. In particular, the
algebras $R_\hbar(\A^{(n)})/(\hbar-a)$ satisfy $s_{2n}$ for all $a\in \K$. So $R_\hbar(\A^{(n)})$
also satisfies $s_{2n}$. Let us note that $R_\hbar(\A^{(n)})$ is flat over $\K[\hbar]$.

Let $\m_x$ be the maximal ideal of $x$ in $B$ and let $\widetilde{\m}_x$ be the inverse image
of $\m_x$ in $R_\hbar(\A^{(n)})$.  Consider the completion $R^\wedge_\hbar(\A^{(n)})$ of $R_\hbar(\A^{(n)})$ w.r.t
$\widetilde{\m}_x$, i.e., $R^\wedge_\hbar(\A^{(n)}):=\varprojlim R_\hbar(\A^{(n)})/\widetilde{\m}_x^n$.

\begin{Lem}
$R^\wedge_\hbar(\A^{(n)})$ is $\K[[\hbar]]$-flat.
\end{Lem}
\begin{proof}
Consider the $\hbar$-adic completions $\widetilde{\m}'_x,R'_\hbar(\A^{(n)})$ of $\widetilde{\m}_x,R_\hbar(\A^{(n)})$.
The assumptions of   \cite[Lemma 2.4.2]{HC} hold for $\widetilde{\m}'^2_x\subset R'_\hbar(\A^{(n)})$. So the {\it blow-up} algebra $\operatorname{Bl}_{\widetilde{\m}'_x}(R'_\hbar(\A^{(n)})):=\bigoplus_{i\geqslant 0}\widetilde{\m}'^{2i}_x$ is Noetherian.
Now we can repeat the argument used in the proof  of  \cite[Proposition 2.4.1]{HC}.
\end{proof}

It follows from the construction that $R^\wedge_\hbar(\A^{(n)})$ satisfies $s_{2n}$. Also note that
$R^\wedge_\hbar(\A^{(n)})/(\hbar)=B_x^\wedge:=\varprojlim B/\m_x^k$.
Let us check that there are elements $\widetilde{a},\widetilde{b}\in R^\wedge_{\hbar}(\A^{(n)})$ with $[\widetilde{a},\widetilde{b}]=\hbar^d$.

Recall that the symplectic leaf of $\Spec(B)$ passing through $x$ has positive dimension.
Therefore, \cite[Proposition 3.3]{Kaledin}, implies that $B_x^\wedge$ can be decomposed
into the completed tensor product of the algebra $\K[[a,b]]$ with $\{a,b\}=1$ and of some
other Poisson algebra $B'$. Lift $a,b$ to some elements $\widetilde{a},\widetilde{b}'$ of
$R^\wedge_\hbar(\A^{(n)})$. Then $[\widetilde{a},\widetilde{b}']\in \hbar^d+\hbar^{d+1}R^\wedge_{\hbar}(\A^{(n)})$.
The map $\{a,\bullet\}:B^\wedge_x\rightarrow B^\wedge_x$ is surjective.
This observation easily implies that we can find an element $y\in \hbar R^\wedge_\hbar(\A^{(n)})$ such that
$\widetilde{a}$ and $\widetilde{b}=\widetilde{b}'+\hbar y$ will satisfy $[\widetilde{a},\widetilde{b}]=\hbar^d$.

Now consider the Weyl algebra $\W_{\hbar}$ with generators $u,v$ and the relation
$[u,v]=\hbar^d$ (recall that previously we had $d=2$). We have an obvious
homomorphism $\W_\hbar\rightarrow R^\wedge_\hbar(\A^{(n)})$. This homomorphism is injective because
$R^\wedge_\hbar(\A^{(n)})$ is $\K[[\hbar]]$-flat. So $\W_{\hbar}$ satisfies $s_{2n}$. Hence the Weyl algebra
$\W=\W_{\hbar}/(\hbar-1)$ also satisfies $s_{2n}$.  This is
impossible, for example, because the center of $\W$ is trivial, (see, for instance, \cite[Proposition 13.6.11]{MR}).
%
\end{proof}

\begin{Cor}\label{Cor:3.11}
Suppose that the scheme $\Spec(A)$ has only one 0-dimensional symplectic leaf.
Then $\I_n$ has finite codimension in $\A$ for any $n$. In particular, $\A$ has finitely many
irreducible representations of any given dimension.
\end{Cor}
\begin{proof}
The ideal $I$ is just the maximal ideal of a point in $A$. So $\gr\I_n$ is of finite codimension in
$A$. It follows that $\I_n$ is of finite codimension in $\A$.
\end{proof}

%
%
%

\end{document}